\documentclass{elsarticle}

\usepackage{lineno,hyperref}

\modulolinenumbers[5]

\journal{Journal of Differential Equations}

\usepackage{amsopn}

\usepackage{lipsum}
\usepackage{amsfonts,amsthm}
\usepackage{tikz,pgfplots,subcaption}
\pgfplotsset{compat=newest}

\usepgfplotslibrary{fillbetween}
\usepackage{epstopdf}
\usepackage{algorithmic}
\usepackage{mathtools}
\usepackage[utf8]{inputenc}
\usepackage[page]{appendix}
\usepackage{float}
\ifpdf
  \DeclareGraphicsExtensions{.eps,.pdf,.png,.jpg}
\else
  \DeclareGraphicsExtensions{.eps}
\fi

  \newcommand{\matlab}{MATLAB\textsuperscript{\textregistered}}
\newtheorem{Prop}{Proposition}[section]
\newtheorem{Lem}[Prop]{Lemma}
\newtheorem{Thm}[Prop]{Theorem}

\newtheorem{Def}[Prop]{Definition}
\newtheorem{Rem}[Prop]{Remark}

\newtheorem{Ass}[Prop]{Assumption}

\newcommand{\cB}{\mathcal{B}}

\newcommand{\cA}{\mathcal{A}}

\newcommand{\cL}{\mathcal{L}}

\newcommand{\cD}{\mathcal{D}}

\newcommand{\cF}{\mathcal{F}}

\newcommand{\N}{{\mathbb{N}}}
\newcommand{\Z}{{\mathbb{Z}}}
\newcommand{\R}{{\mathbb{R}}}
\newcommand{\C}{{\mathbb{C}}}
\newcommand{\D}{{\textcolor{blue}{\mathfrak{D}}}}

\newcommand{\T}{{\mathbb{T}}}

\newcommand{\fa}{\mathfrak{a}}

\newcommand{\sg}{(\T_t)_{t\geq0}}

\newcommand{\st}{\ |\ }

\newcommand{\yref}{y_{\rm ref}}
\newcommand{\yrefk}{y_{{\rm ref},k}}
\newcommand{\dyref}{\dot{y}_{\rm ref}}

\newcommand{\ee}{\mathrm{e}}

\DeclareMathOperator{\im}{im}

\DeclareMathOperator{\esssup}{ess\, sup}

\newcommand{\setdef}[2]{\left\{\ #1\ \left|\ \vphantom{#1} #2\ \right.\right\}}
\newcommand{\ddt}{\tfrac{\text{\normalfont d}}{\text{\normalfont d}t}}

\newcommand{\ds}[1]{{\rm \, d} #1 \,}

\newcommand{\scpr}[2]{\left\langle #1,#2\right\rangle}

\DeclareMathOperator{\loc}{loc}
\newcommand*{\QED}{\hfill\ensuremath{\square}}

\begin{document}

\begin{frontmatter}

\title{Funnel control for the monodomain equations with the FitzHugh-Nagumo model}
\tnotetext[mytitlenote]{Thomas Berger acknowledges support by the German Research Foundation (Deutsche Forschungsgemeinschaft) via the grant BE 6263/1-1.}

\author{Thomas Berger}
\address{Institut f\"ur Mathematik\\ Universit\"at Paderborn\\ Warburger Str.~100\\ 33098~Paderborn\\ Germany}
\ead{thomas.berger@math.upb.de}

\author{Tobias Breiten}
 \address{Institute of Mathematics\\ Technical University of Berlin\\ Stra\ss e des 17. Juni 136 \\ 10623 Berlin\\ Germany}
 \ead{tobias.breiten@tu-berlin.de}

\author{Marc Puche}
\address{Fachbereich Mathematik\\ Universit\"at Hamburg, Bundesstra{\ss}e~55\\ 20146~Hamburg\\ Germany}
\ead{marc.puche@uni-hamburg.de}

\author{Timo Reis}
\address{Fachbereich Mathematik\\ Universit\"at Hamburg, Bundesstra{\ss}e~55\\ 20146~Hamburg\\ Germany}
\ead{timo.reis@uni-hamburg.de}

\begin{abstract}
We consider a nonlinear reaction diffusion system of parabolic type known as the monodomain equations, which model the interaction of the electric current in a cell. Together with the FitzHugh-Nagumo model for the nonlinearity they represent defibrillation processes of the human heart. We study a fairly general type with co-located inputs and outputs describing both boundary and distributed control and observation. The control objective is output trajectory tracking with prescribed performance. To achieve this we employ the funnel controller, which is model-free and of low complexity. The controller introduces a nonlinear and time-varying term in the closed-loop system, for which we prove existence and uniqueness of solutions. Additionally, exploiting the parabolic nature of the problem, we obtain H\"older continuity of the state, inputs and outputs. We illustrate our results by a simulation of a standard test example for the termination of reentry waves.
\end{abstract}


\end{frontmatter}

\linenumbers

\section{Introduction}

We study output trajectory tracking for a class of nonlinear reaction diffusion equations such that a prescribed performance of the tracking error is achieved. To this end, we {use} the method of funnel control which was developed in~\cite{IlchRyan02b}, see also the survey~\cite{IlchRyan08}. The funnel controller is a model-free output-error feedback of high-gain type. Therefore, it is inherently robust and of striking simplicity. The funnel controller has been successfully applied e.g.\ in temperature control of chemical reactor models~\cite{IlchTren04}, control of industrial servo-systems~\cite{Hack17} and underactuated multibody systems~\cite{BergOtto19}, speed control of wind turbine systems~\cite{Hack14,Hack15b,Hack17}, current control for synchronous machines~\cite{Hack15a,Hack17}, DC-link power flow control~\cite{SenfPaug14}, voltage and current control of electrical circuits~\cite{BergReis14a}, oxygenation control during artificial ventilation therapy~\cite{PompAlfo14}, control of peak inspiratory pressure~\cite{PompWeye15} and adaptive cruise control~\cite{BergRaue20}.

A funnel controller for a large class of systems described by functional differential equations with arbitrary (well-defined) relative degree {(see~\cite{Isid95} for a definition in the context of nonlinear systems)} has been developed in~\cite{BergHoan18}. It is shown in~\cite{BergPuch20} that this abstract class indeed allows for fairly general infinite-dimensional systems, where the internal dynamics are modeled by a (PDE). In particular, it was shown in~\cite{BergPuch19} that the linearized model of a moving water tank, where sloshing effects appear, belongs to the aforementioned system class. On the other hand, not even every linear, infinite-dimensional system has a well-defined relative degree, in which case the results as in~\cite{BergHoan18,IlchRyan02b} cannot be applied. Instead,  the feasibility of funnel control has to be investigated directly for the (nonlinear) closed-loop system, see~\cite{ReisSeli15b} for a boundary controlled heat equation,~\cite{PuchReis19} for a general class of boundary control systems {and~\cite{Berg20pp} for the Fokker-Planck equation corresponding to a multi-dimensional Ornstein-Uhlenbeck process.}

The nonlinear reaction diffusion system that we consider in the present paper is known as the monodomain model and represents defibrillation processes of the human heart~\cite{Tung78}. The monodomain equations are a reasonable simplification of the well accepted bidomain equations, which arise in cardiac electrophysiology~\cite{SundLine07}. In the monodomain model the dynamics are governed by a parabolic reaction diffusion equation which is coupled with a linear ordinary differential equation that models the ionic current.

It is discussed in~\cite{KuniNagaWagn11} that, under certain initial conditions, reentry phenomena and spiral waves may occur. From a medical point of view, these situations can be interpreted as fibrillation processes of the heart that should be terminated by an external control, for instance by applying an external stimulus to the heart tissue, see~\cite{NagaKuniPlan13}.

The present paper is organized as follows: In Section~\ref{sec:mono_main} we introduce the mathematical framework, which strongly relies on preliminaries on {Robin} elliptic operators. The control objective is presented in Section~\ref{sec:mono_controbj}, where we also state the main result on the feasibility of the proposed controller design in Theorem~\ref{thm:mono_funnel}. The proof of this result is given in Section~\ref{sec:mono_proof_mt} and it uses several auxiliary results derived in Appendices~\ref{sec:mono_prep_proof} and~\ref{sec:mono_prep_proof2}. We illustrate our result by a simulation in Section~\ref{sec:numerics}.

\textbf{Nomenclature}.
The set of bounded operators from $X$ to $Y$ is denoted by $\cL(X,Y)$, $X'$ stands for the dual of a~Banach space $X$, and $B'$ is the dual of an operator $B$.\\
%For $p\in[1,\infty]$, the space of $p$-summable sequences is denoted by $\ell^p$.
For a bounded and measurable set $\Omega\subset\R^d$, $p\in[1,\infty]$ and $k\in\N_0$, $W^{k,p}(\Omega;\R^n)$ denotes the Sobolev space of equivalence classes of $p$-integrable and $k$-times weakly differentiable functions $f:\Omega\to\R^n$, $W^{k,p}(\Omega;\R^n)\cong (W^{k,p}(\Omega))^n$, and the Lebesgue space of equivalence classes of $p$-integrable functions is $L^p(\Omega)=W^{0,p}(\Omega)$. For $r\in(0,1)$ we further set
\[
    W^{r,p}(\Omega) := \setdef{f\in L^p(\Omega)}{ \left( (x,y)\mapsto \frac{|f(x)-f(y)|}{|x-y|^{d/p+r}}\right) \in L^p(\Omega\times\Omega)}.
\]
For a domain $\Omega$ with smooth boundary, $W^{k,p}(\partial\Omega)$ denotes the Sobolev space at the boundary.\\
We identify functions with their restrictions, that is, for instance, if $f\in L^p(\Omega)$  $\Omega_0\subset \Omega$, then the restriction $f|_{\Omega_0}\in L^p(\Omega_0)$ is again {denoted} by~$f$. For an interval $J\subset\R$, a Banach space $X$ and $p\in[1,\infty]$, we denote by $L^p(J;X)$ the vector space of equivalence classes of strongly measurable functions $f:J\to X$ such that $\|f(\cdot)\|_X\in L^p(J)$.
% The space $L^p(J;B)$ is a Banach space with respect to the norm
%\[\|f\|_{L^p(J;B)}:=\begin{cases}
%\left(\int_J\|f(t)\|^p_B\right)^{1/p}&\mbox{if }p\in[1,\infty),\\
%\esssup_{t\in J}\|f(t)\|_B&\mbox{if }p=\infty.
%\end{cases}\]
Note that if $J=(a,b)$ for $a,b\in\R$, the spaces $L^p((a,b);X)$, $L^p([a,b];X)$, $L^p([a,b);X)$ and $L^p((a,b];X)$ coincide, since the points at the boundary have measure zero. We will simply write $L^p(a,b;X)$, also for the case $a=-\infty$ or $b=\infty$. We refer to \cite{Adam75} for further details on Sobolev and Lebesgue spaces.\\
%If $J$ has finite measure, for $p\in[1,\infty)$, $(L^p(J;X))'$ can be isometrically identified with $L^q(J;X')$ with $p^{-1}+q^{-1}=1$. In fact, if $X$ is a Hilbert space, %$L^2(J;X)$ is a Hilbert space with the natural inner product, see \cite[Sec.~IV.1]{Diestel77}.
In the following, let $J\subset\R$ be an interval, $X$ be a Banach space and $k\in\N_0$. Then $C^k(J;X)$ is defined as the space of $k$-times continuously differentiable functions $f:J\to X$. The space of bounded $k$-times continuously differentiable functions with bounded first $k$ derivatives is denoted by $BC^k(J;X)$, and it is a Banach space endowed with the usual supremum norm. The space of bounded and uniformly continuous functions will be denoted by $BUC(J;X)$. The Banach space of H\"older continuous functions $C^{0,r}(J;X)$ with $r\in(0,1)$ is given by
\begin{align*}
C^{0,r}(J;X)&:=\setdef{f\in BC(J;X)}{[f]_{r}:=\sup_{t,s\in J,s<t}\frac{\|f(t)-f(s)\|}{(t-s)^r}<\infty},\\
\|f\|_r&:=\|f\|_\infty+[f]_{r},
\end{align*}
see \cite[Chap.~0]{Luna95}. We like to note that for all $0<r<q<1$ we have that
\[
    C^{0,q}(J;X) \subseteq C^{0,r}(J;X) \subseteq BUC(J;X).
\]
%In the following we refer to .
For $p\in[1,\infty]$, the symbol $W^{1,p}(J;X)$ stands for the Sobolev space of $X$-valued equivalance classes of weakly differentiable and $p$-integrable functions $f:J\to X$ with $p$-integrable weak derivative, i.e., $f,\dot{f}\in L^p(J;X)$. Thereby, integration (and thus weak differentiation) has to be understood in the Bochner sense, see~\cite[Sec.~5.9.2]{Evan10}.
%By \cite[Thm.~5.9.2.2]{Evan10}, for $f\in W^{1,p}(J;X)$ it holds that $f\in C(J;X)$ and for $s,t\in J$, $s\leq t$, we have
%\[f(t)=f(s)+\int_s^t\dot{f}(r)\ds{r}.\]
The spaces $L^p_{\rm loc}(J;X)$ and $W^{1,p}_{\rm loc}(J;X)$ consist of all $f$ whose restriction to any compact interval $K\subset J$ are in $L^p(K;X)$ or $W^{1,p}(K;X)$, respectively.

\section{The FitzHugh-Nagumo model}
\label{sec:mono_main}

Throughout this paper we will frequently use the following assumption. For $d\in\N$ we denote the scalar product in $L^2(\Omega;\R^d)$ by $\scpr{\cdot}{\cdot}$ and the norm in $L^2(\Omega)$ by $\|\cdot\|$.

\begin{Ass}\label{Ass1}
Let $d\le 3$ and $\Omega\subset \R^d$ be a bounded domain with Lipschitz boundary $\partial\Omega$.
%\footnote{By a domain we mean a set which is open and connected.}
%in $, whose boundary $\Gamma$ is satisfies the {\em cone condition}, see \cite[Section 1.2]{Gris85}.
Further, let {$a\in L^\infty(\partial\Omega)$ be nonnegative and} $D\in L^\infty(\Omega;\R^{d\times d})$ be symmetric-valued and satisfy the \emph{ellipticity condition}
\begin{equation}
\exists\,{\beta>0}:\ \ \text{for {a.a.}}\,\zeta\in\Omega\ \forall\, \xi\in\R^d:\ \xi^\top D(\zeta) \xi = \sum_{i,j=1}^d D_{ij}(\zeta)\xi_i\xi_j\geq {\beta} \|\xi\|_{\R^d}^2.\label{eq:ellcond}
\end{equation}
\end{Ass}

To formulate the model of interest, we consider the {Robin elliptic operator $\cA:\cD(\cA)\subset L^2(\Omega)\to L^2(\Omega)$ with
\begin{equation}\label{eq:RobinOp}
\begin{aligned}
&\cD(\cA)=\setdef{\!\!z\in W^{1,2}(\Omega)\!\!}{
{\rm div} D\nabla z\in L^2(\Omega)\,\wedge\,(\nu^\top\cdot D\nabla z+az)|_{\partial\Omega}=0\!\!},\\[1.5ex]
&\cA z={\rm div} D\nabla z,\qquad z\in\cD(\cA),
\end{aligned}
\end{equation}
where $\nu:\partial \Omega\to\R^d$ is the outward normal unit vector. The domain of $\cA$ is well-defined, since any vector field with square integrable weak divergence has a~well-defined normal trace in $W^{1/2,2}(\partial\Omega)'$ (see~\cite[Lem.\ 20.2.]{Tart07}) and any element of $W^{1,2}(\Omega)$
has a~well-defined trace in $W^{1/2,2}(\partial\Omega)$. Note that the Neumann elliptic operator is a~special case of $\cA$ emerging from setting $a=0$.}

The model for the interaction of the electric current in a cell is
\begin{equation}\label{eq:FHN_model}
\begin{aligned}
\dot v(t)&=\cA v(t)+p_3(v(t))-u(t)+I_{s,i}(t)+\cB I_{s,e}(t),\quad&v(0)&=v_0,\\
\dot u(t)&=c_5v(t)-c_4u(t),&u(0)&=u_0,\\
y(t) &= \cB'v(t),
\end{aligned}
\end{equation}
where
\[p_3(v):=-c_1v+c_2v^2-c_3v^3,\]
with constants $c_i>0$ for $i=1,\dots,5$, initial values $v_0,u_0\in L^2(\Omega)$, the {Robin} elliptic operator $\cA:\cD(\cA)\subset L^2(\Omega)\to L^2(\Omega)$ and control operator $\cB\in\cL(\R^m,W^{1,2}(\Omega)')$, where $W^{1,2}(\Omega)'$ is the dual of $W^{1,2}(\Omega)$ with respect to the pivot space $L^{2}(\Omega)$; consequently, $\cB'\in\cL(W^{1,2}(\Omega),\R^m)$.

System~\eqref{eq:FHN_model} is known as the FitzHugh-Nagumo model for the ionic current~\cite{Fitz61}, where
\[I_{ion}(u,v)=p_3(v)-u.\]
The functions  $I_{s,i}\in L^2_{{\rm loc}}(0,T;L^2(\Omega))$, $I_{s,e}\in L^2_{{\rm loc}}(0,T;\R^m)$
are the intracellular and extracellular stimulation currents, respectively. In particular, $I_{s,e}$ is the control input of the system, whereas $y$ is the output.

{Next we record some properties of the Robin elliptic operator which are frequently used throughout this article.
\begin{Rem}\label{Rem:Aop_n}
If Assumption~\ref{Ass1} holds, then the Robin elliptic operator~$\cA$ on~$\Omega$ has the following properties:
\begin{enumerate}[a)]
\item\label{item:Aop3} It follows from \cite[Prop.~3.10]{Nitt11} that there exists some $\nu\in(0,1)$ with $\cD(\cA)\subset C^{0,\nu}(\Omega)$. In particular, $\cD(\cA)\subset L^{\infty}(\Omega)$.
\item
With $\cA$ we may associate the bilinear form
\begin{equation}
\fa:W^{1,2}(\Omega)\times W^{1,2}(\Omega)\to\R,\ (z_1,z_2)\mapsto\scpr{\nabla z_1}{D\nabla z_2}+\scpr{z_1}{az_2}_{L^2(\partial\Omega)},\label{eq:sesq}
\end{equation}
where $\scpr{\cdot}{\cdot}_{L^2(\partial\Omega)}$ denotes the standard inner product in $L^2(\partial\Omega)$.
The relation between the operator~$\cA$ and the form~$\fa$ is revealed by the properties
\begin{equation}\label{eq:faprop}
\begin{aligned}
&\bullet\    \cD(\cA)\!=\!\setdef{\!\!z_2\!\in\! W^{1,2}(\Omega)\!\!}{\!\exists\, y_2\!\in\! L^2(\Omega)\, \forall\, z_1\!\in\! W^{1,2}(\Omega):\,\fa(z_1,z_2)\!=\!-\!\scpr{z_1}{y_2}\!\!\!}\\[1.5ex]
&\bullet\   \forall\, z_1\in W^{1,2}(\Omega)\ \forall\, z_2\in \cD(\cA):\ \fa(z_1,z_2)=-\scpr{z_1}{\cA z_2}.
\end{aligned}
\end{equation}
Furthermore, it follows from Kato's first representation theorem~\cite[Sec.\ VI.2, Thm~2.1]{Kato80} that~$\cA$ is uniquely determined by the properties~\eqref{eq:faprop}, and it is moreover closed and densely defined. The property $\fa(z_1,z_2)={\fa(z_2,z_1)}$ for all $z_1,z_2\in W^{1,2}(\Omega)$
further implies that $\cA$ is self-adjoint.
\item\label{item:Aop4}
The Rellich-Kondrachov theorem~\cite[Thm.~6.3]{Adam75} implies that $\cA$ has compact resolvent. Moreover, we have $\scpr{z}{\cA z}\leq0$ for all $z\in \cD(\cA)$. Combining these findings with self-adjointness of $\cA$ we may infer that
 there exists a real-valued nonnegative and monotonically increasing sequence $(\alpha_j)_{j\in\N_0}$ without accumulation points, such that
 the spectrum of $\cA$ reads $\sigma(\cA)=\setdef{-\alpha_j}{j\in\N_0}$, and there exists an
an orthonormal basis $(\theta_j)_{j\in\N_0}$ of $L^2(\Omega)$, such that
\begin{equation}
    \forall\,x\in\cD(\cA):\ \cA x=-\sum_{j=0}^\infty\alpha_j\scpr{x}{\theta_j}\theta_j.\label{eq:spectr}
\end{equation}
By~\cite[Prop.~3.2.9]{TucsWeis09} the domain of $\cA$ reads
\begin{equation}
\cD(\cA)=\setdef{\sum_{j=0}^\infty \lambda_j \theta_j}{(\lambda_j)_{j\in\N_0}\text{ with } {\sum_{j=0}^\infty (1+\alpha_j^2)} |\lambda_j|^2<\infty}.\label{eq:spectrda}
\end{equation}
\end{enumerate}
\end{Rem}}

Next we introduce the solution concept.

\begin{Def}\label{def:solution}
Let Assumption~\ref{Ass1} hold and $\cA$ be the {Robin} elliptic operator {as in~\eqref{eq:RobinOp}}, let $\cB\in \cL(\R^m,W^{1,2}(\Omega)')$, and
$u_0,v_0\in L^2(\Omega)$ be given. Further, let $T\in(0,\infty]$ and $I_{s,i}\in L^2_{{\rm loc}}(0,T;L^2(\Omega))$, $I_{s,e}\in L^2_{{\rm loc}}(0,T;\R^m)$. A triple of functions $(u,v,y)$ is called \emph{solution} of~\eqref{eq:FHN_model} on $[0,T)$, if
	\begin{enumerate}[(i)]
%		\item $u,v\in C(0,T;L^2(\Omega))$ with $(u(0),v(0))=(u_0,v_0)$;
		\item $v\in {L^2_{\loc}(0,T;W^{1,2}(\Omega))}\cap C([0,T);L^{2}(\Omega))$ with $v(0)=v_0$;
        \item $u\in C([0,T);L^{2}(\Omega))$ with $u(0)=u_0$;
%		\item $\big(t,Cv(t) - \yref(t)\big)\in\cF_\varphi$ for all $t\in(0,T)$;
		\item for all $\chi\in L^2(\Omega)$, $\theta\in W^{1,2}(\Omega)$, the scalar functions $t\mapsto\scpr{u(t)}{\chi}$, $t\mapsto\scpr{v(t)}{\theta}$ are weakly differentiable on $[0,T)$, and for almost all $t\in (0,T)$ we have
		\begin{equation}\label{eq:solution}
		\begin{aligned}
		{\textstyle\ddt}\scpr{v(t)}{\theta}&=-\fa(v(t),\theta)+\scpr{p_3(v(t))-u(t)+I_{s,i}(t)}{\theta}+\scpr{I_{s,e}(t)}{\cB'\theta}_{\R^m},\\
		{\textstyle\ddt} \scpr{u(t)}{\chi}&=\scpr{c_5v(t)-c_4u(t)}{\chi},\\
        y(t)&=\cB'v(t),
%		I_{s,e}(t)&=-\frac{k_0}{1-\varphi(t)^2\|Cv(t)-\yref(t)\|^2_{\R^m}}(Cv(t)-\yref(t))
		\end{aligned}
		\end{equation}
		where $\fa:W^{1,2}(\Omega)\times W^{1,2}(\Omega)\to\R$ is the {bilinear} form defined in \eqref{eq:sesq}.
	\end{enumerate}
%	A solution of~\eqref{eq:FHN_feedback} in~$[0,\infty)$ is called \emph{global}.
\end{Def}

\begin{Rem}\label{rem:openloop}
	\hspace{1em}
\begin{enumerate}[a)]
\item\label{rem:openloop1} Weak differentiability of $t\mapsto\scpr{u(t)}{\chi}$, $t\mapsto\scpr{v(t)}{\theta}$ for all $\chi\in L^2(\Omega)$, $\theta\in W^{1,2}(\Omega)$ on $(0,T)$
further leads to $v\in {W^{1,2}_{\loc}(0,T;W^{1,2}(\Omega)')}$ and $u\in {W^{1,2}_{\loc}(0,T;L^{2}(\Omega))}$.
\item\label{rem:openloop2} The Sobolev embedding theorem \cite[Thm.~5.4]{Adam75} implies that the embedding $W^{1,2}(\Omega)\hookrightarrow L^6(\Omega)$ is bounded. This guarantees that $p_3(v)\in {L^2_{\loc}(0,T;L^2(\Omega))}$, whence the first equation in \eqref{eq:solution} is well-defined.
{
\item\label{rem:openloop3} For later use we investigate when the operator $\cB$ has trivial kernel. By $\cB\in \cL(\R^m,W^{1,2}(\Omega)')$, this operator can be regarded to be composed of $m$ bounded linear functionals $b_1,\ldots,b_m$. More precisely, for $e_i\in\R^m$ being the $i$-th unit vector, we set $b_i:=\cB e_i\in W^{1,2}(\Omega)'$. Then $\cB$ has the representation
    \[\forall\,J_1,\ldots,J_m\in\R:\ \cB(J_1,\ldots,J_m)=J_1 b_1+\ldots+J_mb_m.\]
    As a~consequence, we have that $\ker\cB=\{0\}$ if, and only if, the functionals $b_1,\ldots,b_m$ are linearly independent.}
%\item By Proposition~\ref{prop:Aop2}, we have $X_{1/2}=W^{1,2}(\Omega)$. In particular, for any $\alpha\in[0,1/2]$, $\cB\in\cL(\R^m,X_{-\alpha})$ fulfills
%$\cB\in\cL(\R^m,W^{1,2}(\Omega)')$. We will see later that $\alpha<1/2$ can be exploited to gain more smoothness of the solution.
\end{enumerate}
\end{Rem}

{In the following subsections we discuss two important control frameworks, which can be modelled by a suitable choice of an operator~$\cB$.}

\subsection{Distributed control}\label{Ssec:openloop3}

{
Finite-dimensional distributed control typically means that each point in space is influenced by the same value of the temporal control function according to some predefined spatial shape function.
Hence, in the context of the model~\eqref{eq:FHN_model}, $\cB I_{s,e}(t)$ should be a function of the control input~$I_{s,e}(t)$.
That is, for instance, the case, if $\im B\subset L^2(\Omega)$. In the context of Remark~\ref{rem:openloop}\,\ref{rem:openloop3}), this means that for all $i=1,\ldots,m$ there exist $w_i\in L^{2}(\Omega)$ such that $b_i=\cB e_i=w_i$. In this case, $\cB$ has the form
\[
   \forall\, J_1,\ldots,J_m\in\R:\ \cB(J_1,\ldots,J_m)=J_1 w_1+\ldots+J_mw_m.
\]
By setting $w=(w_1,\ldots,w_m)\in L^2(\Omega)^m$, the output is given by
\[
    y(t)=\cB'v(t) = \int_\Omega (v(t))(\xi)\cdot w(\xi)\, {\rm d}\xi.
\]
A typical situation is that $w_1,\ldots,w_m$ are indicator functions on some subsets of $\Omega$; such choices have been considered in~\cite{KuniSouz18} for instance. Note that $\ker \cB = \{0\}$ if, and only if, the functions $w_1,\ldots,w_m$ are linearly independent in $L^2(\Omega)$, cf.~Remark~\ref{rem:openloop}\,\ref{rem:openloop3}).  For instance, this is satisfied, if $w_1,\ldots,w_m$ are indicator functions on disjoint subsets of $\Omega$.
}

%{
%For later use we investigate when the operator $\cB$ has trivial kernel. We show that, if $w$ is continuous and there exist %$\xi_1,\ldots,\xi_m\in\Omega$ such that $w(\xi_1),\ldots,w(\xi_m)$ are linearly independent, then $\ker \cB = \{0\}$. First observe that, under %this assumption, there exists a matrix $M \in\R^{m\times m}$ such that $M [w(\xi_1),\ldots,w(\xi_m)] = I_m$, i.e.,
%\[
%    \forall\, j=1,\ldots,m:\ M w(\xi_j) = e_j,
%\]
%where $e_j$ denotes the $j$-th unit vector in $\R^m$. If now $\cB u = 0 \in L^2(\Omega)$ for some $u\in\R^m$, then continuity of~$w$ gives that %$w(\xi)^\top u = 0$ for all $\xi\in\Omega$, and in particular
%\[
%    \forall\, j=1,\ldots,m:\ 0 = (M^{-1} M w(\xi_j))^\top u = (M^{-1} e_j)^\top u,
%\]
%by which $u=0$. This shows $\ker \cB = \{0\}$.
%}

%{
%Let us now consider the case that the above condition is not satisfied. For simplicity assume that $m=2$ and $w(\xi_1)$ and $w(\xi_2)$ are %linearly dependent for all $\xi_1,\xi_2\in\Omega$, by which $w(\xi) = \psi(\xi) W$ for some $\psi\in L^2(\Omega)$ and $W\in\R^2$. Clearly, any %$u\in\R^2$ can be written as $u = W u_1 + u_2$ for some $u_1\in\R$ and $u_2\in \ker W^\top \subseteq\R^2$. Then
%\[
%    \mathcal{B} u = \langle w, u\rangle_{\R^m} = \psi \langle W, W u_1 + u_2\rangle_{\R^2}  = \psi \|W\|_{\R^2}^2 u_1,
%\]
%which depends only on~$u_1$, and hence $\mathcal{B}$ can be rewritten as an operator $\mathcal{B}:\R \to L^2(\Omega)$, which satisfies the %condition $\ker \mathcal{B} = \{0\}$. For $m>2$ a similar argument applies.
%}

\subsection{Boundary control}\label{Ssec:openloop4}

{
Boundary control means that the value at the boundary is determined by the control function. In this case, $\cB$ takes values in a distribution space. More precisely, $\cB$ corresponds to a~tuple of linear functionals which assign to $z\in W^{1,2}(\Omega)$ a value which is defined in terms of the trace of $z$ on the boundary $\partial \Omega$. That is, for $\ell\in(0,1/2]$, $w_1,\ldots,w_m\in W^{\ell,2}(\partial\Omega)'$ and $J_1,\ldots,J_m\in\R$ we consider the linear functional
\[
    \cB (J_1,\ldots,J_m)\in W^{\ell+1/2,2}(\Omega)',\;\ z\mapsto \sum_{i=1}^m J_i
    \scpr{w_i}{{\rm tr}(z)}_{{W^{\ell,2}(\partial\Omega)',W^{\ell,2}(\partial\Omega)}},
\]
where
\[
    {\rm tr}:\ z\mapsto \left. z\right|_{\partial\Omega}
%&&X_{1/4+\varepsilon}&\to L^2(\partial\Omega)
\]
is the trace operator, that satisfies ${\rm tr}\in \cL(W^{\ell+1/2,2}(\Omega),W^{\ell,2}(\partial\Omega))$ for all $\ell\in(0,1/2]$ by the trace theorem~\cite[Thm.~9.2.1]{Agr15}. In particular, by the continuity of the embedding $W^{1,2}(\Omega)\subset W^{\ell+1/2,2}(\Omega)$, the mapping
\[
    \cB: (J_1,\ldots,J_m) \mapsto \cB (J_1,\ldots,J_m)
\]
is linear and bounded from $\R^m$ to $W^{1,2}(\Omega)'$. Note that the kernel of $\cB$ is trivial if, and only if, the boundary functionals $w_1,\ldots,w_m$ are linearly independent, cf.~Remark~\ref{rem:openloop}\,\ref{rem:openloop3}).\\
In the context of the model~\eqref{eq:FHN_model}, the operator~$\cB$ corresponds to a~Robin boundary control
\[
   \left. \nu^\top\cdot D\nabla v(t) + a \, v(t)\right|_{\partial\Omega} =  \scpr{w}{I_{s,e}(t)}_{\R^m},
\]
where $w = (w_1,\ldots,w_m)$. In this case, the output is given by the evaluation of the Dirichlet boundary values of $v(t)$ at $w_1,\ldots,w_m$. More precisely,
\[
    y(t)= \cB' v(t) = \begin{pmatrix}\scpr{w_1}{{\rm tr} (v(t))}_{{W^{\ell,2}(\partial\Omega)',W^{\ell,2}(\partial\Omega)}}\\\vdots\\
    \scpr{w_m}{{\rm tr} (v(t))}_{{W^{\ell,2}(\partial\Omega)',W^{\ell,2}(\partial\Omega)}}\end{pmatrix}.
\]
By taking into account that
$W^{\ell,2}(\partial\Omega)\subset L^2(\partial\Omega)$ via a~canonical embedding, a~special case is where $w=(w_1,\ldots,w_m)\in L^2(\partial\Omega)^m$, in which, for $J\in\R^m$, $\cB$ reads
\[
    \cB J \in W^{1,2}(\Omega)',\ z\mapsto \int_{\partial \Omega} \langle w(\xi),J\rangle_{\R^m}\, {\rm tr}(z(\xi))\, {\rm d} \sigma,
\]
and the output is given by the weighted integral of the Dirichlet boundary values
\[
    y(t)=\int_{\partial\Omega} w(\xi)\cdot(v(t))(\xi)\, {\rm d}\sigma.
\]
Again, $\ker \cB = \{0\}$ if, and only if, the boundary functionals $w_1,\ldots,w_m$ are linearly independent.  This is satisfied, for instance, if $w_1,\ldots,w_m$ are indicator functions on disjoint subsets of $\partial \Omega$.
}

\section{Control objective}\label{sec:mono_controbj}

The objective is that the output~$y$ of the system~\eqref{eq:FHN_model} tracks a given reference signal which is $y_{\rm ref}\in W^{1,\infty}(0,\infty;\R^m)$ with a prescribed performance of the tracking error $e:= y- y_{\rm ref}$, that is~$e$ evolves within the performance funnel
\[
\cF_\varphi := \setdef{ (t,e)\in[0,\infty)\times\R^m}{ \varphi(t) \|e\|_{\R^m} < 1}
\]
defined by a function~$\varphi$ belonging to
\[
\Phi_\gamma:=\setdef{\varphi\in W^{1,\infty}(0,\infty;\R) }{\varphi|_{[0,\gamma]}\equiv0,\ \forall\delta>0, \inf_{t>\gamma+\delta}\varphi(t)>0
},
\]
for some $\gamma>0$.

\begin{figure}[h]
	\begin{minipage}{0.45\textwidth}
		\begin{tikzpicture}[domain=0.001:4,scale=1.2] % Zeichenbereich
		% Gitter zeichnen
		%\draw[very thin,color=gray] (-0.1,-1.1) grid (3.9,3.9);
		% Füllung zeichnen
		\fill[color=blue!20,domain=0.47:4] (0,0)-- plot (\x,{min(2.2,1/\x+2*exp(-3))})--(4,0)-- (0,0);
		\fill[color=blue!20] (0,0) -- (0,2.2) -- (0.47,2.2) -- (0.47,0) -- (0,0);
		%-- plot (\x,{-2*exp(-0.5*\x)})--(4,0);
		% Achsen zeichnen
		\draw[->] (-0.2,0) -- (4.3,0) node[right] {$t$};
		\draw[->] (0,-0.2) -- (0,2.5) node[above] {};
		% Funktionen zeichnen
		\draw[color=blue,domain=0.47:4] plot (\x,{min(2.2,1/\x+2*exp(-3))}) node[above] {$1/\varphi(t)$};
		%\draw[color=blue] plot (\x,{-2*exp(-0.5*\x)});
		\draw[smooth,color=red,style=thick] (1,0) node[below] {$\|e(t)\|_{\R^m}$}
		plot coordinates{(0,0.8)(0.5,1.2)(1,0)}--
		plot coordinates{(1,0)(1.25,0.6)(2,0.2)(3,0.3)(4,0.05)} ;
		\end{tikzpicture}
		\caption{Error evolution in a funnel $\mathcal F_{\varphi}$ with boundary~$1/\varphi(t)$.}
		\label{Fig:monodomain_funnel}
	\end{minipage}
	\quad
	\begin{minipage}{0.5\textwidth}
		The situation is illustrated in Fig.~\ref{Fig:monodomain_funnel}. The funnel boundary given by~$1/\varphi$ is unbounded in a small interval $[0,\gamma]$ to allow for an arbitrary initial tracking error. Since $\varphi$ is bounded there exists $\lambda>0$ such that $1/\varphi(t) \ge \lambda$ for all $t>0$. Thus, we seek practical tracking with arbitrary small accuracy $\lambda>0$, but asymptotic tracking is not required in general.
	\end{minipage}
\end{figure}
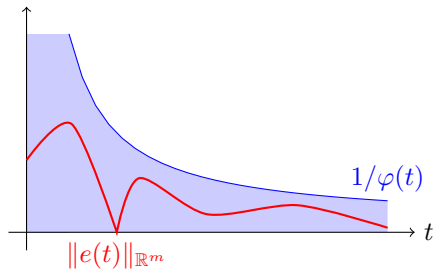
The funnel boundary is not necessarily monotonically decreasing, while in most situations
it is convenient to choose a monotone funnel. Sometimes, widening the funnel over some later time interval might be beneficial, for instance in the presence of periodic disturbances or strongly varying reference signals. For typical choices of funnel boundaries see e.g.~\cite[Sec.~3.2]{Ilch13}.\\
A controller which achieves the above described control objective is the funnel controller. In the present paper, it suffices to restrict ourselves to the simple version developed in~\cite{IlchRyan02b}, which is the feedback law
\begin{equation}\label{eq:monodomain_funnel_controller}
I_{s,e}(t)=-\frac{k_0}{1-\varphi(t)^2\|\cB' v(t)-\yref(t)\|^2_{\R^m}}(\cB'v(t)-\yref(t)),
\end{equation}
where $k_0>0$ is some constant used for scaling and agreement of physical units. Note that, by $\varphi|_{[0,\gamma]}\equiv0$, the controller satisfies
\[
    \forall\,  t\in[0,\gamma]:\ I_{s,e}(t)=-k_0(\cB'v(t)-\yref(t)).
\]
{Inserting the feedback law~\eqref{eq:monodomain_funnel_controller} into the system~\eqref{eq:FHN_model}, we obtain the closed-loop system
\begin{equation}\label{eq:FHN_feedback}
\begin{aligned}
\dot v(t)&=\cA v(t)+p_3(v)(t)-u(t)+I_{s,i}(t) -\frac{k_0 \cB (\cB' v(t)-\yref(t))}{1-\varphi(t)^2\|\cB'v(t)-\yref(t)\|^2_{\R^m}},\\
\dot u(t)&=c_5v(t)-c_4u(t),
\end{aligned}
\end{equation}
for which we seek to show existence and uniqueness of global solutions~-- this is the subject of the main result Theorem~\ref{thm:mono_funnel} below. Note that} the system~\eqref{eq:FHN_feedback} is a nonlinear and non-autonomous PDE and any solution needs to satisfy that the tracking error evolves in the prescribed performance funnel~$\cF_\varphi$. Therefore, existence and uniqueness of solutions is a nontrivial problem and even if a solution exists on a finite time interval $[0,T)$, it is not clear that it can be extended to a global solution.

{We introduce} the following weak solution framework.

\begin{Def}\label{def:solution_feedback}
Use the assumptions from Definition~\ref{def:solution}. Furthermore, let $k_0>0$, $\yref\in W^{1,\infty}(0,\infty;\R^m)$, $\gamma>0$ and $\varphi\in\Phi_\gamma$. A triple of functions $(u,v,y)$ is called \emph{solution} of {system~\eqref{eq:FHN_feedback}} on $[0,T)$, if $(u,v,y)$ satisfies the conditions~(i)--(iii) from Definition~\ref{def:solution} with $I_{s,e}$ as in~\eqref{eq:monodomain_funnel_controller}.
\end{Def}

\begin{Rem}\
	\begin{enumerate}[a)]
		\item
{To be precise,} $(u,v,y)$ is a {solution} {of~\eqref{eq:FHN_feedback} on $[0,T)$} if, and only if,
	\begin{enumerate}[(i)]
%		\item $u,v\in C(0,T;L^2(\Omega))$ with $(u(0),v(0))=(u_0,v_0)$;
		\item $v\in {L^2_{\loc}(0,T;W^{1,2}(\Omega))}\cap C([0,T);L^{2}(\Omega)))$ with $v(0)=v_0$;
        \item $u\in C([0,T);L^{2}(\Omega))$ with $u(0)=u_0$;
%		\item $\big(t,Cv(t) - \yref(t)\big)\in\cF_\varphi$ for all $t\in(0,T)$;
		\item for all $\chi\in L^2(\Omega)$, $\theta\in W^{1,2}(\Omega)$, the scalar functions $t\mapsto\scpr{u(t)}{\chi}$, $t\mapsto\scpr{v(t)}{\theta}$ are weakly differentiable on $[0,T)$, and it holds that, for almost all $t\in (0,T)$,
		\begin{equation}\label{eq:solution_cl}
		\begin{aligned}
		{\textstyle\ddt}\scpr{v(t)}{\theta}&=-\fa(v(t),\theta)+\scpr{p_3(v(t))-u(t)+I_{s,i}(t)}{\theta}\\&\qquad-{\frac{k_0 \scpr{\cB'v(t)-\yref(t)}{\cB'\theta}_{\R^m}}{1-\varphi(t)^2\|\cB' v(t)-\yref(t)\|^2_{\R^m}}},\\[2mm]
		{\textstyle\ddt}\scpr{u(t)}{\chi}&=\scpr{c_5v(t)-c_4u(t)}{\chi},\\
        y(t)&=\cB'v(t).
%		I_{s,e}(t)&=-\frac{k_0}{1-\varphi(t)^2\|Cv(t)-\yref(t)\|^2_{\R^m}}(Cv(t)-\yref(t))
		\end{aligned}
		\end{equation}
	\end{enumerate}
		\item For global solutions it is desirable that $I_{s,e}\in L^\infty(\delta,\infty;\R^m)$ for all $\delta>0$. Note that this is equivalent to
		\[\limsup_{t\to\infty}\varphi(t)^2\|\cB'v(t)-\yref(t)\|^2_{\R^m}<1.\]
        {Furthermore, the output $y$ is a signal that is measured and the control input $I_{s,e}$ is a signal used to manipulate the system, which hence must be generated by a certain device. For both measurement and generation of signals to be feasible it is desirable to have a certain regularity.}
	\end{enumerate}
\end{Rem}

In the following we state the main result of the present paper. We will show that the closed-loop system~\eqref{eq:FHN_feedback}
has a unique global solution so that all signals remain bounded. Furthermore, the tracking error stays uniformly away from the funnel boundary.
We further show that we gain more regularity of the solution, if $\cB\in\cL(\R^m,W^{r,2}(\Omega)')$ for some $r\in [0,1)$ or even
$\cB\in\cL(\R^m,W^{1,2}(\Omega))$. Recall that $\cB\in\cL(\R^m,W^{r,2}(\Omega)')$ if, and only if, $\cB'\in\cL(W^{r,2}(\Omega),\R^m)$. Furthermore, for any $r\in(0,1)$ we have the inclusions
\[
    \cL(\R^m,W^{1,2}(\Omega)) \subset \cL(\R^m,L^2(\Omega)) \subset \cL(\R^m,W^{r,2}(\Omega)') \subset \cL(\R^m,W^{1,2}(\Omega)').
\]

\begin{Thm}\label{thm:mono_funnel}
Use the assumptions from Definition~\ref{def:solution_feedback}. Furthermore, assume that $\ker\cB=\{0\}$ and $I_{s,i}\in L^\infty(0,\infty;L^2(\Omega))$. Then there exists a unique solution of~\eqref{eq:FHN_feedback} on $[0,\infty)$ and we have
	\begin{enumerate}[(i)]
		\item\label{thm:mono_funnel1} $u,\dot{u},v\in BC([0,\infty);L^2(\Omega))$;
		\item\label{thm:mono_funnel2} for all $\delta>0$ we have
		\begin{align*}
			v&\in BUC([\delta,\infty);W^{1,2}(\Omega))\cap C^{0,1/2}([\delta,\infty);L^{2}(\Omega)),\\
			y,I_{s,e}&\in BUC([\delta,\infty);\R^m);
		\end{align*}
        \item\label{thm:mono_funnel3} $\exists\,\varepsilon_0>0\ \forall\,\delta>0\ \forall\, t\geq\delta:\ \varphi(t)^2\|\cB'v(t)-\yref(t)\|^2_{\R^m}\leq1-\varepsilon_0.$
		\end{enumerate}
        Furthermore,
		\begin{enumerate}[a)]
			\item if additionally $\cB\in\cL(\R^m,W^{r,2}(\Omega)')$ for some $r\in (0,1)$, then for all $\delta>0$ we have that
			\begin{align*}
			v\in C^{0,1-r/2}([\delta,\infty);L^{2}(\Omega)),\quad
			y,I_{s,e}\in C^{0,1-r}([\delta,\infty);\R^m).
			\end{align*}
			\item if additionally $\cB\in\cL(\R^m,L^2(\Omega))$, then for all $\delta>0$ and all $\lambda\in(0,1)$ we have
			\begin{align*}
			v\in C^{0,\lambda}([\delta,\infty);L^{2}(\Omega)),\quad
			y,I_{s,e}\in C^{0,\lambda}([\delta,\infty);\R^m).
			\end{align*}
			\item if additionally $\cB\in\cL(\R^m,W^{1,2}(\Omega))$, then for all $\delta>0$ we have $y,I_{s,e}\in W^{1,\infty}([\delta,\infty);\R^m)$.
	\end{enumerate}
\end{Thm}

\begin{Rem}\label{rem:main}
	\hspace{1em}
\begin{enumerate}[a)]
\item\label{rem:main1} The condition $\ker \cB=\{0\}$ is equivalent to $\im \cB'$ being dense in $\R^m$. The latter is equivalent to $\im \cB'=\R^m$ by the finite-dimensionality of $\R^m$.\\
    Note that surjectivity of $\cB'$ is mandatory for tracking control, since it is necessary that any reference signal $\yref\in W^{1,\infty}(0,\infty;\R^m)$ can actually be generated by the output $y(t)=\cB' v$. This property is sometimes called \emph{right-invertibility}, see e.g.~\cite[Sec.~8.2]{TrenStoo01}.
\item\label{rem:main2} {If the input operator corresponds to Robin boundary control as discussed in Section~\ref{Ssec:openloop4}, then $\ker \cB = \{0\}$ and $\cB\in\cL(\R,W^{\ell+1/2,2}(\Omega)')$ for some $\ell\in(0,1/2]$, thus the assertions of Theorem~\ref{thm:mono_funnel}\,a) hold.}
\item\label{rem:main3} {If the input operator corresponds to distributed control as discussed in Section~\ref{Ssec:openloop3}, then $\ker \cB = \{0\}$ and $\cB\in\cL(\R,L^2(\Omega))$, thus the assertions of Theorem~\ref{thm:mono_funnel}\,b) hold.}
{\item\label{rem:main4} The proof of Theorem~\ref{thm:mono_funnel}, carried out in Section~\ref{sec:mono_proof_mt}, exploits the properties of the Robin elliptic operator~$\cA$, which is self-adjoint, nonpositive, has compact resolvent and satisfies $\cD(\cA)\subset L^{\infty}(\Omega)$ (see Remark~\ref{Rem:Aop_n}\,\ref{item:Aop3}). Furthermore, the associated bilinear form~$\fa$ is defined on a subset of $L^6(\Omega)$ (see Remark~\ref{rem:openloop}\,\ref{rem:openloop2}). In principle, the Robin elliptic operator can be replaced by an arbitrary operator~$\cA$ with the aforementioned properties. In this case, Theorem~\ref{thm:mono_funnel} is still valid~-- with the slight modification that the expression $W^{1,2}(\Omega)$ in~\eqref{thm:mono_funnel2} has to be replaced with the interpolation space $(L^2(\Omega),\cD(\cA))_{1/2}$ (see Definition~\ref{def:interp}).\\
    For instance, an elliptic operator ${\rm div} D\nabla$ with domain including homogeneous Dirichlet boundary conditions exhibits the above mentioned properties as long as~$\partial\Omega$ and~$D$ are sufficiently smooth.}
\end{enumerate}
\end{Rem}

%Before we begin to develop the necessary results to prove Theorem \ref{thm:mono_funnel}, we show the simulated system \eqref{eq:FHN_feedback}. We have used the parameter values [...]. The results are given in Fig.

%\ref{fig:distcont}.
%
%\begin{figure}[ht!!]
%	\centering
%	\includegraphics[width=350pt]{funnel_FHN.pdf}
%	\caption{Numerical simulation of (\ref{eq:FHN_feedback}). $\overline{I}_{ion}$ corresponds to the averaged quantity in $\Omega$ of $I_{ion}$.}\label{fig:distcont}
%\end{figure}

\section{Proof of Theorem~\ref{thm:mono_funnel}}
\label{sec:mono_proof_mt}

The proof is inspired by~the results of \cite{Jack90} on existence and uniqueness of (non-controlled) FitzHugh-Nagamo equations, which is based on a spectral approximation and subsequent convergence proofs by using arguments from~\cite{Lion69}. We divide the proof in two major parts. First, we show that there exists a unique solution on the interval $[0,\gamma]$. After that we show that the solution also exists on $(\gamma,\infty)$, is continuous at $t=\gamma$ and has the desired properties.

\subsection{Solution on $[0,\gamma]$}
\label{ssec:mono_proof_tleqgamma}

Assuming that $t\in[0,\gamma]$, we have that $\varphi(t)\equiv0$ so that we need to show existence of a pair of functions $(v,u)$ with the properties as in Definition~\ref{def:solution}~(i)--(iii), where~\eqref{eq:solution} simplifies to
\begin{equation}\label{eq:weak_uv_delta}
\begin{aligned}
\ddt\scpr{v(t)}{\theta}&=-\fa(v(t),\theta)+\scpr{p_3(v(t))-u(t)+I_{s,i}(t)}{\theta}+\scpr{I_{s,e}(t)}{\cB'\theta}_{\R^m},\\
\ddt\scpr{u(t)}{\chi}&=\scpr{c_5 v(t)-c_4u(t)}{\chi},\\
I_{s,e}(t)&=-k_0(\cB' v(t)-\yref(t)),\\
y(t)&= \cB' v(t).
\end{aligned}
\end{equation}
Recall that $\fa:W^{1,2}(\Omega)\times W^{1,2}(\Omega)\to\R$ is the {bilinear} form \eqref{eq:sesq}.

\emph{Step 1: We show existence and uniqueness of a solution.}\\
\emph{Step 1a: We show existence of a local solution on $[0,\gamma]$.} To this end, let $(\theta_i)_{i\in\N_0}$ be the eigenfunctions of $-\cA$ and $\alpha_i$ be the corresponding eigenvalues, with $\alpha_i\geq0$ for all $i\in\N_0$. Recall that $(\theta_i)_{i\in\N_0}$ {forms} an orthonormal basis of $L^2(\Omega)$ by Remark~\ref{Rem:Aop_n}\,\ref{item:Aop4}). Hence, with $a_i :=  \scpr{v_0}{\theta_i}$ and $b_i :=  \scpr{u_0}{\theta_i}$ for $i\in\N_0$ and
\[
v_0^n:= \sum_{i=0}^na_{i}\theta_i,\quad u_0^n:= \sum_{i=0}^nb_{i}\theta_i,\quad n\in\N,
\]
we have that $v^n_0\to v_0$ and $u^n_0\to u_0$ strongly in $L^2(\Omega)$.\\
Fix $n\in\N_0$ and let $\gamma_i:=\cB'\theta_i$ for $i=0,\dots,n$. Consider, for $j=0,\ldots,n$, the differential equations
\begin{equation}\label{eq:muj-nuj-0gamma}
\begin{aligned}
\dot{\mu}_j(t)&=-\alpha_j\mu_j(t)-\nu_j(t)-\scpr{k_0\left(\sum_{i=0}^n\gamma_i\mu_i(t)-\yref(t)\right)}{\gamma_j}_{\R^m}+\scpr{I_{s,i}(t)}{\theta_j} \\
&\quad +\scpr{p_3\left(\sum_{i=0}^n\mu_i(t)\theta_i\right)}{\theta_j},\\
\dot{\nu}_j(t)&=-c_4\nu_j(t)+c_5\mu_j(t),\qquad\qquad \text{with}\ \mu_j(0)=a_j,\ \nu_j(0)=b_j,
\end{aligned}
\end{equation}
defined on $\D:=[0,\infty)\times\R^{2(n+1)}$. {Given that} the functions {defining the system of ODEs \eqref{eq:muj-nuj-0gamma}} are continuous, it follows from ODE theory, see e.g.~\cite[\S~10, Thm.~XX]{Walt98}, that there exists a weakly differentiable solution $(\mu^n,\nu^n)=(\mu_0,\dots,\mu_n,\nu_0,\dots,\nu_n):[0,T_n)\to\R^{2(n+1)}$ of~\eqref{eq:muj-nuj-0gamma} such that $T_n\in(0,\infty]$ is maximal. Furthermore, the closure of the graph of~$(\mu^n,\nu^n)$ is not a compact subset of $\D$.\\
Now, set $v_n(t):=\sum_{i=0}^n{\mu}_i(t)\theta_i$ and $u_n(t):=\sum_{i=0}^n{\nu}_i(t)\theta_i$. {We intend to show  that $(v_n)$ and $(u_n)$ have subsequences which weakly converge to solutions of~\eqref{eq:FHN_feedback} on $[0,\gamma]$.} Invoking~\eqref{eq:muj-nuj-0gamma} and using the functions $\theta_j$ we have that for $j=0,\dots,n$ the functions $(v_n,u_n)$ satisfy
	\begin{equation}\label{eq:weak_delta}
	\begin{aligned} \scpr{\dot{v}_n(t)}{\theta_j}&=-\fa(v_n(t),\theta_j)- \scpr{u_n(t)}{\theta_j}+\scpr{p_3(v_n(t))}{\theta_j}+\scpr{I_{s,i}(t)}{\theta_j}\\
	&\quad -\scpr{k_0(\cB' v_n(t)-\yref(t))}{\cB'\theta_j}_{\R^m},\\
    \scpr{\dot{u}_n(t)}{\theta_j} &= -c_4\scpr{u_n(t)}{\theta_j}+c_5\scpr{v_n(t)}{\theta_j}.
	\end{aligned}
	\end{equation}
\emph{Step 1b: We show {the} boundedness of $(v_n,u_n)$.} Consider the Lyapunov function candidate
\begin{equation}\label{eq:Lyapunov}
V:L^2(\Omega)\times L^2(\Omega)\to\R,\ (v,u)\mapsto\frac12(c_5\|v\|^2+\|u\|^2).
\end{equation}
Observe that, since {the family} $(\theta_i)_{i\in\N_0}$ {is} orthonormal, we have $\|v_n\|^2 = \sum_{j=0}^n \mu_j^2$ and $\|u_n\|^2 = \sum_{j=0}^n \nu_j^2$. Hence we find that, for all $t\in[0, T_n)$,
\begin{align*}
\ddt V(v_n(t),u_n(t)) &\stackrel{\eqref{eq:muj-nuj-0gamma}}{=} c_5\sum_{j=0}^n\mu_j(t)\dot{\mu}_j(t)+\sum_{j=0}^n\nu_j(t)\dot{\nu}_j(t)\\ &=-c_5\sum_{j=0}^{n}\alpha_j\mu_j(t)^2-c_4\sum_{j=0}^{n}\nu_j(t)^2\\
&\quad -c_5\scpr{k_0\left(\sum_{i=0}^n\gamma_i\mu_i(t)-\yref(t)\right)}{\sum_{j=0}^n\gamma_j\mu_j(t)}_{\R^m}\\
&\quad +c_5\scpr{p_3\left(v_n(t)\right)}{v_n(t)}+c_5\scpr{I_{s,i}(t)}{v_n(t)}
\end{align*}
hence, omitting the argument~$t$ for brevity in the following,
\begin{equation}\label{eq:Lyapunov_delta}
\begin{aligned}
\ddt V(v_n,u_n)=&-c_5\fa(v_n,v_n)-c_4\|u_n\|^2+c_5\scpr{I_{s,i}}{v_n}\\
&-c_5k_0\|\overline{e}_n\|_{\R^m}^2+c_5k_0\scpr{\overline{e}_n}{\yref}_{\R^m}+c_5\scpr{p_3(v_n)}{v_n},
\end{aligned}
\end{equation}
where
\[\overline{e}_n(t):=\sum_{i=0}^n\gamma_i\mu_i(t)-\yref(t)=\cB' v_n(t)-\yref(t).\]
Before proceeding, recall Young's inequality for products, i.e., for $a,b\ge 0$ and $p,q\ge 1$ such that $1/p + 1/q = 1$ we have that
\[
a b \le \frac{a^p}{p} + \frac{b^q}{q},
\]
which will be frequently used in the following. Note that
\begin{align*}
\scpr{p_3(v_n)}{v_n}&=-c_1\|v_n\|^2+c_2\scpr{v_n^2}{v_n}-c_3\|v_n\|^4_{L^4},\\
c_2|\scpr{v_n^2}{v_n}|&=|\scpr{\epsilon v_n^3}{\epsilon^{-1}c_2}|\leq \frac{3\epsilon^{4/3}}{4}\|v_n\|^4_{L^4}+\frac{c_2^4}{4\epsilon^4}|\Omega|,
\end{align*}
where the latter follows from Young's inequality with $p=\tfrac43$ and $q=4$. Choosing $\epsilon=\left(\tfrac23 c_3\right)^{\tfrac34}$ we obtain
\[\scpr{p_3(v_n)}{v_n}\leq \frac{27 c_2^4}{32 c_3^3}|\Omega|-c_1\|v_n\|^2-\frac{c_3}{2}\|v_n\|_{L^4}^4.\]
Moreover,
\[\scpr{\overline{e}_n}{\yref}_{\R^m}\leq\frac{1}{2}\|\overline{e}_n\|_{\R^m}^2+\frac{1}{2}\|\yref\|_{\R^m}^2{ \le \frac{1}{2}\|\overline{e}_n\|_{\R^m}^2+\frac{1}{2}\|\yref\|_{\infty}^2}\]
and
\[\scpr{I_{s,i}}{v_n}\leq \frac{c_1}{2}\|v_n\|^2+\frac{1}{2c_1}\|I_{s,i}\|^2 {\le  \frac{c_1}{2}\|v_n\|^2+\frac{1}{2c_1}\|I_{s,i}\|^2_{2,\infty}},\]
{where $\|I_{s,i}\|_{2,\infty} = \esssup_{t\ge 0} \left( \int_\Omega |I_{s,i}(\zeta,t)|^2 \ds{\zeta}\right)^{1/2}$, so} that {with
\[C_\infty:=\frac{k_0c_5}{2}\|\yref\|_{\infty}^2+{\frac{c_5}{2c_1}}\|I_{s,i}\|^2_{2,\infty}+\frac{27 c_2^4}{32 c_3^3}|\Omega|,\]
we may further estimate~\eqref{eq:Lyapunov_delta}} by
\begin{align*}
	\ddt V(v_n,u_n) \leq&-c_5\fa(v_n,v_n){-c_4 \| u_n\|^2}-\frac{c_1c_5}{2}\|v_n\|^2-\frac{c_5k_0}{2}\|\overline{e}_n\|_{\R^m}^2\\
&-\frac{c_3 c_5}{2}\|v_n\|_{L^4}^4+{C_\infty}\\
	\leq& -c_5\fa(v_n,v_n)-\frac{c_1c_5}{2}\|v_n\|^2-\frac{c_5k_0}{2}\|\overline{e}_n\|_{\R^m}^2-\frac{c_3 c_5}{2}\|v_n\|_{L^4}^4+{C_\infty}.
\end{align*}
{Then} we obtain that, for all $t\in[0, T_n)$,
\begin{align*}
 &V(v_n(t),u_n(t))+c_5\int_0^t\fa(v_n(s),v_n(s))\ds{s}+\frac{c_1 c_5}{2}\int_0^t\|v_n(s)\|^2\ds{s}\\
	&+\frac{c_5k_0}{2}\int_0^t\|\overline{e}_n(s)\|_{\R^m}^2\ds{s} +\frac{c_3 c_5}{2}\int_0^t\|v_n(s)\|_{L^4}^4\ds{s}\leq V(v_0^n,u_0^n)+C_\infty t.
\end{align*}
Since $(u_n^0,v_n^0)\to (u_0,v_0)$ strongly in $L^2(\Omega)$ and {since} we have for all $p\in L^2(\Omega)$ that {$ \left\|\sum_{i=0}^n \scpr{p}{\theta_i} \theta_i\right\|^2 \le \| p\| ^2$,}
% \[
% \left\|\sum_{i=0}^n \scpr{p}{\theta_i} \theta_i\right\|^2=\sum_{i=0}^n \scpr{p}{\theta_i}^2\leq\sum_{i=0}^\infty \scpr{p}{\theta_i}^2=\left\|\sum_{i=1}^\infty \scpr{p}{\theta_i}\theta_i\right\|^2=\|p\|^2,
% \]
it follows that, for all $t\in[0, T_n)$,
\begin{equation}\label{eq:bound_1}
\begin{aligned}
&c_5\|v_n(t)\|^2+\|u_n(t)\|^2+ 2c_5\int_0^t\fa(v_n(s),v_n(s))\ds{s}+ {c_1 c_5}\int_0^t\|v_n(s)\|^2\ds{s} \\
&+ {c_5k_0}\int_0^t\|\overline{e}_n(s)\|_{\R^m}^2\ds{s} +{c_3 c_5}\int_0^t\|v_n(s)\|_{L^4}^4\ds{s}\leq 2C_\infty t+c_5\|u_0\|^2+\|v_0\|^2.
\end{aligned}
\end{equation}
\emph{Step 1c: We show that $T_n = \infty$.} Assume that $T_n<\infty$, then it follows from~\eqref{eq:bound_1} together {with~\eqref{eq:ellcond} and~\eqref{eq:sesq}} that $(v_n,u_n)$ is bounded, thus the solution $(\mu^n,\nu^n)$ of~\eqref{eq:muj-nuj-0gamma} is bounded on $[0,T_n)$. But this implies that the closure of the graph of $(\mu^n,\nu^n)$ is a compact subset of $\D$, a contradiction. Therefore, $T_n=\infty$ and in particular the solution is defined for all $t\in[0,\gamma]$.\\
\emph{Step 1d: We show convergence of $(v_n,u_n)$ to  a solution of~\eqref{eq:weak_uv_delta} on $[0,\gamma]$.} First note that it follows from~\eqref{eq:bound_1} that
\begin{equation}\label{eq:uv_bound_delta}
\forall\, t\in[0,\gamma]:\ \|v_n(t)\|^2\leq C_v,\quad \|u_n(t)\|^2\leq C_u
\end{equation}
for some $C_v, C_u>0$. From \eqref{eq:bound_1} and condition~\eqref{eq:ellcond} in Assumption~\ref{Ass1}  it follows that there is a constant $C_\delta>0$ such that
\[\int_0^\gamma\|\nabla v_n(t)\|^2\ds{t}\leq\delta^{-1}\int_0^\gamma\fa(v_n(t),v_n(t))\ds{t}\leq C_\delta.\]
This together with~\eqref{eq:bound_1}  and~\eqref{eq:uv_bound_delta}  implies that there exist constants $C_1,C_2>0$ with
\begin{equation}\label{eq:extra_bounds_delta}
\|v_n\|^4_{L^4(0,\gamma;L^{4}(\Omega))}\leq C_1,\quad \|v_n\|_{L^2(0,\gamma;W^{1,2}(\Omega))}\leq C_2.
\end{equation}
Note that \eqref{eq:extra_bounds_delta} directly implies that
\begin{equation}\label{eq:vn2}
\begin{aligned}
\|v_n^2\|^2_{L^2(0,\gamma;L^{2}(\Omega))} = {\|v_n\|^4_{L^4(0,\gamma;L^{4}(\Omega))}}\leq&\, C_1,\\
\|v_n^3\|_{L^{4/3}(0,\gamma;L^{4/3}(\Omega))} =\, \left( \|v_n^2\|^2_{L^2(0,\gamma;L^{2}(\Omega))}\right)^{3/4} \le &\, C_1^{3/4}.
\end{aligned}\end{equation}
Multiplying the second equation in~\eqref{eq:weak_delta} by $\dot{\nu}_j$ and {summing over} $j\in\{0,\ldots,n\}$ leads to
\begin{align*}
    \|\dot{u}_n\|^2 &=  -\frac{c_4}{2}  \ddt \|u_n\|^2 + c_5 \scpr{v_n}{\dot{u}_n}\\
    &\le -\frac{c_4}{2} \ddt\|u_n\|^2 + \frac{c_5^2}{2}\|v_n\|^2+\frac{1}{2}\|\dot{u}_n\|^2,
\end{align*}
thus
\[\|\dot{u}_n\|^2\leq-c_4 \ddt\|u_n\|^2+c_5^2\|v_n\|^2.\]
Upon integration over $[0,\gamma]$ and using \eqref{eq:uv_bound_delta} this yields that
\[\int_0^\gamma\|\dot{u}_n(t)\|^2\ds{t}\leq c_4 C_u + c_5^2\int_0^\gamma\|v_n(t)\|^2\ds{t} \le \hat{C}_3\]
for some $\hat{C}_3>0$, where the last inequality is a consequence of~\eqref{eq:bound_1}. This together with \eqref{eq:uv_bound_delta} implies that there is $C_3>0$ such that $\|u_n\|_{W^{1,2}(0,\gamma;L^2(\Omega))}\leq C_3$.

Now, let $P_n$ be the orthogonal projection of $L^2(\Omega)$ onto the subspace generated by the set $\setdef{\theta_i}{i=1,\dots,n}$. Consider
\[
    \|v\|_{W^{1,2}}= {\big( \|v\|^2 + \fa(v,v)\big)^{1/2} = }\left(\sum_{i=0}^n(1+\alpha_i)|\scpr{v}{\theta_i}|^2\right)^{1/2},\quad {v\in P_n\big(L^2(\Omega)\big),}
\]
{which is -- by Remark \ref{rem:X_alpha} -- indeed a~norm
on (the projection of) $W^{1,2}(\Omega)$ which is equivalent to the standard norm on $W^{1,2}(\Omega)$ by the properties of~$D$ and~$a$ in Assumption~\ref{Ass1}.} By duality we have that
\[
    \|\hat v\|_{(W^{1,2})'}=\left(\sum_{i=0}^n(1+\alpha_i)^{-1}|\scpr{\hat v}{\theta_i}|^2\right)^{1/2}
\]
is a norm on {(the projection of)} $W^{1,2}(\Omega)'$, cf.~\cite[Prop.~3.4.8]{TucsWeis09}. Note that we can consider $P_n:W^{1,2}(\Omega)'\to W^{1,2}(\Omega)'$, which is a bounded linear operator with norm one, independent of $n$. Using this together with the fact that the injection from $L^2(\Omega)$ into $W^{1,2}(\Omega)'$ is continuous and $\cA\in\cL(W^{1,2}(\Omega),W^{1,2}(\Omega)')$, we can rewrite the weak formulation~\eqref{eq:weak_delta} as
\begin{equation}\label{eq:approx_dual}
\dot{v}_n=P_n\cA v_n+P_np_3(v_n)-P_nu_n+P_nI_{s,i}-P_n\cB k_0(\cB' v_n-\yref).
\end{equation}
Since $v_n\in L^2(0,\gamma;W^{1,2}(\Omega))$ and hence, by the Sobolev embedding theorem, $v_n\in L^2(0,\gamma;L^p(\Omega))$ for all $2\leq p\leq6$, we find that $p_3(v_n)\in L^2(0,\gamma;L^2(\Omega))$. We also have $\cA v_n\in L^2(0,\gamma;W^{1,2}(\Omega)')$ and $\cB k_0(\cB' v_n-\yref)\in L^2(0,\gamma;W^{1,2}(\Omega)')$ so that by using the {estimates~\eqref{eq:bound_1}--\eqref{eq:vn2} together with~\eqref{eq:approx_dual}}, there exists $C_4>0$ independent of $n$ and $t$ with
\[\|\dot{v}_n\|_{L^2(0,\gamma;W^{1,2}(\Omega)')}\leq C_4.\]
Now, by Lemma~\ref{lem:weak_convergence} we have that there exist subsequences of $(u_n)$, $(v_n)$ and $(\dot v_n)$, resp., again denoted in the same way, for which
\begin{equation}\label{eq:convergence_subseq}
\begin{aligned}
u_n\to u&\in W^{1,2}(0,\gamma;L^2(\Omega))\mbox{ weakly},\\
u_n\to u&\in W^{1,\infty}(0,\gamma;L^{2}(\Omega))\mbox{ weak}^\star,\\
v_n\to v&\in L^2(0,\gamma;W^{1,2}(\Omega))\mbox{ weakly},\\
v_n\to v&\in L^\infty(0,\gamma;L^{2}(\Omega))\mbox{ weak}^\star,\\
v_n\to v&\in L^4(0,\gamma;L^{4}(\Omega))\mbox{ weakly},\\
\dot{v}_n\to \dot{v}&\in L^2(0,\gamma;W^{1,2}(\Omega)')\mbox{ weakly}.
\end{aligned}
\end{equation}
Moreover, let
%\marginpar{DAS steht nicht Theorem \ref{thm:comp_inj}!}
$p_0=p_1=2$ and $X=W^{1,2}(\Omega)$, $Y=L^2(\Omega)$, $Z=W^{1,2}(\Omega)'$. Then, \cite[Chap.~1, Thm.~5.1]{Lion69} implies that
\[W:=\setdef{u\in L^{p_0}(0,\gamma;X) }{ \dot{u}\in L^{p_1}(0,\gamma;Z) }\]
with norm $\|u\|_{L^{p_0}(0,\gamma;X)}+\|\dot{u}\|_{L^{p_1}(0,\gamma;Z)}$ has a compact injection into $L^{p_0}(0,\gamma;Y)$, so that the weakly convergent sequence $v_n\to v\in W$ converges strongly in $L^2(0,\gamma;L^2(\Omega))$ by \cite[Lem.~1.6]{HinzPinn09}. Further, $(u(0),v(0))=(u_0,v_0)$ and by $v\in W^{1,2}(0,\gamma;L^2(\Omega))$, $v\in L^2(0,\gamma;W^{1,2}(\Omega))$ and $\dot{v}\in L^2(0,\gamma;W^{1,2}(\Omega)')$ it follows that $u,v\in C([0,\gamma];L^2(\Omega))$, see for instance \cite[Thm.~1.32]{HinzPinn09}. Moreover, note that $\cB' v-\yref\in L^2(0,\gamma;\R^m)$. Hence, $(u,v)$ is a solution of \eqref{eq:FHN_feedback} in $[0,\gamma]$ and
{
\begin{equation}\label{eq:strong_delta}
\begin{aligned}
\dot{v}(t)&=\cA v(t)+p_3(v(t))-u(t)+I_{s,i}(t)-\cB k_0(\cB' v(t)-\yref(t)) ,\ \ &v(0)&=v_0,\\
\dot{u}(t)&=c_5v(t)-c_4 u(t),&u(0)&=u_0,\\
\end{aligned}
\end{equation}}
is satisfied in $W^{1,2}(\Omega)'$. Moreover, by~\eqref{eq:vn2},~\cite[Chap.~1, Lem.~1.3]{Lion69} and $v_n\to v$ in $L^4(0,\gamma;L^{4}(\Omega))$ we have that $v_n^3\to v^3$ weakly in $L^{4/3}(0,\gamma;L^{4/3}(\Omega))$ and $v_n^2\to v^2$ weakly in $L^2(0,\gamma;L^{2}(\Omega))$.\\
\emph{Step 1e: We show uniqueness of the solution $(v,u)$.} To this end, we separate the linear part of~$p_3$ so that
\[p_3(v)=-c_1v-c_3\hat{p}_3(v),\quad \hat{p}_3(v)\coloneqq v^2\left(v-c\right),\ \quad c\coloneqq c_2/c_3.\]
Assume that $(v_1,u_1)$ and $(v_2,u_2)$ are two solutions of~\eqref{eq:FHN_feedback} on $[0,\gamma]$ with the same initial values, $v_1(0) = v_2(0) = v_0$ and $u_1(0) = u_2(0) = u_0$. Let $t_0\in(0,\gamma]$ be given. Let $Q_0\coloneqq(0,t_0)\times\Omega$. Define
\[\Sigma(t,\zeta):= |v_1(t,\zeta)|+|v_2(t,\zeta)|,\]
and let
\[Q^\Lambda:= \{(t,\zeta)\in Q_0\st\Sigma(t,\zeta)\leq\Lambda\},\quad \Lambda>0.\]
Note that, by convexity of the map $x\mapsto x^p$ on $[0,\infty)$ for $p>1$, we have that
\[
 \forall\, a,b\ge 0:\ \big(\tfrac12 a+ \tfrac12 b)^p \le \tfrac12 a^p + \tfrac12 b^p.
\]
Therefore, since $v_1,v_2\in L^4(0,\gamma;L^{4}(\Omega))$, we find that $\Sigma\in L^4(0,\gamma;L^{4}(\Omega))$. Hence, by the monotone convergence theorem, for all $\epsilon>0$ we may choose $\Lambda$ large enough such that
\[\int_{Q_0\setminus Q^\Lambda}|\Sigma(\zeta,t)|^4\,{\ds{\zeta}\ds{t}}<\epsilon.\]
Note that without loss of generality we may assume that $\Lambda>\frac{c}{3}$. Let $V:= v_2-v_1$ and $U:= u_2-u_1$, then, by~\eqref{eq:FHN_feedback},
\begin{align*}
  \dot V &= (\cA-c_1I) V -c_3(\hat{p}_3(v_2) - \hat{p}_3(v_1)) - U - k_0 \cB\cB' V,\\
  \dot U &= c_5 V - c_4 U.
\end{align*}
By~\cite[Thm.~1.32]{HinzPinn09}, we have for all $t\in(0,\gamma)$ that
{
\begin{align*}
 \tfrac{1}{2}\ddt\|V(t)\|^2&=\scpr{\dot{V}(t)}{V(t)}_{{W^{1,2}(\Omega)',W^{1,2}(\Omega)}},\\ \tfrac{1}{2}\ddt\|U(t)\|^2&=\scpr{\dot{U}(t)}{U(t)}_{{W^{1,2}(\Omega)',W^{1,2}(\Omega)}},
 \end{align*}
}
% \[\tfrac{1}{2}\ddt\|V(t)\|^2=\scpr{\dot{V}(t)}{V(t)},\quad \tfrac{1}{2}\ddt\|U(t)\|^2=\scpr{\dot{U}(t)}{U(t)},\]
thus we may compute that
\begin{align*}
    \tfrac{c_5}{2}\ddt\|V\|^2+\tfrac{1}{2}\ddt\|U\|^2&= \scpr{(\cA-c_1I) V - U - k_0 \cB\cB' V}{c_5 V}_{{W^{1,2}(\Omega)',W^{1,2}(\Omega)}} \\& \quad  - c_4 \|U\|^2 + c_5 \langle U,V \rangle -c_5 c_3\scpr{ \hat{p}_3(v_2) - \hat{p}_3(v_1)}{V}\\
    &= c_5 \scpr{(\cA-c_1I) V}{V}_{{W^{1,2}(\Omega)',W^{1,2}(\Omega)}} - c_5 k_0 \scpr{\cB' V}{\cB' V}_{\R^m}\\ &\quad - c_4 \|U\|^2
    - c_5c_3 \scpr{\hat{p}_3(v_2)-\hat{p}_3(v_1)}{V}\\
    &\le-c_5c_3 \scpr{\hat{p}_3(v_2)-\hat{p}_3(v_1)}{V}.
\end{align*}
Integration over $[0,t_0]$ and using $(U(0),V(0))=(0,0)$ leads to
\begin{equation}\label{eq:aux}
\begin{aligned}
\tfrac{c_5}{2}\|V(t_0)\|^2+\tfrac{1}{2}\|U(t_0)\|^2
&{\le}-c_5c_3\int_0^{t_0}\int_\Omega (\hat{p}_3(v_2(\zeta,t))-\hat{p}_3(v_1(\zeta,t)))V(\zeta,t)\ds{\zeta}\ds{t}\\
&=-c_5c_3\int_{Q^\Lambda} (\hat{p}_3(v_2(\zeta,t))-\hat{p}_3(v_1(\zeta,t)))V(\zeta,t)\ds{\zeta}\ds{t}\\
&\quad-c_5c_3\int_{Q_0\setminus Q^\Lambda} (\hat{p}_3(v_2(\zeta,t))-\hat{p}_3(v_1(\zeta,t)))V(\zeta,t)\ds{\zeta}\ds{t}.
\end{aligned}
\end{equation}

Note that on $Q^\Lambda$ we have $-\Lambda\leq v_1\leq\Lambda$ and $-\Lambda\leq v_2\leq\Lambda$. Let $a,b\in[-\Lambda,\Lambda]$, then the mean value theorem implies
\begin{align*}
(\hat{p}_3(b)-\hat{p}_3(a))(b-a)=\hat{p}_3'(\xi)(b-a)^2
\end{align*}
for some $\xi\in(-\Lambda,\Lambda)$. Since $\hat{p}_3'(\xi)=3\xi^2-2c\xi$ has a minimum at
\[\xi^\ast=\frac{c}{3}\]
we have that
\[(\hat{p}_3(b)-\hat{p}_3(a))(b-a)=\hat{p}_3'(\xi)(b-a)^2\geq-\frac{c^2}{3}(b-a)^2.\]
Using that {in inequality \eqref{eq:aux}} leads to
\begin{align*}
\tfrac{c_5}{2}\|V(t_0)\|^2+\tfrac{1}{2}\|U(t_0)\|^2&\leq c_5c_3\frac{c^2}{3}\int_{Q^\Lambda} V(\zeta,t)^2\ds{\zeta}\ds{t}\\
&\quad-c_5c_3\int_{Q_0\setminus Q^\Lambda} (\hat{p}_3(v_2(\zeta,t))-\hat{p}_3(v_1(\zeta,t)))V(\zeta,t)\ds{\zeta}\ds{t}\\
&\leq c_5c_3\frac{c^2}{3}\int_{Q_0} V(\zeta,t)^2\ds{\zeta}\ds{t}\\
&\quad+c_5c_3\int_{Q_0\setminus Q^\Lambda} |\hat{p}_3(v_2(\zeta,t))-\hat{p}_3(v_1(\zeta,t))||V(\zeta,t)|\ds{\zeta}\ds{t}\\
&\le c_5c_3\frac{c^2}{3}\int_0^{t_0}\|V(t)\|^2\ds{t}+2c_5c_3\int_{Q_0\setminus Q^\Lambda}|\Sigma(\zeta,t)|^4 {\ds{\zeta}\ds{t}} \\
&\leq  c_3\frac{c^2}{3}\int_0^{t_0}c_5\|V(t)\|^2+\|U(t)\|^2\ds{t}+2c_5c_3\epsilon.
\end{align*}
Since $\epsilon>0$ was arbitrary we may infer that
\[\tfrac{c_5}{2}\|V(t_0)\|^2+\tfrac{1}{2}\|U(t_0)\|^2\leq \frac{2c_3c^2}{3}\int_0^{t_0}\tfrac{c_5}{2}\|V(t)\|^2+\tfrac12\|U(t)\|^2\ds{t}.\]
Hence, by Gronwall's lemma and $U(0)=0,V(0)=0$ it follows that $U(t_0)=0$ and $V(t_0)=0$. Since $t_0$ was arbitrary, this shows that $v_1 = v_2$ and $u_1 = u_2$ on $[0,\gamma]$.

\emph{Step 2: We show that for all $\epsilon\in(0,\gamma)$ and all $t\in[\epsilon,\gamma]$ we have $v(t)\in W^{1,2}(\Omega)$. {In particular, this guarantees $v(\gamma)\in W^{1,2}(\Omega)$, which is required as an initial condition in the second part of the proof in Section~\ref{ssec:mono_proof_tgeqgamma}.}}\\
Fix $\epsilon\in(0,\gamma)$. First we show that $v\in BUC([\epsilon,\gamma];W^{1,2}(\Omega))$. Multiplying the first equation in~\eqref{eq:weak_delta}  by $\dot{\mu}_j$ and {summing over} $j\in\{0,\dots,n\}$ we obtain
\begin{align*}
\|\dot{v}_n\|^2&=-\tfrac{1}{2}\ddt \fa(v_n,v_n)-\scpr{u_n}{\dot{v}_n}+\scpr{p_3(v_n)}{\dot{v}_n}+\scpr{I_{s,i}}{\dot{v}_n}\\
&\quad -k_0\scpr{\cB' v_n-\yref}{\cB'\dot{v}_n}_{\R^m}\\
&=-\tfrac{1}{2}\ddt \fa(v_n,v_n)-\scpr{u_n}{\dot{v}_n}+\scpr{p_3(v_n)}{\dot{v}_n}+\scpr{I_{s,i}}{\dot{v}_n}\\
&\quad -k_0\scpr{\cB' v_n-\yref}{\cB'\dot{v}_n-\dyref}_{\R^m} - k_0\scpr{\cB' v_n-\yref}{\dyref}_{\R^m}
\end{align*}
Furthermore, we may derive that
\begin{align*}
  \ddt v_n^4 &= 4 v_n^3 \dot v_n = -\frac{4}{c_3}\big( p_3(v_n) - c_2 v_n^2 + c_1 v_n\big) \dot v_n,\quad \text{thus}\\
  p_3(v_n) \dot v_n &= -\frac{c_3}{4} \ddt v_n^4 + c_2 v_n^2 \dot v_n - c_1 v_n \dot v_n,
\end{align*}
and this implies, for any {$\eta>0$},
\begin{align*}
    \scpr{p_3(v_n)}{\dot v_n} &\le -\frac{c_3}{4} \ddt \|v_n\|_{L^4}^4 + c_2  \scpr{v_n^2}{\dot v_n} - c_1 \scpr{v_n}{\dot v_n}\\
    &\le -\frac{c_3}{4} \ddt \|v_n\|_{L^4}^4 + \frac{c_2}{2} \left( \eta \|v_n\|_{L^4}^4 + \frac{1}{\eta} \|\dot v_n\|^2\right) + \frac{c_1}{2} \left(\eta \|v_n\|^2 + \frac{1}{\eta} \|\dot v_n\|^2\right)\\
    &\stackrel{\eqref{eq:uv_bound_delta}}{\le} -\frac{c_3}{4} \ddt \|v_n\|_{L^4}^4 + \frac{c_2}{2} \left( \eta \|v_n\|_{L^4}^4 + \frac{1}{\eta} \|\dot v_n\|^2 \right) + \frac{c_1}{2} \left(\eta C_v + \frac{1}{\eta} \|\dot v_n\|^2\right).
\end{align*}
Moreover, we find that, recalling $\overline{e}_n = \cB' v_n-\yref$
\begin{align*}
\scpr{u_n}{\dot{v}_n} &\le \frac{\eta}{2} \|u_n\|^2 + \frac{1}{\eta} \|\dot v_n\|^2 \stackrel{\eqref{eq:uv_bound_delta}}{\le} \frac{\eta C_u}{2} + \frac{1}{\eta} \|\dot v_n\|^2,\\
\scpr{I_{s,i}}{\dot{v}_n} &\le \frac{\eta}{2} \|I_{s,i}\|_{2,\infty}^2 + \frac{1}{\eta} \|\dot v_n\|^2,\\
\scpr{\overline{e}_n}{\dyref}_{\R^m} &\le \frac{1}{2} \|\overline{e}_n\|^2_{\R^m} + \frac12 \|\dyref\|_\infty^2.
\end{align*}
Therefore, choosing $\eta$ large enough, we obtain that there exist constants $Q_1,Q_2>0$ independent of $n$ such that
\begin{align*}
\|\dot{v}_n\|^2\le&-\frac{1}{2}\ddt \fa(v_n,v_n)-\frac{c_3}{4}\ddt\|v_n\|_{L^4}^4-\frac{k_0}{2}\ddt\|\overline{e}_n\|_{\R^m}^2+\frac{1}{2}\|\dot{v}_n\|^2\\
&+Q_1\|v_n\|^4_{L^4}+Q_2+\frac{k_0}{2}\|\overline{e}_n\|_{\R^m}^2,
\end{align*}
thus,
\begin{equation}\label{eq:est-norm-dotvn}
\begin{aligned}
\|\dot{v}_n\|^2&+\ddt\left(\fa(v_n,v_n)+\frac{c_3}{2}\|v_n\|_{L^4}^4+k_0\|\overline{e}_n\|_{\R^m}^2\right)\\
\le &\ 2Q_1\|v_n\|^4_{L^4}+2Q_2+k_0\|\overline{e}_n\|_{\R^m}^2.
\end{aligned}
\end{equation}
As a consequence, we find that for all $t\in[0,\gamma]$ we have
\begin{align*}
t\|\dot{v}_n(t)\|^2&+\ddt\left(t\fa(v_n(t),v_n(t))+\frac{c_3t}{2}\|v_n(t)\|_{L^4}^4+k_0t\|\overline{e}_n(t)\|_{\R^m}^2\right)\\
\stackrel{\eqref{eq:est-norm-dotvn}}{\le}&\ \left(2Q_1 t +\frac{c_3}{2}\right)\|v_n(t)\|^4_{L^4}+\fa(v_n(t),v_n(t)) +2Q_2 t+k_0(t+1)\|\overline{e}_n(t)\|_{\R^m}^2.
\end{align*}
Since $t\|\dot{v}_n(t)\|^2\geq0$ and $t\le \gamma$ for all $t\in[0,\gamma]$, it follows that
\begin{align*}
\ddt&\left(t\fa(v_n(t),v_n(t))+\frac{c_3t}{2}\|v_n(t)\|_{L^4}^4+k_0t\|\overline{e}_n(t)\|_{\R^m}^2\right)\\
\le&\ \left(2Q_1\gamma+\frac{c_3}{2}\right)\|v_n(t)\|^4_{L^4}+\fa(v_n(t),v_n(t)) +2Q_2\gamma+k_0(\gamma+1)\|\overline{e}_n(t)\|_{\R^m}^2.
\end{align*}
Integrating the former and using \eqref{eq:bound_1}, there exist $P_1,P_2>0$ independent of $n$ such that for $t\in[0,\gamma]$ we have
\begin{align*}
t\fa(v_n(t),v_n(t))+\frac{c_3t}{2}\|v_n(t)\|_{L^4}^4+k_0t\|\overline{e}_n(t)\|_{\R^m}^2
\le P_1+P_2t.
\end{align*}
Thus, there exist constants $C_5,C_6>0$ independent of $n$ such that
\[\forall\, t\in[0,\gamma]:\ t\fa(v_n(t),v_n(t))\leq C_5\ \wedge\ t\|\overline{e}_n(t)\|_{\R^m}\leq C_6.\]
Hence, for all $\epsilon\in(0,\gamma)$, it follows from the above estimates together with~\eqref{eq:bound_1} that $v_n\in L^\infty(\epsilon,\gamma;W^{1,2}(\Omega))$ and $\overline{e}_n\in L^\infty(\epsilon,\gamma;\R^m)$, so that in addition to~\eqref{eq:convergence_subseq}, from Lemma~\ref{lem:weak_convergence} we further have that there exists a subsequence such that
\[v_n\to v\in L^\infty(\epsilon,\gamma;W^{1,2}(\Omega))\mbox{ weak}^\star\]
and $\cB'v\in L^\infty(\epsilon,\gamma;\R^m)$ for all $\epsilon\in(0,\gamma)$, hence $I_{s,e}\in L^2(0,\gamma;\R^m)\cap L^\infty(\epsilon,\gamma;\R^m)$. By the Sobolev embedding theorem, $W^{1,2}(\Omega)\hookrightarrow L^p(\Omega)$ for $2\leq p\leq 6$ we have that $p_3(v)\in L^\infty(\epsilon,\gamma;L^2(\Omega))$. Moreover, since \eqref{eq:strong_delta} holds, we can rewrite it as
\[\dot{v}(t)=(\cA-c_1 I) v(t)+I_r(t)+\cB I_{s,e}(t),\]
where $I_r:=c_2v^2-c_3v^3-u+I_{s,i}\in L^2(0,\gamma;L^2(\Omega))\cap L^\infty(\epsilon,\gamma;L^2(\Omega))$ and Proposition~\ref{prop:hoelder} (recall that $W^{1,2}(\Omega)' = X_{-1/2}$ and hence $\cB\in\cL(\R^m,X_{-1/2})$) with $c=c_1$ implies that $v\in BUC([\epsilon,\gamma];W^{1,2}(\Omega))$. Hence, for all $\epsilon\in(0,\gamma)$, $v(t)\in W^{1,2}(\Omega)$ for $t\in[\epsilon,\gamma]$, so that in particular $v(\gamma)\in W^{1,2}(\Omega)$.

\subsection{Solution on $(\gamma,\infty)$}
\label{ssec:mono_proof_tgeqgamma}

The crucial step in this part of the proof is to show that the error remains uniformly bounded away from the funnel boundary while $v\in L^\infty(\gamma,\infty;W^{1,2}(\Omega))$. The proof is divided into several steps.

\emph{Step 1: We show existence of an approximate solution by means of a time-varying state-space transformation.}\\
Again, let $(\theta_i)_{i\in\N_0}$ be the eigenfunctions of $-\cA$ and let $\alpha_i$ be the corresponding eigenvalues, with $\alpha_i\geq0$ for all $i\in\N_0$. Recall that $(\theta_i)_{i\in\N_0}$ {forms} an orthonormal basis of $L^2(\Omega)$ by Remark~\ref{Rem:Aop_n}\,\ref{item:Aop4}). Let $(u_\gamma,v_\gamma):=(u(\gamma),v(\gamma))$, $a_i :=  \scpr{v_\gamma}{\theta_i}$ and $b_i :=  \scpr{u_\gamma}{\theta_i}$ for $i\in\N_0$ and
\[
v_\gamma^n:= \sum_{i=0}^na_{i}\theta_i,\quad u_\gamma^n:= \sum_{i=0}^nb_{i}\theta_i,\quad n\in\N.
\]
Then we have that $v^n_\gamma\to v_\gamma$ strongly in $W^{1,2}(\Omega)$ and $u^n_\gamma\to u_\gamma$ strongly in $L^2(\Omega)$. As stated in Remark~\ref{rem:main}\,\ref{rem:main1}) we have that $\ker\cB=\{0\}$ implies $\cB' \cD(\cA)=\R^m$. As a~consequence, there exist $q_1,\ldots,q_m\in\cD(\cA)$ such that $\cB'q_k=e_k$, {where $e_k$ denotes the $k$-th unit vector in $\R^m$ for $k=1,\dots,m$.}  By Remark~\ref{Rem:Aop_n}\,\ref{item:Aop3}), we further have  $q_k\in C^{0,\nu}(\Omega)$ for some $\nu\in(0,1)$.\\
Note that $U\coloneqq\bigcup_{n\in\N} U_n$, where $U_n = \mathrm{span}\{\theta_i\}_{i=0}^n$, satisfies $\overline{U}=W^{1,2}(\Omega)$ with the respective norm. Moreover, $\overline{\cB' U}=\R^m$. Since $\R^m$ is complete and finite dimensional and $\cB'$ is linear and continuous it follows that $\cB' U=\R^m$. By the surjectivity of $\cB'$ we have that for all $k\in\{1,\dots,m\}$ there exist $n_k\in\N$ and $q_k\in U_{n_k}$ such that $\cB' q_k=e_k$. Thus, there exists $n_0\in\N$ with $q_k\in U_{n_0}$ for all $k=\{1,\dots,m\}$, hence the $q_k$ are a (finite) linear combination of the eigenfunctions $\theta_i$.\\
Define $q\in W^{1,2}(\Omega;\R^m)\cap C^{0,\nu} (\Omega;\R^m)$ by $q(\zeta)=\big(q_1(\zeta),\ldots,q_m(\zeta)\big)^\top$ and $q\cdot\yref$ by
\[
(q\cdot\yref) (t,\zeta) :=  \sum_{k=1}^mq_k(\zeta)\yrefk(t),\quad \zeta\in \Omega,\, t\ge 0.
\]
We may define $q\cdot\dyref$ analogously. Note that we have $(q\cdot\yref)\in BC([0,\infty)\times\Omega)$, because
\[
|(q\cdot\yref) (t,\zeta)| \le \sum_{k=1}^m \|q_k\|_\infty\, \|y_{{\rm ref},k}\|_\infty
\]
for all $\zeta\in \Omega$ and $t\ge 0$, where we write $\|\cdot\|_\infty$ for the supremum norm. We define $q_{k,j}:= \scpr{q_k}{\theta_j}$ for $k=1,\ldots,m$, $j\in\N_0$ and $q_k^n:=\sum_{j=0}^nq_{k,j}$ for $n\in\N_0$. Similarly, $q^n:=(q_1^n,\dots,q_m^n)^\top$, $n\in\N$, {satisfies $q^n = q$ for all $n\ge n_0$ since $q_k\in U_{n_0}$ for all $k=1,\ldots,m$, thus} $q^n\to q$ strongly in $W^{1,2}(\Omega)$.\\
Since $\cB':W^{r,2}(\Omega)\to\R^m$ is continuous for some $r\in[0,1]$, it follows that for all $\theta\in W^{r,2}(\Omega)$ there exists $\Gamma_r>0$ such that
\[\|\cB'\theta\|_{\R^m}\leq\Gamma_r\|\theta\|_{W^{r,2}}.\]
For $n\in\N_0$, let
\[\kappa_n\coloneqq \big((n+1) \Gamma_r (1+\|v_\gamma^n-q^n\cdot\yref(\gamma)\|^2_{W^{r,2}})\big)^{-1}.\]
Note that for $v_\gamma\in W^{1,2}(\Omega)$ it holds that $\kappa_n>0$ for all $n\in\N_0$, $(\kappa_n)_{n\in\N_0}$ is bounded by $\Gamma_r^{-1}$ (and monotonically decreasing) and $\kappa_n\to0$ as $n\to\infty$ and by construction
\[\forall\,n\in\N_0:\ \kappa_n\|\cB' (v_\gamma^n-q^n\cdot\yref(\gamma))\|_{\R^m}<1.\]
Consider a modification of $\varphi$ induced by $\kappa_n$, namely
\[\varphi_n\coloneqq\varphi+\kappa_n,\quad n\in\N_0.\]
It is clear that for each $n\in\N_0$ we have $\varphi_n\in W^{1,\infty}([\gamma,\infty);\R)$, the estimates $\|\varphi_n\|_\infty\leq\|\varphi\|_\infty+\Gamma_r^{-1}$ and $\|\dot{\varphi}_n\|_\infty=\|\dot{\varphi}\|_\infty$ are independent of $n$, and
$\varphi_n\to\varphi\in\Phi_\gamma$ uniformly. Moreover, $\inf_{t>\gamma}\varphi_n(t)>0$.\\
Now, fix $n\in\N_0$. For $t\ge \gamma$, define, {in the respective spaces},
\begin{align*}
\phi(e)&:= \frac{k_0}{1-\|e\|_{\R^m}^2}e,\quad e\in\R^m,\ \|e\|_{\R^m}<1,\\
\omega_0(t)&:= \dot{\varphi}_n(t)\varphi_n(t)^{-1},\\
F(t,z)&:= \varphi_n(t)f_{-1}(t)+\varphi_n(t)f_0(t)+f_1(t)z+\varphi_n(t)^{-1}f_2(t)z^2\\
&\quad \ -c_3\varphi_n(t)^{-2}z^3,\quad z\in\R,\\
f_{-1}(t)&:= I_{s,i}(t)+\sum_{k=1}^m\yrefk(t)\cA q_k,\\
f_0(t)&:= -q\cdot(\dyref(t)+c_1\yref(t))+c_2(q\cdot\yref(t))^2-c_3(q\cdot\yref(t))^3,\\
f_1(t)&:= (q\cdot\yref(t))(2c_2-3c_3(q\cdot\yref(t))),\\
f_2(t)&:= c_2-3c_3(q\cdot\yref(t)),\\
g(t)&:= c_5(q\cdot\yref(t)).
\end{align*}
We have that $f_{-1}\in L^\infty(\gamma,\infty;L^2(\Omega))$, since
\begin{align*}
\|f_{-1}\|_{2,\infty}  &:=  \esssup_{t\ge \gamma} \left(\int_\Omega f_{-1}(\zeta,t)^2 \ds{\lambda} \right)^{1/2}\\
&\le \|I_{s,i}\|_{2,\infty} + \sum_{k=1}^m \|\yrefk\|_\infty\, \|A q_k\|_{L^2} < \infty.
\end{align*}
Furthermore, we have that $f_0\in L^\infty((\gamma,\infty)\times\Omega)$, because
\begin{align*}
|f_0(\zeta,t)|&\leq\ (\|\dyref\|_\infty+c_1\|\yref\|_\infty)\sum_{k=1}^m\|q_k\|_\infty+c_2\|\yref\|_\infty^2\left(\sum_{k=1}^m\|q_k\|_\infty\right)^2\\
&\quad +c_3\|\yref\|^3_\infty\left(\sum_{k=1}^m\|q_k\|_\infty\right)^3\text{ for a.a.\ $(\zeta,t)\in \Omega\times [\gamma,\infty)$},
\end{align*}
whence
\[\|f_0\|_{\infty,\infty}:= \esssup_{t\geq\gamma,\zeta\in\Omega}|f_0(\zeta,t)|<\infty.\]
Similarly $\|f_1\|_{\infty,\infty}<\infty$, $\|f_2\|_{\infty,\infty}<\infty$ and $\|g\|_{\infty,\infty}<\infty$.\\
Consider the system of $2(n+1)$ ODEs
\begin{equation}\label{eq:appr_ODE}
\begin{aligned}
\dot{\mu}_j(t)&=-\alpha_j\mu_j(t)-(c_1-\omega_0(t))\mu_j(t)-\nu_j(t)- \scpr{\phi\left(\sum_{i=0}^n \cB'\theta_i \mu_i(t)\right)}{\cB'\theta_j}_{\R^m} \\
&\quad +\scpr{F\left(t,\sum_{i=0}^n\mu_i(t)\theta_i\right)}{\theta_j},\\
\dot{\nu}_j(t)&=-(c_4-\omega_0(t))\nu_j(t)+c_5\mu_j(t)+\varphi_n(t)\scpr{g(t)}{\theta_j}
\end{aligned}
\end{equation}
defined on
\[
\D:= \setdef{(t,\mu_0,\dots,\mu_n,\nu_0,\dots,\nu_n)\in[\gamma,\infty)\times\R^{2(n+1)} }{ \left\|\sum_{i=0}^n\gamma_i\mu_i\right\|_{\R^m}<1 },
\]
with initial value
\[\mu_j(\gamma)=\kappa_n\left(a_j-\sum_{k=1}^mq_{k,j}\yrefk(\gamma)\right),\quad \nu_j(\gamma)=\kappa_n b_j,\quad j\in\N_0.\]
{Given that} the functions {defining the system of ODEs~\eqref{eq:appr_ODE}} are continuous, the set~$\D$ is relatively open in $[\gamma,\infty)\times\R^{2(n+1)}$ and by construction the initial condition satisfies $(\gamma,\mu_0(\gamma),\dots,\mu_n(\gamma),\nu_0(\gamma),\dots,\nu_n(\gamma))\in\D$ it follows from ODE theory, see e.g.~\cite[\S~10, Thm.~XX]{Walt98}, that there exists a weakly differentiable solution
\[(\mu^n,\nu^n)=(\mu_0,\dots,\mu_n,\nu_0,\dots,\nu_n):[\gamma,T_n)\to\R^{2(n+1)}\] such that $T_n\in(\gamma,\infty]$ is maximal. Furthermore, the closure of the graph of~$(\mu^n,\nu^n)$ is not a compact subset of $\D$.\\
With that, we may define
\[z_n(t):=\sum_{i=0}^n\mu_i(t)\theta_i,\quad w_n(t):=\sum_{i=0}^n\nu_i(t)\theta_i,\quad e_n(t):= \sum_{i=0}^n\cB'\theta_i\mu_i(t),\quad t\in[\gamma,T_n)\]
and note that
\[z_\gamma^n:=z_n(\gamma)=\kappa_n(v_\gamma^n-q^n\cdot\yref(\gamma)),\quad w_\gamma^n:=w_n(\gamma)=\kappa_n u_\gamma^n.\]
From the orthonormality of the $\theta_i$ we have that
\begin{equation}\label{eq:weak}
	\begin{aligned} \scpr{\dot{z}_n(t)}{\theta_j}&=-\fa(z_n(t),\theta_j)-(c_1-\omega_0(t))\scpr{z_n(t)}{\theta_j} - \scpr{w_n(t)}{\theta_j}\\
	&\quad -\scpr{\phi\left(\cB' z_n(t)\right)}{\cB' \theta_j}_{\R^m} +\scpr{F\left(t,z_n(t)\right)}{\theta_j},\\ \scpr{\dot{w}_n(t)}{\theta_j} &= -(c_4-\omega_0(t))\scpr{w_n(t)}{\theta_j}+c_5\scpr{z_n(t)}{\theta_j}+\varphi_n\scpr{g(t)}{\theta_j}.
\end{aligned}
	\end{equation}
Define now
\begin{equation}\label{eq:transformation}
\begin{aligned}
v_n(t)&\coloneqq \varphi_n(t)^{-1}z_n(t)+q^n\cdot\yref(t),\\
u_n(t)&\coloneqq \varphi_n(t)^{-1}w_n(t),\\
\tilde{\mu}_i(t) &\coloneqq \varphi_n(t)^{-1}\mu_i(t) +\sum_{k=1}^m q_{k,i}\yrefk(t),\\
\tilde{\nu}_i(t) &\coloneqq \varphi_n(t)^{-1}\nu_i(t),
\end{aligned}
\end{equation}
then $v_n(t)=\sum_{i=0}^n\tilde{\mu}_i(t)\theta_i$ and $u_n(t)=\sum_{i=0}^n\tilde{\nu}_i(t)\theta_i$. With this transformation we obtain that $(v_n, u_n)$ satisfies, for all $\theta\in W^{1,2}(\Omega)$, $\chi\in L^2(\Omega)$ and all $t\in[\gamma,T_n)$ that
\begin{equation*}
\begin{aligned}
\ddt\scpr{v_n(t)}{\theta}=&-\fa(v_n(t),\theta)+\scpr{p_3(v_n(t)+(q-q^n)\cdot\yref(t))-u_n(t)}{\theta}\\
&+\scpr{I_{s,i}(t)-(q-q^n)\cdot\dyref(t)+\sum_{k=1}^m\yrefk(t)\cA(q_k-q_k^n) }{\theta}\\
&+\scpr{I_{s,e}^n(t)}{\cB'\theta}_{\R^m},\\
\ddt\scpr{u_n(t)}{\chi}=&\scpr{c_5(v_n(t)+(q-q^n)\cdot\yref(t))-c_4u_n(t)}{\chi},\\
I_{s,e}^n(t)=&-\frac{k_0}{1-\varphi_n(t)^2\|\cB' (v_n(t)-q^n\cdot\yref(t))\|^2_{\R^m}}(\cB'( v_n(t)-q^n\cdot\yref(t))),
\end{aligned}
\end{equation*}
with $(u_n(\gamma),v_n(\gamma))=(u_\gamma,v_\gamma)$. Since there exists some $n_0\in\N$ with $q^{n}=q$ for all $n\geq n_0$, we have for all $n\geq n_0$, $\theta\in W^{1,2}(\Omega)$ and $\chi\in L^2(\Omega)$ that
\begin{equation}\label{eq:weak_uv}
\begin{aligned}
\ddt\scpr{v_n(t)}{\theta}=&-\fa(v_n(t),\theta)+\scpr{p_3(v_n(t))-u_n(t)}{\theta}\\
&+\scpr{I_{s,i}(t)}{\theta}+\scpr{I_{s,e}^n(t)}{\cB'\theta}_{\R^m},\\
\ddt\scpr{u_n(t)}{\chi}=&\scpr{c_5v_n(t)-c_4u_n(t)}{\chi},\\
I_{s,e}^n(t)=&-\frac{k_0}{1-\varphi_n(t)^2\|\cB' v_n(t)-\yref(t)\|^2_{\R^m}}(\cB' v_n(t)-\yref(t)),
\end{aligned}
\end{equation}
%Note that if the $q_k$ are a (finite) linear combination of eigenfunctions $\theta_i$, then for some $n_0\in\N$, $q=q^n$ for all $n\geq n_0$.
%\pagebreak[2]

\emph{Step 2: We show boundedness of $(z_n,w_n)$ in terms of~$\varphi_n$.}\\
Consider again the Lyapunov function \eqref{eq:Lyapunov} and observe that $\|z_n(t)\|^2 = \sum_{j=0}^n \mu_j(t)^2$ and $\|w_n(t)\|^2 = \sum_{j=0}^n \nu_j(t)^2$. We find that, for all $t\in[\gamma, T_n)$,
\begin{align*}
\ddt V(z_n(t),w_n(t)) &= c_5\sum_{j=0}^n\mu_j(t)\dot{\mu}_j(t)+\sum_{j=0}^n\nu_j(t)\dot{\nu}_j(t)\\ &=-c_5\sum_{j=0}^{n}\alpha_j\mu_j(t)^2- c_5 (c_1-\omega_0(t))\sum_{j=0}^{n}\mu_j(t)^2\\
&\quad-(c_4-\omega_0(t))\sum_{j=0}^{n}\nu_j(t)^2 -c_5\scpr{\phi(e_n(t))}{e_n(t)}_{\R^m} \\
&\quad +\varphi_n (t) \scpr{g(t)}{\sum_{i=0}^n\nu_i(t)\theta_i}\\
&\quad+c_5\scpr{F\left(t,\sum_{i=0}^n\mu_i(t)\theta_i\right)}{\sum_{i=0}^n\mu_i(t)\theta_i},
\end{align*}
hence, omitting the argument~$t$ for brevity in the following,
\begin{equation}\label{eq:Lyapunov_boundary_1}
\begin{aligned}
\ddt V(z_n,w_n)=&-c_5\fa(z_n,z_n)-c_5(c_1-\omega_0)\|z_n\|^2-(c_4-\omega_0)\|w_n\|^2\\
&-c_5\frac{k_0\|e_n\|_{\R^m}^2}{1-\|e_n\|_{\R^m}^2}+c_5\scpr{F(t,z_n)}{z_n}+\varphi_n \scpr{g}{w_n}.
\end{aligned}
\end{equation}
Next we use some Young and Hölder inequalities to estimate the term
\begin{align*} \scpr{F(t,z_n)}{z_n}&=\underbrace{\varphi_n(t)\scpr{f_{-1}(t)}{z_n}}_{I_{-1}} +\underbrace{\varphi_n(t)\scpr{f_0(t)}{z_n}}_{I_0}+\underbrace{\scpr{f_1(t)z_n}{z_n}}_{I_1}\\	&\quad+\underbrace{\varphi_n(t)^{-1}\scpr{f_2(t)z_n^2}{z_n}}_{I_2}-c_3\varphi_n(t)^{-2}\underbrace{\scpr{z_n^3}{z_n}}_{=\|z_n\|_{L^4}^4}.
\end{align*}
For the first term we derive, using Young's inequality for products with $p=4/3$ and $q=4$, that
\begin{align*}
I_{-1}&\leq \scpr{\frac{2^{1/2}\varphi_n^{3/2} |I_{s,i}|}{c_3^{1/4}}}{\frac{c_3^{1/4}|z_n|}{2^{1/2}\varphi_n^{1/2}}}+ \sum_{k=1}^m\scpr{\frac{(4m)^{1/4}\varphi_n^{3/2}\|\yref\|_\infty |Aq_k|}{c_3^{1/4}}}{\frac{c_3^{1/4}|z_n|}{(4m)^{1/4}\varphi_n^{1/2}}}\\
&\leq\frac{2^{2/3} 3\varphi_n^2\|I_{s,i}\|_{2,\infty}^{4/3}|\Omega|^{1/3}}{4c_3^{1/3}}+ \sum_{k=1}^m\frac{3(4m)^{1/3}\varphi_n^2\|\yref\|_\infty^{4/3}\|Aq_k\|^{4/3}|\Omega|^{1/3}}{4c_3^{1/3}}+ \frac{c_3\|z_n\|^4_{L^4}}{8\varphi_n^2}
\end{align*}
and with the same choice we obtain for the second term
\[I_0\leq\scpr{\frac{2^{1/4}\varphi_n^{3/2}\|f_0\|_{\infty,\infty}}{c_3^{1/4}}}{\frac{c_3^{1/4}|z_n|}{2^{1/4}\varphi_n^{1/2}}}\leq \frac{2^{1/3}3\varphi_n^2\|f_0\|_{\infty,\infty}^{4/3}|\Omega|}{4c_3^{1/3}}+\frac{c_3\|z_n\|^4_{L^4}}{8\varphi_n^2}.\]
Using $p=q=2$ we find that the third term satisfies
\[ I_1\leq\scpr{\frac{2\varphi_n\|f_1\|_{\infty,\infty}}{\sqrt{c_3}}}{\frac{\sqrt{c_3}|z_n|^2}{2\varphi_n}}\leq\frac{2\varphi_n^2\|f_1\|_{\infty,\infty}^2|\Omega|}{c_3}+ \frac{c_3\|z_n\|^4_{L^4}}{8\varphi_n^2},
\]
and finally, with $p=4$ and $q=4/3$,
\begin{align*} I_2&\leq\scpr{\varphi_n^{-1}\|f_2\|_{\infty,\infty}}{|z_n|^3}= \scpr{\frac{3^{3/2}\varphi_n^{1/2}\|f_2\|_{\infty,\infty}}{c_3^{3/4}}}{\left|\frac{c_3^{1/4}z_n}{\varphi_n^{1/2}\sqrt{3}}\right|^3}\\
&\leq\frac{9^3\varphi_n^2\|f_2\|^4_{\infty,\infty}|\Omega|}{4c_3^3}+\frac{c_3}{12\varphi_n^2}\|z_n\|^4_{L^4}.
\end{align*}
Summarizing, we have shown that
\[\scpr{F(t,z_n)}{z_n}\leq K_0\varphi_n^2-\frac{13c_3}{24\varphi_n^2}\|z_n\|^4_{L^4}\leq K_0\varphi_n^2-\frac{c_3}{2\varphi_n^2}\|z_n\|^4_{L^4},\]
where
\begin{align*}
K_0\coloneqq  & \ \frac{2^{2/3} 3\|I_{s,i}\|_{2,\infty}^{4/3}|\Omega|^{1/3}}{4c_3^{1/3}}+ \sum_{k=1}^m\frac{3(4m)^{1/3}\|\yref\|_\infty^{4/3}\|Aq_k\|^{4/3}|\Omega|^{1/3}}{4c_3^{1/3}}\\
&+ \frac{2^{1/3}3\|f_0\|_{\infty,\infty}^{4/3}|\Omega|}{4c_3^{1/3}}+\frac{2\|f_1\|_{\infty,\infty}^2|\Omega|}{c_3}+\frac{9^3\|f_2\|^4_{\infty,\infty}|\Omega|}{4c_3^3}.
\end{align*}
Finally, using Young's inequality with $p=q=2$, we estimate the last term in~\eqref{eq:Lyapunov_boundary_1} as follows
\[\varphi_n \scpr{g}{w_n}\leq\frac{\varphi_n^2\|g\|_{\infty,\infty}^2|\Omega|}{2c_4}+\frac{c_4}{2}\|w_n\|^2.\]
We have thus obtained the estimate
\begin{equation}\label{eq:Lyapunov_boundary_2}
\begin{aligned}
\ddt V(z_n,w_n)\leq&-(\sigma-2\omega_0) V(z_n,w_n)\\
&-c_5\fa(z_n,z_n)-c_5\frac{k_0\|e_n\|_{\R^m}^2}{1-\|e_n\|_{\R^m}^2}-\frac{c_3c_5}{2\varphi_n^{2}}\|z_n\|_{L^4}^4+\varphi_n^2K_1,
\end{aligned}
\end{equation}
where
\begin{align*}
\sigma\coloneqq  2\min\{c_1,c_4\},\quad  K_1\coloneqq  c_5K_0+\frac{\|g\|_{\infty,\infty}^2|\Omega|}{2c_4}.
\end{align*}
In particular, we have the
% conservative
estimate
\[
\ddt V(z_n,w_n) \le -(\sigma-2\omega_0) V(z_n,w_n) +\varphi_n^2K_1
\]
on $[\gamma,T_n)$, which implies that
\[
V(z_n(t),w_n(t)) \le \ee^{-K(t,\gamma)} V(z_n(\gamma),w_n(\gamma)) + \int_\gamma^t \ee^{-K(t,s)} \varphi_n(s)^2 K_1 \ds{s},
\]
where
\[
K(t,s) = \int_s^t \sigma -2\omega_0(\tau) \ds{\tau} = \sigma(t-s) - 2\ln \varphi_n(t) + 2 \ln \varphi_n(s),\quad \gamma\le s\le t < T_n.
\]
Therefore, invoking $\varphi_n(\gamma)=\kappa_n$, for all $t\in[\gamma, T_n)$ we have
\begin{align*}
&c_5\|z_n(t)\|^2+\|w_n(t)\|^2 = 2V(z_n(t),w_n(t))\\
&\leq 2\ee^{-\sigma (t-\gamma)}\frac{\varphi_n(t)^2}{\kappa_n^2}V(z_n(\gamma),w_n(\gamma))+\frac{2K_1}{\sigma}\varphi_n(t)^2\\
&= \varphi_n(t)^2\left((c_5\|v_\gamma^n-q^n\cdot\yref(\gamma)\|^2+\|u_\gamma^n\|^2)\ee^{-\sigma (t-\gamma)}+2K_1\sigma^{-1}\right)\\
&\le \varphi_n(t)^2\left(c_5\|v_\gamma-q\cdot\yref(\gamma)\|^2+\|u_\gamma\|^2+2K_1\sigma^{-1}\right).
\end{align*}
Thus there exist $M,N>0$ which are independent of $n$ and $t$ such that
\begin{equation}\label{eq:L2_bound}
\forall\, t\in[\gamma,T_n):\
\|z_n(t)\|^2\leq M\varphi_n(t)^2\ \ \text{and}\ \
\|w_n(t)\|^2\leq N\varphi_n(t)^2,
\end{equation}
and, as a consequence,
\begin{equation}\label{eq:L2_bound_uv}
\forall\, t\in[\gamma,T_n):\ \|v_n(t)-q^n\cdot\yref(t)\|^2\leq M\ \ \text{and}\ \ \|u_n(t)\|^2\leq N.
\end{equation}

\emph{Step 3: We show $T_n=\infty$ and that $e_n$ is uniformly bounded away from~1 on~$[\gamma,\infty)$.}\\
\emph{Step 3a: We derive some estimates for $\ddt\|z_n\|^2$ and for an integral involving $\|z_n\|^4_{L^4}$.} In a similar way in which we have derived~\eqref{eq:Lyapunov_boundary_2} we can obtain the estimate
\begin{equation}\label{eq:energy_boundary_z}
\begin{aligned}
\tfrac{1}{2}\ddt\|z_n\|^2\leq&-\fa(z_n,z_n)-(c_1-\omega_0)\|z_n\|^2+\|z_n\|\|w_n\|\\
&-\frac{k_0\|e_n\|^2_{\R^m}}{1-\|e_n\|^2_{\R^m}}-\frac{c_3}{2\varphi_n^2}\|z_n\|^4_{L^4}+K_0\varphi_n^2.
\end{aligned}
\end{equation}
Using \eqref{eq:L2_bound} and $-c_1\|z_n\|^2\le 0$ leads to
\begin{align*}
\tfrac{1}{2}\ddt\|z_n\|^2\leq&-\fa(z_n,z_n)-\frac{k_0\|e_n\|^2_{\R^m}}{1-\|e_n\|^2_{\R^m}}-\frac{c_3}{2\varphi_n^2}\|z_n\|^4_{L^4}\\
&+\|\dot{\varphi}\|_\infty M\varphi_n+(K_0+\sqrt{MN})\varphi_n^2.
\end{align*}
Hence,
\begin{equation}\label{eq:energy_boundary_2}
\begin{aligned}
\tfrac{1}{2}\ddt\|z_n\|^2\leq&-\fa(z_n,z_n)-\frac{k_0\|e_n\|^2_{\R^m}}{1-\|e_n\|^2_{\R^m}}-\frac{c_3}{2\varphi_n^2}\|z_n\|^4_{L^4}+K_1\varphi_n+K_2\varphi_n^2
\end{aligned}
\end{equation}
on $[\gamma,T_n)$, where $K_1\coloneqq  M\|\dot{\varphi}\|_\infty$ and $K_2\coloneqq  K_0+\sqrt{MN}$.
Observe that
\[
\frac{c_3}{2}\varphi_n^{-3}\|z_n\|^4_{L^4}\leq-\frac{\varphi_n^{-1}}{2}\ddt\|z_n\|^2+K_3,
\]
where $K_3\coloneqq  K_1+K_2\|\varphi\|_\infty$. Therefore,
\begin{align*}
&\frac{c_3}{2}\int_\gamma^t\ee^s\varphi_n(s)^{-3}\|z_n(s)\|^4_{L^4}\ds{s}\\
&\le K_3(\ee^t-\ee^\gamma)-\frac{1}{2}\int_\gamma^t\ee^s\varphi_n(s)^{-1}\ddt\|z_n(s)\|^2\ds{s}\\
&=  K_3(\ee^t-\ee^\gamma)-\frac{1}{2}\left(\ee^t\varphi_n(t)^{-1}\|z_n(t)\|^2-\frac{\|z_\gamma^n\|^2}{\kappa_n}\ee^\gamma\right)\\
&\quad +\frac{1}{2}\int_\gamma^t\ee^s\varphi_n(s)^{-2}(\varphi_n(s)-\dot{\varphi}_n(s))\|z_n(s)\|^2\ds{s}\\
&\le \frac{\ee^t}{2}(2K_3+(\|\varphi\|_\infty+\Gamma_r^{-1}+\|\dot{\varphi}\|_\infty)M)+\kappa_n\ee^\gamma(\|v_\gamma\|^2+\|q\cdot\yref(\gamma)\|^2),
\end{align*}
and hence there exist $D_0,D_1>0$ independent of $n$ and $t$ such that
\begin{equation}
\label{eq:expz4_boundary}
\forall\, t\in[\gamma,T_n):\ \int_\gamma^t\ee^s\varphi_n(s)^{-3}\|z_n(s)\|^4_{L^4}\ds{s}\leq D_1\ee^t+\kappa_n D_0.
\end{equation}
\emph{Step 3b: We derive an estimate for $\|\dot z_n\|^2$.} Multiplying the first equation in~\eqref{eq:weak} by $\dot{\mu}_j$ and {summing over} $j\in\{0,\ldots,n\}$ we obtain
\begin{align*}
\|\dot{z}_n\|^2=&-\frac{1}{2}\ddt\fa(z_n,z_n)-\frac{c_1}{2}\ddt\|z_n\|^2+\frac{k_0}{2}\ddt\ln(1-\|e_n\|^2_{\R^m})\\
&+ \scpr{\omega_0z_n+F\left(t,z_n\right)-w_n}{\dot{z}_n}.
\end{align*}
We can estimate the last term above by
\begin{align*}
\scpr{\omega_0z_n}{\dot{z}_n}\leq&\ \frac{7}{2}\|\dot{\varphi}\|_\infty^2\varphi_n^{-2}\|z_n\|^2+\frac{1}{14}\|\dot{z}_n\|^2 \stackrel{\eqref{eq:L2_bound}}{\leq} \frac{7}{2}\|\dot{\varphi}\|_\infty^2M+\frac{1}{14}\|\dot{z}_n\|^2,\\
\scpr{-w_n}{\dot{z}_n}\leq&\ \frac{7}{2}\|w_n\|^2+\frac{1}{14}\|\dot{z}_n\|^2,\\
\scpr{F\left(t,z_n\right)}{\dot{z}_n}\leq&\ \frac{7}{2}\varphi_n^2\left(m\sum_{k=1}^m \|\yrefk\|_\infty^2 \|\cA q_k\|^2+\|I_{s,i}\|^2_{2,\infty}+\|f_0\|_{\infty,\infty}^2|\Omega|\right)\\
&+\frac{7}{2}\|f_1\|^2_{\infty,\infty}\|z_n\|^2+\frac{7}{2}\varphi_n^{-2}\|f_2\|_{\infty,\infty}^2\|z_n\|^4_{L^4}\\
&+\frac{5}{14}\|\dot{z}_n\|^2-\frac{c_3}{4\varphi_n^2}\ddt\|z_n\|^4_{L^4}.
\end{align*}
Inserting these inequalities, substracting $\tfrac12 \|\dot{z}_n\|^2$ and then multiplying by~$2$ gives
\begin{equation*}\label{eq:energy_boundary_1}
\begin{aligned}
\|\dot{z}_n\|^2=&-\ddt\fa(z_n,z_n)-c_1\ddt\|z_n\|^2+k_0\ddt\ln(1-\|e_n\|^2_{\R^m})-\frac{c_3}{2\varphi_n^2}\ddt\|z_n\|^4_{L^4}\\
&+7\varphi_n^2\left(m\sum_{k=1}^m\|\yrefk\|_\infty^2 \|\cA q_k\|^2+\|I_{s,i}\|^2_{2,\infty}\!+\|f_0\|_{\infty,\infty}^2|\Omega|+\|f_1\|^2_{\infty,\infty}M\!+\!N\!\right)\\
&+7\|\dot{\varphi}\|_\infty^2M+7\varphi_n^{-2}\|f_2\|_{\infty,\infty}^2\|z_n\|^4_{L^4}.
\end{aligned}
\end{equation*}
Now we add and subtract $\frac{1}{2}\ddt\|z_n\|^2$, thus we obtain
\begin{align*}
\|\dot{z}_n\|^2\leq&-\ddt\fa(z_n,z_n)-\left(c_1+\frac{1}{2}\right)\ddt\|z_n\|^2+k_0\ddt\ln(1-\|e_n\|^2_{\R^m}) -\frac{c_3}{2\varphi_n^2}\ddt\|z_n\|^4_{L^4}\\
&+7(\|\varphi\|_\infty+\Gamma_r^{-1})^2\left(m\sum_{k=1}^m\|\yrefk\|_\infty^2\|\cA q_k\|^2+\|I_{s,i}\|^2_{2,\infty}+\|f_0\|_{\infty,\infty}^2|\Omega|\right.\\
&\left.\phantom{\sum_{i=0}^n}\hspace*{-6mm}+\|f_1\|^2_{\infty,\infty}M\right)+7(N(\|\varphi\|_\infty+\Gamma_r^{-1})^2 +\|\dot{\varphi}\|_\infty^2M)+7\varphi_n^{-2}\|f_2\|_{\infty,\infty}^2\|z_n\|^4_{L^4}\\
&+\frac{1}{2}\ddt\|z_n\|^2.
\end{align*}
By the product rule we have
\[-\frac{c_3}{2\varphi_n^2}\ddt\|z_n\|^4_{L^4}=-\ddt\left(\frac{c_3}{2\varphi_n^2}\|z_n\|^4_{L^4}\right)- c_3\varphi_n^{-3}\dot{\varphi_n}\|z_n\|^4_{L^4},\]
thus we find that
\begin{equation}\label{eq:energy_boundary_3}
\begin{aligned}
\|\dot{z}_n\|^2&+\ddt\fa(z_n,z_n)-k_0\ddt\ln(1-\|e_n\|^2_{\R^m})+\ddt\left(\frac{c_3}{2\varphi_n^2}\|z_n\|^4_{L^4}\right)\\
\leq&-\left(c_1+\frac{1}{2}\right)\ddt\|z_n\|^2+E_1+E_2\varphi_n^{-3}\|z_n\|^4_{L^4}+\frac{1}{2}\ddt\|z_n\|^2,
\end{aligned}
\end{equation}
where
\begin{align*}
E_1&\coloneqq  \ 7(\|\varphi\|_\infty+\Gamma_r^{-1})^2\left(m\sum_{k=1}^m\|\yrefk\|_\infty^2\|Aq_k\|^2+ \|I_{s,i}\|^2_{2,\infty}+\|f_0\|_{\infty,\infty}^2|\Omega|\right.\\
&\quad\ \ \left.\phantom{\sum_{i=0}^n}\hspace*{-6mm}+\|f_1\|^2_{\infty,\infty}M\right) +7\big(N(\|\varphi\|_\infty+\Gamma_r^{-1})^2+\|\dot{\varphi}\|_\infty^2M\big),\\
E_2&\coloneqq  7\|f_2\|_{\infty,\infty}^2(\|\varphi\|_\infty+\Gamma_r^{-1})+c_3\|\dot{\varphi}\|_\infty
\end{align*}
are independent of $n$ and $t$.\\
\emph{Step 3c: We show uniform boundedness of~$e_n$.} Using~\eqref{eq:energy_boundary_2} in~\eqref{eq:energy_boundary_3} we obtain
\begin{align*}
\|\dot{z}_n\|^2+\dot\rho_n\leq&-\left(c_1+\frac{1}{2}\right)\ddt\|z_n\|^2+E_1+E_2\varphi_n^{-3}\|z_n\|^4_{L^4}\\
&-\fa(z_n,z_n)-\frac{k_0\|e_n\|^2_{\R^m}}{1-\|e_n\|^2_{\R^m}}-\frac{c_3}{2\varphi_n^2}\|z_n\|^4_{L^4}+K_1\varphi_n+K_2\varphi_n^2\\
=&-\left(c_1+\frac{1}{2}\right)\ddt\|z_n\|^2+E_2\varphi_n^{-3}\|z_n\|^4_{L^4}\\
&-\fa(z_n,z_n)-\frac{k_0}{1-\|e_n\|^2_{\R^m}}-\frac{c_3}{2\varphi_n^2}\|z_n\|^4_{L^4}+\Lambda,
\end{align*}
where
\begin{align*}
\rho_n&\coloneqq  \fa(z_n,z_n)-k_0\ln(1-\|e_n\|^2_{\R^m})+\frac{c_3}{2\varphi_n^2}\|z_n\|^4_{L^4},\\
\Lambda&\coloneqq  E_1+K_1(\|\varphi\|_\infty+\Gamma_r^{-1})+K_2(\|\varphi\|_\infty+\Gamma_r^{-1})^2+k_0,
\end{align*}
and we have used the equality
$$\frac{\|e_n\|^2_{\R^m}}{1-\|e_n\|^2_{\R^m}}=-1+\frac{1}{1-\|e_n\|^2_{\R^m}}.$$
Adding and subtracting $k_0\ln(1-\|e_n\|^2_{\R^m})$ leads to
\begin{align}
\|\dot{z}_n\|^2+\dot\rho_n\leq&-\rho_n-\left(c_1+\frac{1}{2}\right)\ddt\|z_n\|^2+E_2\varphi_n^{-3}\|z_n\|^4_{L^4}\notag\\
&-k_0\left(\frac{1}{1-\|e_n\|^2_{\R^m}}+\ln(1-\|e_n\|^2_{\R^m})\right)+\Lambda\notag\\
\leq&-\rho_n-\left(c_1+\frac{1}{2}\right)\ddt\|z_n\|^2+E_2\varphi_n^{-3}\|z_n\|^4_{L^4}+\Lambda, \label{eq:L2zdot_boundary}
\end{align}
where for the last inequality we have used that
\[\forall\, p\in(-1,1):\ \frac{1}{1-p^2}\geq\ln\left(\frac{1}{1-p^2}\right) = -\ln(1-p^2).\]
We may now use the integrating factor $\ee^t$ to obtain
\[
\ddt \left(\ee^t\rho_n\right) = \ee^t(\rho_n + \dot \rho_n) \leq -\ee^t\left(c_1+\frac{1}{2}\right)\ddt\|z_n\|^2+E_2\ee^t\varphi_n^{-3}\|z_n\|^4_{L^4}+\Lambda \ee^t \underset{\le 0}{\underbrace{- \ee^t \|\dot z_n\|^2}}.
\]
Integrating and using \eqref{eq:expz4_boundary} yields that for all $t\in[\gamma,T_n)$ we have
\begin{align*}
\ee^t\rho_n(t)-\rho_n(\gamma)\ee^\gamma\leq&\ (E_2D_1+\Lambda)\ee^t+\kappa_n E_2D_0-\int_\gamma^t\ee^s\left(c_1+\frac{1}{2}\right)\ddt\|z_n(s)\|^2\ds{s}\\
\leq&\ (E_2D_1+\Lambda)\ee^t+\kappa_n E_2D_0+\left(c_1+\frac{1}{2}\right)\|z_\gamma^n\|^2\ee^\gamma\\
&+\left(c_1+\frac{1}{2}\right)\int_\gamma^t\ee^s\|z_n(s)\|^2\ds{s}\\
\stackrel{\eqref{eq:L2_bound}}{\leq}&\ (E_2D_1+\Lambda)\ee^t+\kappa_n E_2D_0+\left(c_1+\frac12\right)\kappa_n^2\ee^\gamma(\|v_\gamma-q\cdot\yref(\gamma)\|^2)\\
&+\left(c_1+\frac{1}{2}\right)(\|\varphi\|_\infty+\Gamma_r^{-1})^2M\ee^t.
\end{align*}
Thus, there exit $\Xi_1,\Xi_2,\Xi_3>0$ independent of $n$ and $t$, such that
\[\rho_n(t)\leq\rho_n(\gamma)\ee^{-(t-\gamma)}+\Xi_1+\kappa_n(\Xi_2+\kappa_n\Xi_3)\ee^{-(t-\gamma)}.\]
Invoking the definition of $\rho_n$ and that $\ee^{-(t-\gamma)}\leq1$ for $t\geq\gamma$ we find that
\begin{equation}\label{eq:rho}
\forall\, t\in[\gamma,T_n):\  \rho_n(t)\leq \rho_{n}^0+\Xi_1+\kappa_n\Xi_2+\kappa_n^2\Xi_3,
\end{equation}
where
\begin{align*}
\rho_{n}^0\coloneqq  &\ \kappa_n^2\fa(v_\gamma^n\!-\!q^n\cdot\yref(\gamma),v_\gamma^n\!-\!q^n\cdot\yref(\gamma))\!-\!k_0\ln(1\!-\!\kappa_n^2\|\cB' (v_\gamma^n\!-\!q^n\cdot\yref(\gamma))\|^2_{\R^m})\\
&+\kappa_n^2\|v_\gamma^n-q^n\cdot\yref(\gamma)\|_{L^4}^4 = \rho_n(\gamma).
\end{align*}
Note that by construction of $\kappa_n$ and the Sobolev embedding theorem, $(\rho_n^0)_{n\in\N}$ is bounded, $\rho_n^0\to0$ as $n\to\infty$, so that $\rho_n^0$ can be bounded independently of $n$.\\
Again using the definition of $\rho_n$ and~\eqref{eq:rho} we find that
\begin{align*}
k_0 \ln\left(\frac{1}{1-\|e_n\|^2_{\R^m}}\right)  = \rho_n - \fa(z_n,z_n)- \frac{c_3}{2\varphi_n^2}\|z_n\|^4_{L^4}
\le \rho_{n}^0+\Xi_1+\kappa_n\Xi_2+\kappa_n^2\Xi_3,
\end{align*}
and hence
\[\frac{1}{1-\|e_n\|^2_{\R^m}}\leq\exp\left(\frac{1}{k_0}\left(\rho_{n}^0+\Xi_1+\kappa_n\Xi_2+\kappa_n^2\Xi_3\right)\right)=:\varepsilon(n).\]
We may thus conclude that
\begin{equation}\label{eq:err_bounded_away}
\forall\, t\in[\gamma,T_n):\ \|e_n(t)\|^2_{\R^m}\leq1-\varepsilon(n),
\end{equation}
or, equivalently,
\begin{equation}\label{eq:err_bdd}
\forall\, t\in[\gamma,T_n):\ \varphi_n(t)^2\|\cB'( v_n(t)-q^n\cdot\yref(t))\|^2_{\R^m}\leq1-\varepsilon(n).
\end{equation}
Moreover, from~\eqref{eq:rho}, the definition of $\rho$, $k_0\ln(1-\|e_n\|^2_{\R^m})\leq0$ and Assumption~\ref{Ass1} we have that
\begin{align*}
\delta\|\nabla z_n\|^2+\frac{c_3}{2\varphi_n^2}\|z_n\|^4_{L^4}
\leq\rho_{n}^0+\Xi_1+\kappa_n\Xi_2+\kappa_n^2\Xi_3.
\end{align*}
Reversing the change of variables leads to
\begin{equation}\label{eq:potential}
\begin{aligned}
\forall\, t\in[\gamma,T_n):\ \delta\varphi_n(t)^2\|\nabla(v_n(t)-q^n\cdot\yref(t))\|^2&+\varphi_n(t)^2\|v_n(t)-q^n\cdot\yref(t)\|_{L^4}^4\\
&\leq\rho_{n}^0+\Xi_1+\kappa_n\Xi_2+\kappa_n^2\Xi_3,
\end{aligned}
\end{equation}
which implies that for all $t\in[\gamma,T_n)$ we have $v_n(t)\in W^{1,2}(\Omega)$.\\
\emph{Step 3d: We show that~$T_n=\infty$.} Assuming $T_n<\infty$ it follows from~\eqref{eq:err_bounded_away} that the graph of the solution~$(\mu^n,\nu^n)$ from Step~2 would be a compact subset of~$\D$, a contradiction. Therefore, we have $T_n=\infty$.

\emph{Step 4: We show convergence of the approximate solution, uniqueness and regularity of the solution in $[\gamma,\infty)\times \Omega$.}\\
\emph{Step~4a: we prove some inequalities for later use.} From~\eqref{eq:rho} we have that, on $[\gamma,\infty)$,
\[\varphi_n^{-2}\|z_n\|_{L^4}^4\leq\rho_{n}^0+\Xi_1+\kappa_n\Xi_2 + \kappa_n^2\Xi_3.\]
Using a similar procedure as for the derivation of~\eqref{eq:expz4_boundary} we may obtain the estimate
\begin{equation}\label{eq:noexpz4_boundary}
\forall\, t\ge 0:\ \int_\gamma^t\varphi_n(s)^{-3}\|z_n(s)\|^4_{L^4}\ds{s}\leq\kappa_n d_0+d_1t
\end{equation}
for $d_0,d_1>0$ independent of $n$ and $t$.
Further, we can integrate \eqref{eq:L2zdot_boundary} on the interval $[\gamma,t]$ to obtain, invoking $\rho_n(t)\ge 0$ and~\eqref{eq:noexpz4_boundary},
\[
\int_\gamma^t\|\dot{z}_n(s)\|^2\ds{s}\leq\rho_{n}^0+ \left(c_1+\frac12\right)\kappa_n^2(\|v_\gamma-q\cdot\yref(\gamma)\|^2)+E_2(\kappa_n d_0+d_1t)+\Lambda t
\]
for all $t\ge \gamma$. Hence, there exist $S_0,S_1,S_2>0$ independent of $n$ and $t$ such that
\begin{equation}\label{eq:intzdot}
\forall\, t\ge \gamma:\ \int_\gamma^t\|\dot{z}_n(s)\|^2\ds{s}\leq\rho_{n}^0+S_0\kappa_n+S_1\kappa_n^2+S_2t.
\end{equation}
This implies existence of $S_3,S_4>0$ such that
\begin{equation}\label{eq:dot_var_v}
\forall\, t\ge \gamma:\ \int_\gamma^t\left\|\ddt(\varphi_nv_n)\right\|^2\ds{s}\leq\rho_{n}^0+S_0\kappa_n+S_1\kappa_n^2+S_3t+S_4.
\end{equation}
In order to improve~\eqref{eq:noexpz4_boundary}, we observe that from~\eqref{eq:energy_boundary_z} it follows
\begin{align*}
\tfrac{1}{2}\ddt\|z_n\|^2\leq&-\fa(z_n,z_n)-(c_1-\omega_0)\|z_n\|^2+\|z_n\|\|w_n\|\\
&-\frac{k_0\|e_n\|^2_{\R^m}}{1-\|e_n\|^2_{\R^m}}-\frac{c_3}{2\varphi_n^2}\|z_n\|^4_{L^4}+K_0\varphi_n^2\\
\leq&\ \omega_0\|z_n\|^2-\frac{c_3}{2\varphi_n^2}\|z_n\|^4_{L^4}+K_2\varphi_n^2 -\fa(z_n,z_n)-\frac{k_0\|e_n\|^2_{\R^m}}{1-\|e_n\|^2_{\R^m}},
\end{align*}
which gives
\[\ddt\varphi_n^{-2}\|z_n\|^2\leq 2K_2-c_3\varphi_n^{-4}\|z_n\|^4_{L^4}-2\varphi_n^{-2}\fa(z_n,z_n)-\frac{2k_0\varphi_n^{-2}\|e_n\|^2_{\R^m}}{1-\|e_n\|^2_{\R^m}}.\]
This implies that for all $t\ge \gamma$ we have
\begin{equation}\label{eq:vphi4z4}
\begin{aligned}
\int_\gamma^t c_3\varphi_n(s)^{-4}\|z_n(s)\|^4_{L^4}
&+2\varphi_n(s)^{-2}\fa(z_n(s),z_n(s))+\frac{2k_0\varphi_n(s)^{-2}\|e_n(s)\|^2_{\R^m}}{1-\|e_n(s)\|^2_{\R^m}}\ds{s}\\
&\leq2K_2t+\|v_\gamma-q\cdot\yref(\gamma)\|^2,
\end{aligned}
\end{equation}
which is bounded independently of $n$. This shows that for all $t\ge \gamma$ we have
\begin{equation}\label{eq:L4}
\begin{aligned}
c_3\int_\gamma^t\|v_n(s)-q^n\cdot\yref(s)\|^4_{L^4}\ds{s}+\int_\gamma^t2\fa(v_n(s)-q^n\cdot \yref(s),v_n(s)-q^n\cdot \yref(s))\ds{s}\\
+\int_\gamma^t\frac{2k_0\|\cB' (v_n(s)-q^n\cdot\yref(s))\|^2_{\R^m}}{1-\varphi_n(s)^{2}\|\cB' (v_n(s)-q^n\cdot\yref(s))\|^2_{\R^m}}\ds{s}\leq2K_2t+\|v_\gamma-q\cdot\yref(\gamma)\|^2.
\end{aligned}
\end{equation}
In order to prove that $\|\dot{w}_n\|^2$ is bounded independently of $n$ and $t$, a last calculation is required. Multiply the second equation in~\eqref{eq:weak} by $\dot{\nu}_j$ and sum over $j$ to obtain
\[\|\dot{w}_n\|^2=-(c_4-\omega_0)\scpr{w_n}{\dot{w}_n}+c_5\scpr{z_n}{\dot{w}_n}+\varphi_n \scpr{g}{\dot{w}_n}.\]
Using $(\omega_0 - c_4) w_n = (\dot\varphi_n - c_4\varphi_n) \varphi_n^{-1} w_n$ and the inequalities
\begin{align*} -(c_4-\omega_0)\scpr{w_n}{\dot{w}_n}&\leq\frac{3}{2} \| \dot\varphi - c_4\varphi\|_\infty^2 \varphi_n^{-2}\|w_n\|^2+\frac{\|\dot{w}_n\|^2}{6}\\
&\leq\frac{3}{2}(\| \dot{\varphi}\|_\infty + c_4(\|\varphi\|_\infty+\Gamma_r^{-1}))^2 N+\frac{\|\dot{w}_n\|^2}{6},\\
c_5\scpr{z_n}{\dot{w}_n}&\leq\frac{3c_5^2}{2}\|z_n\|^2+\frac{1}{6}\|\dot{w}_n\|^2\\
&\leq\frac{3c_5^2M}{2}(\|\varphi\|_\infty+\Gamma_r^{-1})^2+\frac{1}{6}\|\dot{w}_n\|^2,\\
\varphi_n \scpr{g}{\dot{w}_n}& \leq\frac{3}{2}(\|\varphi\|_\infty+\Gamma_r^{-1})^2\|g\|^2_{\infty,\infty}|\Omega|+\frac{1}{6}\|\dot{w}_n\|^2,
\end{align*}
it follows that for all $t\ge\gamma$ we have
\begin{equation}\label{eq:dotw}
\begin{aligned}
\|\dot{w}_n(t)\|^2\leq&3\|(\| \dot{\varphi}\|_\infty + c_4(\|\varphi\|_\infty+\Gamma_r^{-1}))^2 N\\&+3c_5^2M(\|\varphi\|_\infty+\Gamma_r^{-1})^2+3(\|\varphi\|_\infty+\Gamma_r^{-1})^2\|g\|^2_{\infty,\infty}|\Omega|,\end{aligned}
\end{equation}
which is bounded independently of $n$ and $t$. Multiplying the second equation in~\eqref{eq:weak} by $\varphi_n^{-1}$ and $\theta_i$ and {summing over} $i\in\{0,\dots,n\}$ leads to
\[\ddt(\varphi_n^{-1}w_n)=-\varphi^{-2}\dot{\varphi}_nw_n+\varphi_n^{-1}\dot{w}_n=-c_4\varphi_n^{-1}w_n+c_5\varphi_n^{-1}z_n+g_n,\]
where
\[g_n\coloneqq  \sum_{i=0}^n\scpr{g}{\theta_i}\theta_i.\]
Taking the norm of the latter gives
\begin{align*}
\left\|\ddt(\varphi_n^{-1}w_n)\right\|&\leq c_4\varphi_n^{-1}\|w_n\|+c_5\varphi_n^{-1}\|z_n\|+\|g_n\|\\
&\leq c_4N+c_5M+\|g\|_{\infty,\infty},
\end{align*}
thus
\begin{equation}\label{eq:dot_u}
\forall\, t\ge \gamma:\ \|\dot{u}_n(t)\|\leq c_4N+c_5M+\|g\|_{\infty,\infty}.
\end{equation}
\emph{Step 4b: We show that $(v_n,u_n)$ converges weakly.} Let $T>\gamma$ be given. Using a similar argument as in Section~\ref{ssec:mono_proof_tleqgamma}, we have that $v_n\in L^2(\gamma,T;W^{1,2}(\Omega))$ and $\dot{v}_n\in L^2(\gamma,T;W^{1,2}(\Omega)')$, since~\eqref{eq:L4} together with~\eqref{eq:err_bdd} implies that $I_{s,e}^n\in L^2(\gamma,T;\R^m)$ and $v_n\in L^2(\gamma,T;W^{1,2}(\Omega))$.\\
Furthermore, analogously to Section~\ref{ssec:mono_proof_tleqgamma}, we have that there exist subsequences such that
\begin{align*}
u_n\to u&\in W^{1,2}(\gamma,T;L^{2}(\Omega))\mbox{ weakly},\\
v_n\to v&\in L^2(\gamma,T;W^{1,2}(\Omega))\mbox{ weakly},\\
\dot{v}_n\to\dot{v}&\in L^2(\gamma,T;(W^{1,2}(\Omega))')\mbox{ weakly},
\end{align*}
so that $u,v\in C([\gamma,T];L^2(\Omega))$. Also $v_n^2\to v^2$ weakly in $L^2((\gamma,T)\times\Omega)$ and $v_n^3\to v^3$ weakly in $L^{4/3}((\gamma,T)\times\Omega)$.\\
We may infer further properties of $u$ and $v$. By \eqref{eq:L2_bound_uv}, \eqref{eq:potential}, \eqref{eq:dot_var_v} \& \eqref{eq:dot_u} we have that $u_n,\dot{u}_n$ lie in a bounded subset of $L^\infty(\gamma,\infty;L^2(\Omega))$ and that $v_n$ lie in a bounded subset of $L^\infty(\gamma,\infty;L^2(\Omega))$. Moreover, $\ddt(\varphi_nv_n)\in L^2_{\rm loc}(\gamma,\infty;L^2(\Omega))$. Then, using Lemma~\ref{lem:weak_convergence}, we find a subsequence such that
\begin{align*}
u_n\to u&\in L^\infty(\gamma,T;L^2(\Omega))\mbox{ weak}^\star,\\
\dot{u}_n\to \dot{u}&\in L^\infty(\gamma,T;L^2(\Omega))\mbox{ weak}^\star,\\
v_n\to v&\in L^\infty(\gamma,T;L^{2}(\Omega))\mbox{ weak}^\star,\\
\varphi_nv_n\to \varphi v&\in L^\infty(\gamma,T;W^{1,2}(\Omega))\mbox{ weak}^\star,\\
\dot{v}_n\to \dot{v}&\in L^2(\gamma,T;W^{1,2}(\Omega)')\mbox{ weakly},\\
\varphi_n\dot{v}_n\to \varphi\dot{v}&\in L^2(\gamma,T;L^2(\Omega))\mbox{ weakly},
\end{align*}
since $\varphi_n\to\varphi$ in $BC([\gamma,T];\R)$. Moreover, by $\inf_{t>\gamma+\delta}\varphi(t)>0$, we also have that $v\in L^\infty(\gamma+\delta,T;W^{1,2}(\Omega))$ and $\dot{v}\in L^2(\gamma+\delta,T;L^2(\Omega))$ for all~$\delta>0$.\\
Further, $\kappa_n,\rho_n^0\to0$ and
\[\varepsilon(n)\underset{n\to\infty}{\to}\varepsilon_0\coloneqq  \exp\left(-k_0^{-1}\Xi_1\right).\]
%Recall that $q^n\to q$ strongly in $W^{1,2}(\Omega;\R^m)$ and $\cB' q_k=\ee_k$.
Thus, by \eqref{eq:L2_bound_uv}, \eqref{eq:err_bdd}, \eqref{eq:potential} \& \eqref{eq:L4} we have $v\in L^4((\gamma,T)\times\Omega)$ and for almost all $t\in[\gamma,T)$ the following estimates hold:
\begin{equation}\label{eq:potential_limit}
\begin{aligned}
&\|v(t)-q\cdot\yref(t)\|\leq \sqrt{M},\\
&\|u(t)\|\leq \sqrt{N},\\
&\varphi(t)^2\|\cB' v(t)-\yref(t)\|^2_{\R^m}\leq1-\varepsilon_0,\\
&\delta\varphi(t)^2\|\nabla(v(t)-q\cdot\yref(t))\|^2+\varphi(t)^2\|v(t)-q\cdot\yref(t)\|_{L^4}^4\leq\Xi_1,\\
&\int_\gamma^t\|v(s)-q\cdot\yref(s)\|^4_{L^4}\ds{s}\leq2K_2t+\|v_\gamma-q\cdot\yref(\gamma)\|^2.
\end{aligned}
\end{equation}
%since from \eqref{eq:potential} and $\kappa_n,\rho_n^0\to0$ it follows
%\begin{align*}
%	\forall\, t\in[\gamma,T):\ \delta\varphi(t)^2\|\nabla(v(t)-q\cdot\yref(t))\|^2&+\varphi(t)^2\|v(t)-q\cdot\yref(t)\|_{L^4}^4\\
%	&\leq\lim_{n\to\infty} \big(\rho_{n}^0+\Xi_1+\kappa_n\Xi_2+\kappa_n^2\Xi_3\big).
%\end{align*}
Moreover, as in Section~\ref{ssec:mono_proof_tleqgamma}, $v_n\to v$ strongly in $L^2(\gamma,T;L^2(\Omega))$ and $u,v\in C([\gamma,T);L^2(\Omega))$ with $(u(\gamma),v(\gamma))=(u_\gamma,v_\gamma)$.\\
%Note that $\|q-q^n\|_{W^{1,2}(\Omega)}\to0$ as $n\to\infty$. Further, $q^n,q\in\cD(\cA)$ and we have that
%\[\scpr{\cA(q_k-q_k^n)}{\theta}=-\fa(q_k-q_k^n,\theta)\leq\|D\|_{L^\infty}\|q_k-q_k^n\|_{W^{1,2}}\|\theta\|_{W^{1,2}}\quad\forall\theta\in W^{1,2}(\Omega),\]
%so that $\cA(q_k-q_k^n)\to0$ weakly.
Hence, for $\chi\in L^2(\Omega)$ and $\theta\in W^{1,2}(\Omega)$ we have that $(u_n,v_n)$ satisfy the integrated version of \eqref{eq:weak_uv}, thus we obtain that for $t\in(\gamma,T)$
\begin{align*}
\scpr{v(t)}{\theta}=&\ \scpr{v_\gamma}{\theta}+\int_\gamma^t-\fa(v(s),\theta)+\scpr{p_3(v(s))-u(s)+I_{s,i}(s)}{\theta}\ds{s},\\
&+\int_\gamma^T\scpr{I_{s,e}(s)}{\cB'\theta}_{\R^m}\ds{s},\\
\scpr{u(t)}{\chi}=&\ \scpr{u_\gamma}{\chi}+\int_\gamma^t\scpr{c_5v(s)-c_4u(s)}{\chi}\ds{s},\\
I_{s,e}(t)=&\ -\frac{k_0}{1-\varphi(t)^2\|\cB' v(t)-\yref(t)\|^2_{\R^m}}(\cB' v(t)-\yref(t))
\end{align*}
by bounded convergence \cite[Thm.~II.4.1]{Dies77}.
%\marginpar{Zitat f\"ur ``bounded convergence''}
Hence, $(u,v)$ is a solution of~\eqref{eq:FHN_feedback} in $(\gamma,T)$. Moreover,~\eqref{eq:strong_delta} also holds in $W^{1,2}(\Omega)'$ for $t\geq\gamma$, that is
\begin{equation}\label{eq:Xminus}
\dot{v}(t)=\cA v(t)+p_3(v(t))+\cB I_{s,e}(t)-u(t)+I_{s,i}(t).
\end{equation}

{\emph{Step 5: We show uniqueness of the solution on $[0,\infty)$.}\\
The proof is similar, but, in an essential step, also different from the proof given in Step~1e of Section~\ref{ssec:mono_proof_tleqgamma}. Let $T>0$ and assume that $(v_1,u_1)$ and $(v_2,u_2)$ are two solutions of~\eqref{eq:FHN_feedback} on $[\gamma,T)$ with $v_1(\gamma) = v_2(\gamma) = v_\gamma$ and $u_1(\gamma) = u_2(\gamma) = u_\gamma$. Choose~$\hat p_3$, $Q_\gamma := (\gamma,T)\times\Omega$, $\Sigma$, $\Lambda$ and $Q^\Lambda$ similar to Step~1e of Section~\ref{ssec:mono_proof_tleqgamma}, where we invoke that $v_1,v_2\in L^4((\gamma,T)\times\Omega)$. Let $V:= v_2-v_1$ and $U:= u_2-u_1$, then, by~\eqref{eq:FHN_feedback},
\begin{align*}
  \dot V &= (\cA-c_1I) V -c_3(\hat{p}_3(v_2) - \hat{p}_3(v_1)) - U \\
  &\quad -k_0 \cB \left(\frac{\cB'v_2-\yref}{1-\varphi^2\|\cB' v_2-\yref\|^2_{\R^m}}- \frac{\cB'v_1-\yref}{1-\varphi^2\|\cB' v_1-\yref\|^2_{\R^m}}\right),\\
  \dot U &= c_5 V - c_4 U.
\end{align*}
Define
\[
    \Xi(t) := \frac{\cB'v_2(t)-\yref(t)}{1-\varphi(t)^2\|\cB' v_2(t)-\yref(t)\|^2_{\R^m}}- \frac{\cB'v_1(t)-\yref(t)}{1-\varphi(t)^2\|\cB' v_1(t)-\yref(t)\|^2_{\R^m}} \in\R^m
\]
for $t\in (\gamma,T)$, then we may compute that
\begin{align*}
     \tfrac{c_5}{2}\ddt\|V\|^2+\tfrac{1}{2}\ddt\|U\|^2 &= \scpr{(\cA-c_1I) V - U}{c_5 V}_{{W^{1,2}(\Omega)',W^{1,2}(\Omega)}} - c_4 \|U\|^2  \\
     &\quad + c_5 \langle U,V \rangle -c_5 c_3\scpr{ \hat{p}_3(v_2) - \hat{p}_3(v_1)}{V} - c_5 k_0 \scpr{\Xi}{\cB' V}_{\R^m}\\
    &\le-c_5c_3 \scpr{\hat{p}_3(v_2)-\hat{p}_3(v_1)}{V} - c_5 k_0 \scpr{\Xi}{\cB' V}_{\R^m}.
\end{align*}
Define
\[
    e_i(t) := \cB'v_i(t)-\yref(t),\quad k_i(t):= \frac{1}{1-\varphi(t)^2\|e_i(t)\|^2_{\R^m}},\quad t\in (\gamma,T),\ i=1,2
\]
and observe that on $(\gamma,T)$ we have
\begin{align*}
  k_1 \ge k_2\quad &\iff\quad 1- \varphi^2  \|e_2\|_{\R^m}^2 \ge 1- \varphi^2  \|e_1\|_{\R^m}^2\\
&\iff\quad \|e_1\|_{\R^m}^2 \ge \|e_2\|_{\R^m}^2,
\end{align*}
by which
\[
    \forall\, t\in(\gamma,T):\ \big(k_1(t) - k_2(t) \big) \big(\|e_1(t) \|_{\R^m}^2 - \|e_2(t) \|_{\R^m}^2\big) \ge 0.
\]
Then we may calculate that
\begin{align*}
  \scpr{\Xi}{\cB' V}_{\R^m} &= \scpr{k_2 e_2 - k_1 e_1}{e_2 - e_1}_{\R^m} \\
  &=  k_2 \|e_2\|_{\R^m}^2 + k_1 \|e_1\|_{\R^m}^2 - (k_1 + k_2) \scpr{e_2}{e_1}_{\R^m} \\
 &\ge   k_2 \|e_2\|_{\R^m}^2 + k_1 \|e_1\|_{\R^m}^2 - \tfrac12 (k_1 + k_2) ( \|e_1\|_{\R^m}^2 +  \|e_2\|_{\R^m}^2) \\
&= \tfrac12 (k_2 - k_1) \|e_2\|_{\R^m}^2 + \tfrac12 (k_1 - k_2) \|e_1\|_{\R^m}^2 \\
&= \tfrac12 (k_1 - k_2) \big(\|e_1\|_{\R^m}^2 - \|e_2\|_{\R^m}^2\big) \ge 0.
\end{align*}
Therefore, we have that
\[
     \tfrac{c_5}{2}\ddt\|V\|^2+\tfrac{1}{2}\ddt\|U\|^2  \le-c_5c_3 \scpr{\hat{p}_3(v_2)-\hat{p}_3(v_1)}{V} .
\]
Then the same arguments as in Step~1e of Section~\ref{ssec:mono_proof_tleqgamma} apply to conclude that $V(t)=0$ and $U(t)=0$ for all $t\in (\gamma,T)$, thus $v_1 = v_2$ and $u_1 = u_2$ on $(\gamma,T)$. Combining this with uniqueness on $[0,\gamma]$ and invoking that $T>0$ was arbitrary we obtain a unique solution on $[0,\infty)$.
}

\emph{Step 6: We show the regularity properties of the solution.}\\
To this end, note that for all $\delta>0$ we have that
\[v\in L^2_{\rm loc}(\gamma,\infty;W^{1,2}(\Omega))\cap L^\infty(\gamma+\delta,\infty;W^{1,2}(\Omega)),\]
so that $I_r\coloneqq I_{s,i}+c_2v^2-c_3v^3-u\in L^2_{\rm loc}(\gamma,\infty;L^{2}(\Omega))\cap L^\infty(\gamma+\delta,\infty;L^2(\Omega))$, and the application of Proposition \ref{prop:hoelder} yields that $v\in BC([\gamma,\infty);L^2(\Omega))\cap BUC((\gamma,\infty);W^{1,2}(\Omega))$. By the uniform continuity of $v$ and the completeness of $W^{1,2}(\Omega)$, $v$ has a limit at $t=\gamma$, see for instance \cite[Thm.~II.13.D]{Simm63}. Thus, $v\in L^\infty(\gamma,\infty;W^{1,2}(\Omega))$. From Section \ref{ssec:mono_proof_tleqgamma} and the latter we have that $v\in L^2_{\rm loc}(0,\infty;W^{1,2}(\Omega))\cap L^\infty(\delta,\infty;W^{1,2}(\Omega))$ for all $\delta>0$, so we have
\begin{align*}
I_{s,e}&\in L^2_{\rm loc}(0,\infty;\R^m)\cap L^\infty(\delta,\infty;\R^m),\\
v&\in L^2_{\rm loc}(0,\infty;W^{1,2}(\Omega))\cap L^\infty(\delta,\infty;W^{1,2}(\Omega))\\
&\quad \ \ \cap BC([0,\infty);L^2(\Omega)) \cap BUC([\delta,\infty);W^{1,2}(\Omega)),
\end{align*}
so that $I_r\coloneqq  I_{s,i}+c_2v^2-c_3v^3-u\in L^2_{\rm loc}(0,\infty;L^{2}(\Omega))\cap L^\infty(\delta,\infty;L^2(\Omega))$.\\
Recall that by assumption we have $\cB\in\cL(\R^m,W^{r,2}(\Omega)')$ for some $r\in [0,1]$. Applying Proposition \ref{prop:hoelder} we have that  for all $\delta>0$ the unique solution of~\eqref{eq:Xminus} satisfies
\begin{equation}\label{eq:sol_reg}
\begin{aligned}
\text{if $r=0$:} &\quad \forall\,\lambda\in(0,1):\ 	v\in C^{0,\lambda}([\delta,\infty);L^2(\Omega)); \\
\text{if $r\in(0,1)$:} &\quad
	v\in C^{0,1-r/2}([\delta,\infty);L^2(\Omega));\\
\text{if $r=1$:} &\quad	v\in C^{0,1/2}([\delta,\infty);L^2(\Omega)).
\end{aligned}
\end{equation}
Since $u,v\in BC([0,\infty);L^2(\Omega))$ and $\dot{u}=c_4v-c_5u$, we also have $\dot{u}\in BC([0,\infty);L^2(\Omega))$.

Now, from~\eqref{eq:sol_reg} and $\cB' \in\cL(W^{r,2}(\Omega),\R^m)$ for $r\in[0,1]$ we obtain that
\begin{itemize}
	\item for $r=0$ and $\lambda\in(0,1)$:\ $y=  \cB' v\in C^{0,\lambda}([\delta,\infty);\R^m)$;
	\item for $r\in(0,1)$:\ $y=  \cB' v\in C^{0,1-r}([\delta,\infty);\R^m)$;
	\item for $r=1$:\ $y=  \cB' v\in BUC([\delta,\infty);\R^m)$.	
\end{itemize}
Further, from \eqref{eq:potential_limit} we have
\[\forall\,t\geq\delta:\ \varphi(t)^{2}\|\cB' v(t)-\yref(t)\|^2_{\R^m}\leq1-\varepsilon_0,\]
hence $I_{s,e}\in L^\infty(\delta,\infty;\R^m)$ and $I_{s,e}$ has the same regularity properties as $y$, since we have that $\varphi\in\Phi_\gamma$ and $\yref\in W^{1,\infty}(0,\infty;\R^m)$. Therefore, we have proved statements (i)--(iii) in Theorem~\ref{thm:mono_funnel} as well as~a) and~b).

It remains to show~c), for which we additionally require that $\cB\in\cL(\R^m,W^{1,2}(\Omega))$. Then there exist $b_1,\ldots,b_m\in W^{1,2}(\Omega)$ such that $(\cB' x)_i=\scpr{x}{b_i}$ for all $i=1,\dots,m$ and $x\in L^2(\Omega)$. Using the $b_i$ in the weak formulation for $i=1,\dots,m$, we have
\[\ddt\scpr{v(t)}{b_i}=-\fa(v(t),b_i)+\scpr{p_3(v(t))-u(t)+I_{s,i}(t)}{b_i}+\scpr{I_{s,e}(t)}{\cB' b_i}_{\R^m}.\]
Since $(\cB' v(t))_i=\scpr{v(t)}{b_i}$, this leads to
\[\ddt(\cB' v(t)_i)=-\fa(v(t),b_i)+\scpr{p_3(v(t))-u(t)+I_{s,i}(t)}{b_i}+\scpr{I_{s,e}(t)}{\cB' b_i}_{\R^m}.\]
Taking the absolute value and using the Cauchy-Schwarz inequality yields
\begin{align*}
\left|\ddt(\cB' v(t))_i\right|\leq&\ \|D\|_{L^\infty}\|v(t)\|_{W^{1,2}}\|b_i\|_{W^{1,2}}+\|p_3(v(t))-u(t)+I_{s,i}(t)\|_{L^2}\|b_i\|_{L^2}\\
&+\|I_{s,e}(t)\|_{\R^m}\|\cB' b_i\|_{\R^m},
\end{align*}
and therefore
\begin{align*}
\forall\, i=1,\ldots,m\ \forall\,\delta>:\ \left\|\ddt(\cB' v)_i\right\|_{L^\infty(\delta,\infty;\R^m)}<\infty,
\end{align*}
by which $y = \cB'v \in W^{1,\infty}(\delta,\infty;\R^m)$ as well as $I_{s,e} \in W^{1,\infty}(\delta,\infty;\R^m)$. This completes the proof of the theorem.
\ensuremath{\QED}

\newlength\fheight
\newlength\fwidth
\setlength\fheight{0.3\linewidth}
\setlength\fwidth{0.9\linewidth}

\section{A numerical example}
\label{sec:numerics}

\begin{figure}[H]
\begin{subfigure}{.5\linewidth}
\begin{center}
\includegraphics[scale=0.4]{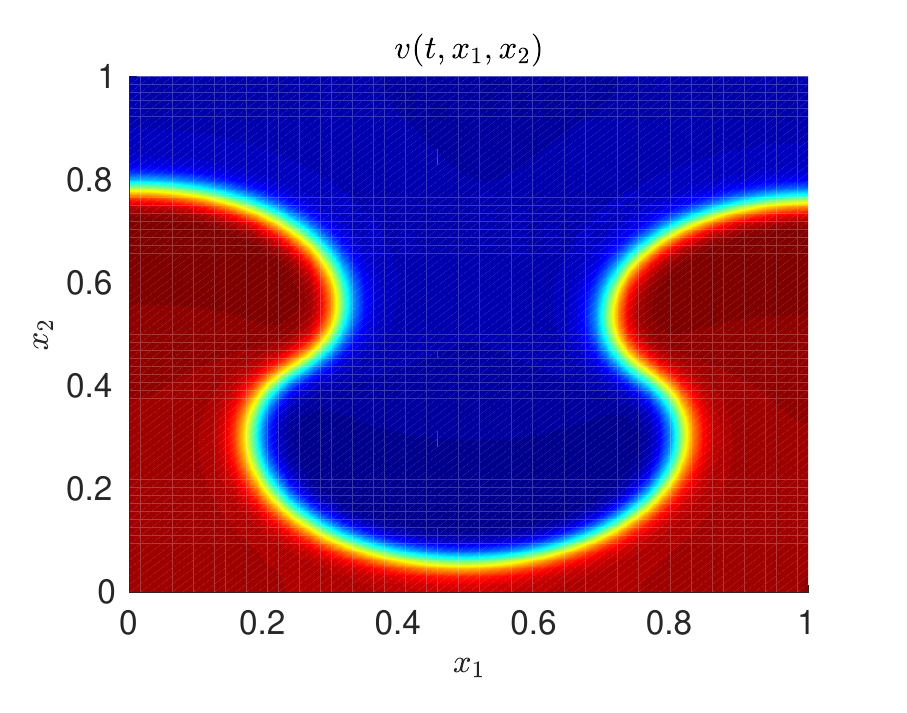}
  %\caption{...}
  %\label{fig:y_funnel_4}
  \end{center}
\end{subfigure}
\begin{subfigure}{.5\linewidth}
\begin{center}
  \includegraphics[scale=0.4]{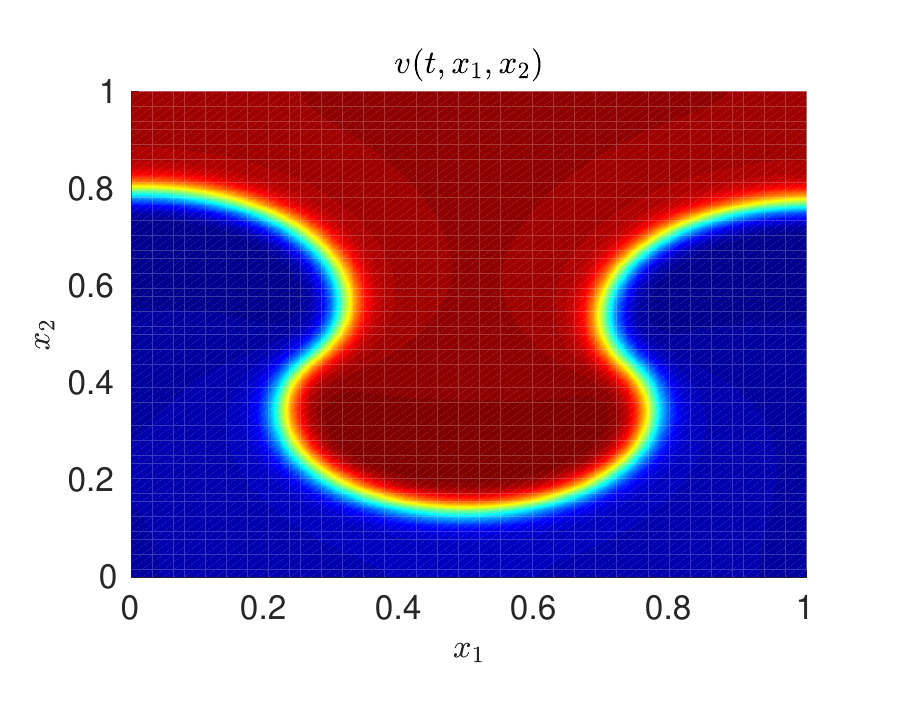}
  %\caption{...}
  %\label{fig:y_funnel_4}
  \end{center}
\end{subfigure}
\caption{Snapshots of reentry waves for $t=100$ (left) and $t=200$ (right).}
\label{fig:reentry_waves}
\end{figure}

\begin{figure}[H]
\begin{center}
  \input{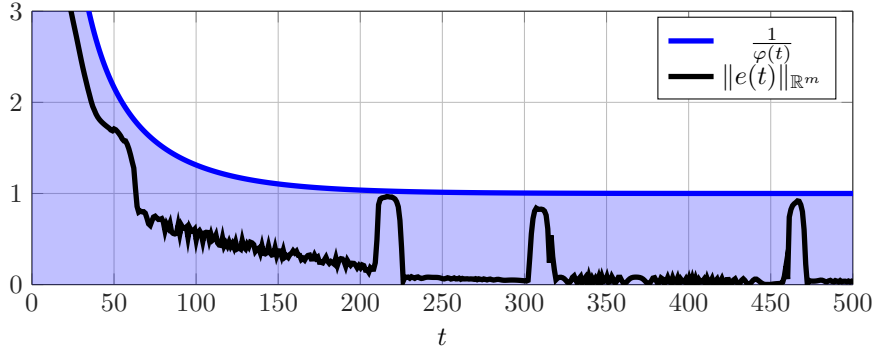}
  \end{center}
\vspace{-5mm}
    \caption{Error dynamics and funnel boundary.}
  \label{fig:funnel_error}
\end{figure}

In this section, we illustrate the practical applicability of the funnel controller by means of a numerical example. The setup chosen here is a standard test example for termination of reentry waves and has been considered similarly e.g.\ in~\cite{BreiKuni17,KuniNagaWagn11}. All simulations are generated on an AMD Ryzen 7 1800X @ 3.68 GHz x 16,
64 GB RAM,  \matlab \;Version 9.2.0.538062 (R2017a). The solutions of the ODE systems are obtained by the \matlab\;routine \texttt{ode23}. The parameters for the FitzHugh-Nagumo model \eqref{eq:FHN_model} used here are as follows:
\begin{align*}
\Omega&=(0,1)^2,\ \ D=\begin{bmatrix} 0.015 & 0 \\ 0 & 0.015 \end{bmatrix}, \ \
\begin{pmatrix}
c_1 \\ c_2 \\ c_3 \\ c_4 \\ c_5
\end{pmatrix}
\approx
\begin{pmatrix}
1.614\\ 0.1403 \\ 0.012\\ 0.00015\\ 0.015
\end{pmatrix}.
\end{align*}
 The spatially discrete system of ODEs corresponds to a finite element discretization with piecewise linear finite elements on a uniform $64\times 64$ mesh.
 %The uncontrolled system is stimulated such that reentry phenomena occur, see Fig.~\ref{fig:reentry_waves}. Let us emphasize that the associated dynamics are obtained from \eqref{eq:FHN_model} with $I_{s,i}=0=I_{s,e}$.
  For the control action, we assume that $\mathcal{B}\in \mathcal{L}(\mathbb R^4,W^{1,2}(\Omega)')$, where the {Robin} control operator is defined by
 \begin{align*}
 \mathcal{B}'z &= \begin{pmatrix} \int_{\Gamma_1} z(\xi)\, \mathrm{d}\sigma,\int_{\Gamma_2} z(\xi)\, \mathrm{d}\sigma,\int_{\Gamma_3} z(\xi)\, \mathrm{d}\sigma ,\int_{\Gamma_4} z(\xi) \,\mathrm{d}\sigma\end{pmatrix}^\top,
     \\ \Gamma_1 &= \{1\}\times [0,1],  \ \ \Gamma_2= [0,1]\times \{1\}, \ \ \Gamma_3 = \{0\}\times [0,1], \ \ \Gamma=[0,1]\times \{0\}.
 \end{align*}
 The purpose of the numerical example is to model a typical defibrillation process as a tracking problem as discussed above. In this context, system \eqref{eq:FHN_model} is initialized with $(v(0),u(0))=(v_0^*,u_0^*)$ and $I_{s,i}=0=I_{s,e}$, where $(v_0^*,u_0^*)$ is an arbitrary snapshot of a reentry wave. The resulting reentry phenomena are shown in Fig.~\ref{fig:reentry_waves} and resemble a dysfunctional heart rhythm which impedes  the intracellular stimulation current~$I_{s,i}$. The objective is to design a stimulation current $I_{s,e}$ such that the dynamics return to a natural heart rhythm modeled by a reference trajectory $y_{\text{ref}}$. The trajectory $y_{\text{ref}} = \cB' v_{\text{ref}}$ corresponds to a solution $(v_{\text{ref}},u_{\text{ref}})$ of~\eqref{eq:FHN_model} with $(v_{\text{ref}}(0),u_{\text{ref}}(0))=(0,0)$, $I_{s,e}=0$ and
 \begin{align*}
  I_{s,i}(t) = 101\cdot w(\xi) (\chi_{[49,51]}(t) + \chi_{[299,301]}(t)),
  \end{align*}
   where the excitation domain of the intracellular stimulation current~$I_{s,i}$ is described by
 \begin{align*}
  w(\xi) = \begin{cases} 1 , & \text{if } (\xi_1-\frac{1}{2})^2+(\xi_2-\frac{1}{2})^2 \le 0.0225, \\
  0 , & \text{otherwise}. \end{cases}
 \end{align*}
  The smoothness of the signal is guaranteed by convoluting the original signal with a triangular function.
The function $\varphi$ characterizing the performance funnel (see Fig.~\ref{fig:funnel_error}) is chosen as
\begin{align*}
\varphi(t)=
\begin{cases}
0, & t \in [0,0.05] ,\\
\mathrm{tanh}(\frac{t}{100}), & t > 0.05 .
\end{cases}
\end{align*}

%The simulations correspond to a finite element discretization of \eqref{eq:FHN_model} with piecewise linear finite elements

\begin{figure}[H]
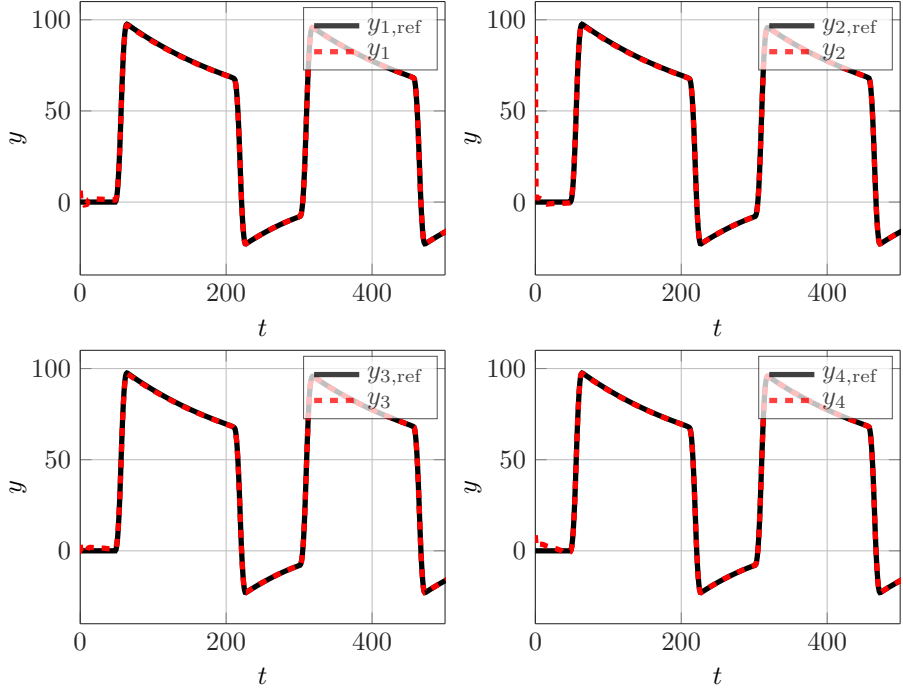

\begin{subfigure}{.47\linewidth}
\begin{center}
  \input{y_funnel_1.tikz}
  %\caption{...}
  %\label{fig:y_funnel_1}
  \end{center}
\end{subfigure}\quad
\begin{subfigure}{.47\linewidth}
\begin{center}
  \input{y_funnel_2.tikz}
  %\caption{...}
  %\label{fig:y_funnel_2}
  \end{center}
\end{subfigure}
\begin{subfigure}{.47\linewidth}
\begin{center}
  \input{y_funnel_3.tikz}
  %\caption{...}
  %\label{fig:y_funnel_3}
  \end{center}
\end{subfigure}\quad
\begin{subfigure}{.47\linewidth}
\begin{center}
  \input{y_funnel_4.tikz}
  %\caption{...}
  %\label{fig:y_funnel_4}
  \end{center}
\end{subfigure}
\caption{Reference signals and outputs of the funnel controlled system.}
\label{fig:y_funnel}
\end{figure}

\begin{figure}[H]
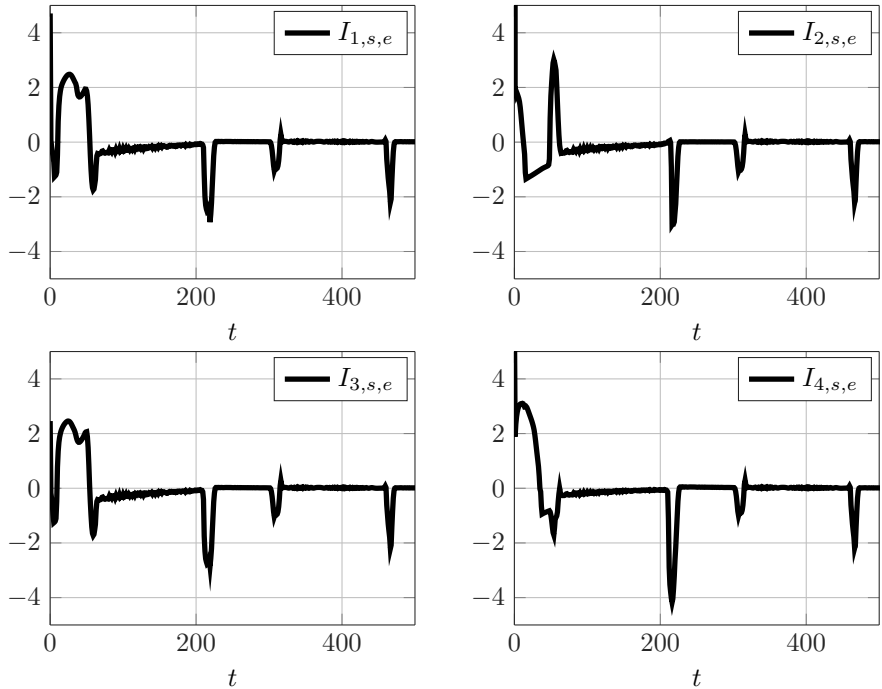

\begin{subfigure}{.48\linewidth}
\begin{center}
  \input{u_funnel_1.tikz}
  %\caption{...}
  %\label{fig:y_funnel_1}
  \end{center}
\end{subfigure}\quad
\begin{subfigure}{.48\linewidth}
\begin{center}
  \input{u_funnel_2.tikz}
  %\caption{...}
  %\label{fig:y_funnel_2}
  \end{center}
\end{subfigure}
\begin{subfigure}{.48\linewidth}
\begin{center}
  \input{u_funnel_3.tikz}
  %\caption{...}
  %\label{fig:y_funnel_3}
  \end{center}
\end{subfigure}\quad
\begin{subfigure}{.48\linewidth}
\begin{center}
  \input{u_funnel_4.tikz}
  %\caption{...}
  %\label{fig:y_funnel_4}
  \end{center}
\end{subfigure}
\caption{Funnel control laws.}
\label{fig:u_funnel}
\end{figure}

 Fig.~\ref{fig:y_funnel} shows the results of the closed-loop system for $(v(0),u(0))=(v_0^*,u_0^*)$ and the control law
\begin{align*}
 I_{s,e}(t)=-\frac{0.75}{1-\varphi(t)^2\|\cB' v(t)-\yref(t)\|^2_{\R^m}}(\cB'v(t)-\yref(t)),
\end{align*}
which is visualized in Fig.~\ref{fig:u_funnel}. Let us note that the sudden changes in the feedback law are due to the jump discontinuities of the intracellular stimulation current $I_{s,i}$ used for simulating a regular heart beat.

\setlength\fheight{0.3\linewidth}
\setlength\fwidth{0.4\linewidth}

We see from Fig.~\ref{fig:y_funnel} that the controlled system tracks the desired reference signal with the prescribed performance. Also note that the performance constraints are not active on the interval $[0,0.05]$. Fig.~\ref{fig:u_funnel} further shows that the tracking is achieved with a comparably small control effort.

\section{Conclusions}
\label{sec:conclusions}

{In this work we have proved existence and uniqueness of global bounded solutions of a reaction diffusion system under funnel control. The considered monodomain equations with the FitzHugh-Nagumo model~\eqref{eq:FHN_model} are a relevant system arising in mathematical biology. The considered input-output configurations allow for both distributed and boundary control and observation. The proposed funnel control feedback law~\eqref{eq:monodomain_funnel_controller} renders the closed-loop system~\eqref{eq:FHN_feedback} a nonlinear and non-autonomous PDE with the requirement that the tracking error associated to any solution evolves in the performance funnel $\cF_\varphi$. In the main result Theorem~\ref{thm:mono_funnel} we have put special emphasis on the regularity properties of the solutions of~\eqref{eq:FHN_feedback}.}

{The present work is the basis for extensions in several directions. Apart from more general reaction diffusion systems with more complex nonlinearities and other nonlinear parabolic equations from mathematical biology (such as the Keller-Segel system modelling chemotactic behavior~\cite{KellSege70}) an important topic for future research is the investigation of the bidomain model of the human heart~\cite{SundLine07}. Since this model is closer to reality, the authors expect that successfully applying funnel control methods will potentially lead to real-world applications such as in implantable cardioverter defibrillators.}

%\appendix
\renewcommand*\appendixpagename{\Large Appendix}
\begin{appendices}

\section{Interpolation spaces}
\label{sec:mono_prep_proof}

We collect some results on interpolation spaces, which are necessary for the proof of Theorem~\ref{thm:mono_funnel}. For a (more) general interpolation theory, we refer to \cite{Luna18}.

\begin{Def}\label{def:interp}
Let $X,Y$ be Hilbert spaces and let $\alpha\in[0,1]$. Consider the function
\[K:(0,\infty)\times (X+Y)\to\R,\ (t,x)\mapsto \underset{x=a+b}{\inf_{a\in X,\, b\in Y,}}\, \|a\|_X+t\|b\|_Y.\]
The {\em interpolation space} $(X,Y)_{\alpha}$ is defined by
\[(X,Y)_{\alpha}:=\setdef{x\in X+Y}{\Big(t\mapsto t^{-\alpha} K(t,x)\Big)\in L^2(0,\infty)},\]
and it is a Hilbert space with the norm
\[\|x\|_{(X,Y)_\alpha}=\|t\mapsto t^{-\alpha} K(t,x)\|_{L^2}.\]
\end{Def}

Note that interpolation can be performed in a more general fashion for Banach spaces $X$, $Y$. More precise, we may utilize the $L^p$-norm of the map $t\mapsto t^{-\alpha} K(t,x)$ for some $p\in[1,\infty)$ instead of the $L^2$-norm in the above definition. However, this does not lead to Hilbert spaces $(X,Y)_\alpha$, not even when~$X$ and~$Y$ are Hilbert spaces.

For a~self-adjoint operator $A:\cD(A)\subset X\to X$, $X$ a Hilbert space and $n\in\N$, we may define the space $X_n:=\cD(A^n)$ by $X_0=X$ and $X_{n+1}:=\setdef{x\in X_n}{Ax\in X_n}$. This is a Hilbert space with norm $\|z\|_{X_{n+1}}=\|-\lambda z+Az\|_{X_n}$, where $\lambda\in\C$ is in the resolvent set of $A$.
Likewise, we introduce $X_{-n}$ as the completion of $X$ with respect to the norm $\|z\|_{X_{-n}}=\|(-\lambda I+A)^{-n}z\|$. Note that $X_{-n}$ is the dual of $X_n$ with respect to the pivot space $X$, cf.~\cite[Sec.~2.10]{TucsWeis09}.
Using interpolation theory, we may further introduce the spaces~$X_\alpha$ for any $\alpha\in\R$ as follows.

\begin{Def}\label{Def:int-space-A}
Let $\alpha\in\R$, $X$ a Hilbert space and $A:\cD(A)\subset X\to X$ be self-adjoint. Further, let $n\in\Z$ be such that $\alpha\in[n,n+1)$.
The space $X_\alpha$ is defined as the interpolation space
\[X_\alpha=(X_{n},X_{n+1})_{\alpha-n}.\]
\end{Def}

The reiteration theorem, see~\cite[Cor.~1.24]{Luna18}, together with~\cite[Prop.~3.8]{Luna18} yields that for all $\alpha\in[0,1]$ and $\alpha_1,\alpha_2\in\R$ with $\alpha_1\leq\alpha_2$ we have that
\begin{equation}\label{eq:Xalpha-reit}
    (X_{\alpha_1},X_{\alpha_2})_{\alpha}=X_{\alpha_1+\alpha (\alpha_2-\alpha_1)}.
\end{equation}
%If $B$ has a spectral decomposition, i.e., for an orthonormal system $(v_j)_{j\in\N_0}$ and a~real sequence $(\beta_j)$ with no accumulation points such that
%\[Bx=\sum_{k=1}\beta_j \scpr{x}{v_j} v_j\]
Next we characterize interpolation spaces associated with the {Robin} elliptic operator $\cA$ {as in~\eqref{eq:RobinOp}}.

{\begin{Rem}
\label{rem:X_alpha}
Let Assumption~\ref{Ass1} hold and $\cA$ be the Robin elliptic operator as in~\eqref{eq:RobinOp}. Further let $X_\alpha$, $\alpha \in\R$, be the corresponding interpolation spaces with, in particular, $X=X_0=L^2(\Omega)$. Then
the equation $X_{1/2}=W^{1,2}(\Omega)$ is an immediate consequence of a~combination of \eqref{eq:faprop} with Kato's second representation theorem~\cite[Sec.~VI.2, Thm.~2.23]{Kato80}.
Further, \eqref{eq:Xalpha-reit} implies that
\[
    X_{r/2}=(X_{0},X_{1/2})_{r}\;\text{  for all $r\in[0,1]$}.
\]
On the other hand,
\cite[Thm.~1.35]{Yagi10} gives
$(L^2(\Omega),W^{1,2}(\Omega))_{r}=W^{r,2}(\Omega)$,
whence we obtain
\begin{equation}
    X_{r/2}=W^{r,2}(\Omega)\;\text{  for all $r\in[0,1]$}.
\label{eq:Xr12}
\end{equation}
In terms of the spectral decomposition \eqref{eq:spectr}, the interpolation space has the representation
\begin{equation}X_\alpha=\setdef{\sum_{j=0}^\infty \lambda_j \theta_j}{(\lambda_j)_{j\in\N_0}\text{ with }\sum_{j=0}^\infty \alpha_j^{2\alpha} |\lambda_j|^2<\infty},\label{eq:Xrspec}\end{equation}
which follows from a~combination of~\cite[Thm.~4.33]{Luna18} with~\cite[Thm.~4.36]{Luna18}.
\end{Rem}}

\section{Abstract Cauchy problems and regularity}
\label{sec:mono_prep_proof2}

We consider mild solutions of certain abstract Cauchy problems and the concept of admissible control operators. This notion is well-known in infinite-dimensional linear systems theory with unbounded control and observation operators and we refer to~\cite{TucsWeis09} for further details.

Let $X$ be a real Hilbert space and recall that a semigroup $\sg$ on $X$ is a $\mathcal{L}(X,X)$-valued map satisfying $\T_0=I_{X}$ and $\T_{t+s}=\T_t \T_s$, $s,t\geq0$, where $I_{X}$ denotes the identity operator, and $t\mapsto \T_t x$ is continuous for every $x\in X$. Semigroups are characterized by their generator~$A$, which is a, not necessarily bounded, operator on~$X$. If $A:\cD(A)\subset X\to X$ is self-adjoint with $\scpr{x}{Ax}\leq0$ for all $x\in\cD(A)$, then
it generates a~contractive, analytic semigroup $\sg$ on $X$, cf.~\cite[Thm.~4.2]{ArenElst12}. Furthermore, if additionally there exists $\omega_0>0$ such that $\scpr{x}{Ax}\leq-\omega_0 \|x\|^2$ for all $x\in\cD(A)$, then the semigroup $\sg$  generated by $A$ satisfies $\|\T_t\|\leq \ee^{-\omega_0 t}$ for all $t\geq0$; the smallest number $\omega_0$ for which this is true is called \emph{growth bound} of $\sg$. We can further conclude from~\cite[Thm.~6.13\,(b)]{Pazy83} that, for all $\alpha\in\R$, $\sg$ restricts (resp.\ extends) to an analytic semigroup $((\T|_{\alpha})_t)_{t\ge0}$ on $X_\alpha$ with same growth bound as $\sg$. Furthermore, we have $\im \T_t\subset X_r$ for all $t>0$ and $r\in\R$, see~\cite[Thm.~6.13(a)]{Pazy83}. In the following we present an estimate for the corresponding operator norm.

\begin{Lem}\label{lem:A_alpha}
	Assume that $A:\cD(A)\subset X\to X$, $X$ a Hilbert space, is self-adjoint and there exists $\omega_0>0$ with $\scpr{x}{Ax}\leq-\omega_0 \|x\|^2$ for all $x\in\cD(A)$. Then there exist $M,\omega>0$ such that  the semigroup $\sg$ generated by $A$ satisfies
	\[
    \forall\, \alpha\in[0,2]\ \forall\, t>0:\ \|\T_t\|_{\cL(X,X_\alpha)}\leq M(1+t^{-\alpha})\mathrm{e}^{-\omega t}.\]
	Thus, for each $\alpha\in[0,2]$ there exists $K>0$ such that
	\[\sup_{t\in[0,\infty)}t^\alpha\|\T_t\|_{\cL(X,X_\alpha)}<K.\]
\end{Lem}

\begin{proof}
	Since $A$ with the above properties generates an exponentially stable analytic semigroup $\sg$, the cases $\alpha\in[0,1]$ and $\alpha=2$ follow from~\cite[Cor.~3.10.8~\&~Lem.~3.10.9]{Staf05}. The result for $\alpha\in[1,2]$ is a consequence of~\cite[Lem~3.9.8]{Staf05} and interpolation between $X_1$ and $X_2$, cf.\ Appendix~\ref{sec:mono_prep_proof}.
\end{proof}

Next we consider the abstract Cauchy problem with source term.

\begin{Def}\label{Def:Cauchy}
	Let $X$ be a Hilbert space, $A:\cD(A)\subset X\to X$ be self-adjoint with $\scpr{x}{Ax}\leq0$ for all $x\in\cD(A)$, $T\in(0,\infty]$, and $\alpha\in[0,1]$. Let $\sg$ be the semigroup on $X$ generated by $A$, and let $B\in\cL(\R^m,X_{-\alpha})$. For $x_0\in X$, $p\in[1,\infty]$, $f\in L^p_{\loc}(0,T;X)$ and $u\in L^p_{\loc}(0,T;\R^m)$, we call~$x:[0,T)\to X$ a \emph{mild solution} of
	\begin{equation}\label{eq:abstract_cauchy}
	\begin{aligned}
	\dot{x}(t)=Ax(t)+f(t)+Bu(t),\quad
	x(0)=x_0
	\end{aligned}
	\end{equation}
	on $[0,T)$, if it satisfies
	\begin{equation}\label{eq:mild_solution}
	\forall\, t\in[0,T):\ x(t)=\T_tx_0+\int_0^t\T_{t-s}f(s)\ds{s}+\int_0^t(\T|_{-\alpha})_{t-s}Bu(s)\ds{s}.
	\end{equation}
%where $((\T|_{-\alpha})_t)_{t\ge0}$ is the extension of $\sg$ to $X_{-\alpha}$.\\
We further call $x:[0,T)\to X$ a \emph{strong solution} of \eqref{eq:abstract_cauchy} on $[0,T)$, if $x$ in~\eqref{eq:mild_solution}  satisfies $x\in C([0,T);X)\cap W^{1,p}_{\rm loc}(0,T;X_{-1})$.
\end{Def}

Definition~\ref{Def:Cauchy} requires that the integral $\int_0^t(\T|_{-\alpha})_{t-s}Bu(s)\ds{s}$ is in $X$, whilst the integrand is not necessarily in $X$. This motivates the definition of admissibility, which is now introduced for self-adjoint $A$. Note that admissibility can  also be defined for arbitrary generators of semigroups, see~\cite{TucsWeis09}.

\begin{Def}\label{Def:Adm}
Let $X$ be a Hilbert space, $A:\cD(A)\subset X\to X$ be self-adjoint with $\scpr{x}{Ax}\leq0$ for all $x\in\cD(A)$, $T\in(0,\infty]$, $\alpha\in[0,1]$ and $p\in[1,\infty]$. Let $\sg$ be the semigroup on $X$ generated by $A$, and let $B\in\cL(\R^m,X_{-\alpha})$. Then $B$ is called an {\em $L^p$-admissible (control operator) for $\sg$}, if for some (and hence any) $t> 0$ we have
$$\forall\,  u\in L^p(0,t;\R^m):\ \Phi_{t}u :=\int_0^t(\T|_{-\alpha})_{t-s}Bu(s)\ds{s} \in X.$$
By a closed graph theorem argument this implies that $\Phi_t\in \cL(L^p(0,t;\R^m),X)$ for all $t> 0$. We call $B$ an  {\em infinite-time $L^p$-admissible (control operator) for $\sg$}, if
\[
    \sup_{t>0} \|\Phi_t\| < \infty.
\]
\end{Def}

In the following we show that for $p\ge 2$ and $\alpha\le 1/2$ any~$B$ is admissible and the mild solution of the abstract Cauchy problem is indeed a strong solution.

\begin{Lem}\label{lem:abstract_solution}
	Let $X$ be a Hilbert space, $A:\cD(A)\subset X\to X$ be self-adjoint with $\scpr{x}{Ax}\leq0$ for all $x\in\cD(A)$, $B\in\cL(\R^m,X_{-\alpha})$ for some $\alpha\in[0,1/2]$, and $\sg$ be the analytic  semigroup generated by $A$. Then for all $p\in[2,\infty]$ we have that~$B$ is $L^p$-admissible for $\sg$.
	
	Furthermore, for all $x_0\in X$, $T\in(0,\infty]$, $f\in L^p_{\loc}(0,T;X)$ and $u\in L^p_{\loc}(0,T;\R^m)$, the function~$x$ in~\eqref{eq:mild_solution} is a strong solution of~\eqref{eq:abstract_cauchy} on $[0,T)$.
\end{Lem}
\begin{proof}
For the case $p=2$, there exists a unique strong solution in $X_{-1}$ (that is, we replace $X$ by $X_{-1}$ and $X_{-1}$ by $X_{-2}$ in the definition) given by~\eqref{eq:mild_solution} and at most one strong solution in $X$, see for instance~\cite[Thm.~3.8.2~(i)~\&~(ii)]{Staf05}, so we only need to check that all the elements are in the correct spaces. Since $A$ is self-adjoint, the semigroup generated by $A$ is self-adjoint as well. Further, by combining~\cite[Prop.~5.1.3]{TucsWeis09} with~\cite[Thm.~4.4.3]{TucsWeis09}, we find that~$B$ is an $L^2$-admissible control operator for~$\sg$. Moreover, by~\cite[Prop.~4.2.5]{TucsWeis09} we have that
	\[\left(t\mapsto\T_t x_0+\int_0^t(\T|_{-\alpha})_{t-s}Bu(s)\ds{s}\right)\in C([0,T);X)\cap W^{1,2}_{\rm loc}(0,T;X_{-1})\]
and from \cite[Thm.~3.8.2~(iv)]{Staf05},
	\[\left(t\mapsto\int_0^t\T_{t-s}f(s)\ds{s}\right)\in C([0,T);X)\cap W^{1,2}_{\rm loc}(0,T;X_{-1}),\]
	whence $x\in C([0,T);X)\cap W^{1,2}_{\rm loc}(0,T;X_{-1})$, which proves that~$x$ is a strong solution of~\eqref{eq:abstract_cauchy} on $[0,T)$.

Since $B$ is $L^2$-admissible, it follows from the nesting property of $L^p$ on finite intervals that~$B$ is an $L^p$-admissible control operator for~$\sg$ for all $p\in[2,\infty]$. Furthermore, for $p>2$, set $\tilde{f}\coloneqq f+Bu$ and apply~\cite[Thm.~3.10.10]{Staf05} with $\tilde{f}\in L^\infty_{\rm loc}(0,T;X_{-\alpha})$ to conclude that~$x$ is a strong solution.
\end{proof}

%Observe that if $x$ is the unique solution of \eqref{eq:abstract_cauchy} under the conditions of Lemma \ref{lem:abstract_solution} and $A$ is the operator associated to a form $a:X_{1/2}\times X_{1/2}\to \R$, then this is also the unique solution of the weak formulation
%\[
%\ddt\scpr{x(t)}{\theta}=-a(x(t),\theta)+\scpr{f(t)}{\theta}+\scpr{u(t)}{B' \theta},\quad\forall\theta\in X_{1/2},
%\]
%which can be easily seen by adapting~\cite[Thm.~1.33]{HinzPinn09} and~\cite[Thm.~3.1.7]{CurtZwar95} to our setting.

Next we show the regularity properties of the solution of~\eqref{eq:abstract_cauchy}, if $A = \cA$ and $B = \cB$ are as in the model~\eqref{eq:FHN_model}. Note that this result also holds when considering some $t_0\ge 0$, $T\in(t_0,\infty]$, and the initial condition $x(t_0)=x_0$ (instead of $x(0)=x^0$) by some straightforward modifications, cf.~\cite[Sec.~3.8]{Staf05}.

\begin{Prop}\label{prop:hoelder}
Let Assumption~\ref{Ass1} hold, $\cA$ be the {Robin} elliptic operator {as in~\eqref{eq:RobinOp}}, $T\in(0,\infty]$, and $c>0$.
% Further let $X=X_0=L^2(\Omega)$ and $X_r$, $r\in\R$, be the interpolation spaces corresponding to~$\cA$ according to Definition~\ref{Def:int-space-A}.
	Define $\cA_0\coloneqq  \cA-cI$ with $\cD(\cA_0)=\cD(\cA)$ and consider $\cB\in\cL(\R^m,W^{r,2}(\Omega)')$ for $r\in[0,1]$, $u\in L_{\rm loc}^2(0,T;\R^m)\cap L^\infty(\delta,T;\R^m)$ and $f\in L_{\rm loc}^2(0,T;L^2(\Omega))\cap L^\infty(\delta,T;L^2(\Omega))$ for all $\delta>0$. Then for all $x_0\in L^2(\Omega)$ and all $\delta>0$, the mild solution of~\eqref{eq:abstract_cauchy} (with $A=\cA_0$ and $B=\cB$) on~$[0,T)$, given by~$x$ as in~\eqref{eq:mild_solution}, satisfies
{\begin{enumerate}[(i)]
		\item if $r=0$, then
		\[\forall\, \lambda\in(0,1):\ x\in BC([0,T);L^2(\Omega))\cap C^{0,\lambda}([\delta,T);L^2(\Omega));\]
		\item if $r\in(0,1)$, then
		\[x\in BC([0,T);L^2(\Omega))\cap C^{0,1-r/2}([\delta,T);L^2(\Omega))\cap C^{0,1-r}([\delta,T);W^{r/2,2}(\Omega));\]
		\item if $r=1$, then
		\[x\in BC([0,T);L^2(\Omega))\cap C^{0,1/2}([\delta,T);L^2(\Omega))\cap BUC([\delta,T);W^{1,2}(\Omega)).\]
	\end{enumerate}}
\end{Prop}
\begin{proof}
{For brevity we set $X=X_0=L^2(\Omega)$, and let $X_\alpha$, $\alpha\in\R$, be the interpolation spaces corresponding to~$\cA$ according to Definition~\ref{Def:int-space-A}. Observe that the Robin elliptic operator satisfies the assumptions of Lemma~\ref{lem:abstract_solution}} with $p=2$, hence~$x$ as in~\eqref{eq:mild_solution} is a strong solution of~\eqref{eq:abstract_cauchy} on $[0,T)$ in the sense of Definition~\ref{Def:Cauchy}. In the following we restrict ourselves to the case $T=\infty$, and the assertions for $T<\infty$ follow from these arguments by considering the restrictions to $[0,T)$.
	Define, for $t\ge 0$, the functions
	\begin{align}
	x_h(t)\coloneqq  \T_tx_0,\quad	x_f(t)\coloneqq  \int_0^t\T_{t-s}f(s)\ds{s},\quad
	x_u(t)\coloneqq  \int_0^t(\T|_{-\alpha})_{t-s}\cB u(s)\ds{s},\label{eq:xhxfxu}
	\end{align}
	so that $x=x_h+x_f+x_u$.

\emph{Step 1}: We show that $x\in BC([0,\infty);X)$. We obtain from Remark~\ref{Rem:Aop_n}~\ref{item:Aop4}) that for all $z\in \cD(\cA)$ we have $\langle z,\cA_0 z\rangle\leq -c\|z\|^2$. The self-adjointness of $\cA$ moreover implies that $\cA_0$ is self-adjoint, whence~\cite[Thm.~4.2]{ArenElst12} gives that $\cA_0$ generates an analytic, contractive semigroup $\sg$ on $X$, which satisfies
\begin{equation}\label{eq:est-sg-exp}
    \forall\, t\ge 0\ \forall\, x\in X:\ \|\T_t x\|\leq \ee^{-ct}\|x\|.
\end{equation}
Since, by Lemma~\ref{lem:abstract_solution}, $x$ is a strong solution, we have $x\in C([0,\infty);X)\cap W^{1,2}_{\rm loc}(0,\infty;X_{-1})$. Further observe that $\cB$ is $L^\infty$-admissible by Lemma~\ref{lem:abstract_solution}. Then it follows from~\eqref{eq:est-sg-exp} and~\cite[Lem.~2.9\,(i)]{JacoNabi18} that~$\cB$ is infinite-time $L^\infty$-admissible, which implies that for $x_u$ as in~\eqref{eq:xhxfxu} we have
\[
    \|x_u\|_\infty \le \left( \sup_{t>0} \|\Phi_t\|\right) \|u\|_\infty < \infty,
\]
thus  $x_u\in BC([0,\infty);X)$. A direct calculation using~\eqref{eq:est-sg-exp} further shows that $x_h,x_f\in BC([0,\infty);X)$, whence $x\in BC([0,\infty);X)$.

\emph{Step 2}: We show~(i). Let $\delta>0$ and set $\tilde{f}:=f+Bu\in L^2_{\rm loc}(0,\infty;X)\cap L^\infty(\delta,\infty;X)$, then we may infer  from~\cite[Props.~4.2.3~\&~4.4.1\,(i)]{Luna95} that
	\[\forall\,\lambda\in(0,1):\ x\in C^{0,\lambda}([\delta,\infty);X).\]
From this together with Step~1 we may infer~(i).

\emph{Step 3}: We show~(ii). Let $\delta>0$, then it follows from~\cite[Props.~4.2.3~\&~4.4.1\,(i)]{Luna95} together with $x_0\in X$ and $f\in L^\infty(\delta,\infty;X)$, that
	\[\begin{aligned}
x_h+x_f&\in C^{0,1-r/2}([\delta,\infty);X_{r/2})\cap C^{1}([\delta,\infty);X)\\
&= C^{0,1-r}([\delta,\infty);X_{r/2})\cap C^{0,1-r/2}([\delta,\infty);X).\end{aligned}\]
Since we have shown in Step~1 that $x\in BC([0,\infty),X)$, it remains to show that
	$x_u\in C^{0,1-r}([\delta,\infty);X_{r/2})\cap C^{0,1-r/2}([\delta,\infty);X)$.\\
To this end, consider the space $Y:=X_{-r/2}$. Then $\sg$ extends to a~semigroup~$\big((\T|_{-r/2})_{t}\big)_{t\ge 0}$ on $Y$ with generator $\cA_{0,r/2}:\cD(\cA_{0,r/2})=X_{-r/2+1}\subset X_{-r/2}=Y$, cf.~\cite[pp.~50]{Luna95}. Now, for $\alpha\in\R$, consider the interpolation spaces $Y_\alpha$ as in Definition~\ref{Def:int-space-A} by means of the operator $\cA_{0,r/2}$. Then it is straightforward to show that $Y_n = \cD(\cA_{0,r/2}^n) = X_{n-r/2}$ for all $n\in\N$ using the representation~\eqref{eq:Xrspec}. Similarly, we may show that $Y_n = X_{n-r/2}$ for all $n\in\Z$. Then the reiteration theorem, see~\cite[Cor.~1.24]{Luna18} and also \eqref{eq:Xalpha-reit}, gives
\[\forall\, \alpha\in\R\,:\quad Y_\alpha=X_{\alpha-r/2}.\]
Since $\cB\in\cL(\R^m,Y)$, \cite[Props.~4.2.3~\&~4.4.1\,(i)]{Luna95} now imply
\[\begin{aligned}
x_u &\in C^{0,1-r}([\delta,\infty);Y_{r})\cap C^{0,1-r/2}([\delta,\infty);Y_{r/2})\\
&= C^{0,1-r}([\delta,\infty);X_{r/2})\cap C^{0,1-r/2}([\delta,\infty);X),
\end{aligned}\]
{which leads to $x\in C^{0,1-r}([\delta,\infty);X_{r/2})\cap C^{0,1-r/2}([\delta,\infty);X)$, and by further using~\eqref{eq:Xr12}, we may conclude statement~(ii).}

\emph{Step 4}: We show~(iii). The proof of $x\in C^{0,1/2}([\delta,\infty);X)$ is analogous to that of $x\in C^{0,1-r/2}([\delta,\infty);X)$ in Step~3. Boundedness and continuity of $x$ on $[0,\infty)$ was proved in Step~1. Since we have
$X_{1/2}=W^{1,2}(\Omega)$ by \eqref{eq:Xr12}, it remains to show that $x$ is uniformly continuous as a~mapping to $X_{1/2}$: Again consider the additive decomposition of $x$ into $x_h$, $x_f$ and $x_u$ as in~\eqref{eq:xhxfxu}. Similar to Step~3 it can be shown that $x_h,x_f\in C^{0,1/2}([\delta,\infty);X_{1/2})$, whence $x_h,x_f\in BUC([\delta,\infty);X_{1/2})$. It remains to show that $x_u\in BUC([\delta,\infty);X_{1/2})$.

Note that Lemma \ref{lem:abstract_solution} gives that $x_\delta\coloneqq x(\delta)\in X$. Then $x_u$ solves $\dot z(t) = \cA_0 z(t) + \cB u(t)$ with $z(\delta) = x_u(\delta)$ and hence, for all $t\geq\delta$ we have
\begin{equation}\label{eq:x-T-delta}
\begin{aligned}
    x_u(t)=&\,\T_{t-\delta}x_u(\delta)+\underbrace{\int_\delta^t(\T|_{-\alpha})_{t-s}\cB u(s)\ds{s}}_{=:x_u^\delta(t)}\\
%    =&\,\T_{t-\delta}x_u(\delta)+\int_0^{t-\delta}(\T|_{-\alpha})_{t-(s+\delta)}\cB u(s+\delta)\ds{s}
    \end{aligned}
\end{equation}
Since $x_u(\delta)\in X$ by Lemma~\ref{lem:abstract_solution}, it remains to show that $x_u^\delta\in BUC([\delta,\infty);X_{1/2})$.
We obtain from Remark~\ref{Rem:Aop_n}\,\ref{item:Aop4}) that $\cA_0$ has an eigendecomposition of type~\eqref{eq:spectr}
	with eigenvalues $(-\beta_j)_{j\in\N_0}$, $\beta_j\coloneqq  \alpha_j+c$, and eigenfunctions $(\theta_j)_{j\in\N_0}$. Moreover, there exist $b_i\in X_{-1/2}$ for $i=1,\dots,m$ such that $\cB \xi = \sum_{i=1}^m b_i \cdot \xi_i$ for all $\xi\in\R^m$. Therefore,
	\begin{align*}
	x^\delta_u(t)&=\int_\delta^t\sum_{j=0}^\infty \ee^{-\beta_j(t-\tau)}\theta_j \sum_{i=1}^m\scpr{b_i \cdot u_i(\tau)}{\theta_j}\ds{\tau}\\
	&=\int_\delta^t\sum_{j=0}^\infty \ee^{-\beta_j(t-\tau)}\theta_j \sum_{i=1}^m u_i(\tau) \scpr{b_i}{\theta_j}\ds{\tau},
	\end{align*}
	where the last equality holds since $u_i(\tau)\in\R$ and can be treated as a constant in~$X$. By considering each of the factors in the sum over $i=1,\dots,m$, we can assume without loss of generality that $m=1$ and $b\coloneqq b_1$, so that
	\[x_u^\delta(t)=\int_\delta^t\sum_{j=0}^\infty \ee^{-\beta_j(t-\tau)} u(\tau) \scpr{b}{\theta_j} \theta_j \ds{\tau}.\]
	Define $b^j \coloneqq  \scpr{b}{\theta_j}$ for $j\in\N_0$. Since $b\in X_{-1/2}$ we have that $\sum_{j=0}^\infty b_j^2/\beta_j$ converges, which implies
	\begin{equation}S\coloneqq  \sum_{j=0}^\infty\frac{(b^j)^2}{\beta_j}<\infty.\label{eq:Sdef}\end{equation}
	Recall that the spaces $X_\alpha$, $\alpha\in\R$, are defined by using $\lambda\in\C$ belonging to the resolvent set of $\cA$, and they are independent of the choice of $\lambda$. Since $c>0$ in the statement of the proposition is in the resolvent set of $\cA$, the spaces $X_\alpha$ coincide for $\cA$ and $\cA_0=\cA-cI$.\\
	Using the diagonal representation from Remark~\ref{rem:X_alpha} and \cite[Prop.~3.4.8]{TucsWeis09}, we may infer that $x_u^\delta(t)\in X_{1/2}$ for {a.a.} $t\geq\delta$, namely,
	\begin{align*}
	\|x_u^\delta(t)\|_{X_{1/2}}^2&\leq\sum_{j=0}^\infty\beta_j(b^j)^2\|u\|^2_{L^\infty(\delta,\infty)}\left(\int_\delta^t\ee^{-\beta_j(t-s)}\ds{s}\right)^2\\
	&=\|u\|^2_{L^\infty(\delta,\infty)}\sum_{j=0}^\infty\frac{(b^j)^2}{\beta_j}\left(1-\ee^{-\beta_j(t-\delta)}\right)^2\\
	&\leq\|u\|^2_{L^\infty(\delta,\infty)}\sum_{j=0}^\infty\frac{(b^j)^2}{\beta_j}<\infty.
	\end{align*}
	Hence,
	\begin{equation}\label{eq:xudelta_X12}
	\|x_u^\delta(t)\|_{X_{1/2}}\leq \|u\|_{L^\infty(\delta,\infty)}\sqrt{S}.
	\end{equation}
	Now let $t>s>\delta$ and $\sigma>0$ such that $t-s<\sigma$. By dominated convergence \cite[Thm.~II.2.3]{Dies77}, summation and integration can be interchanged, so that
	\begin{align*}
	\|x_u^\delta&(t)-x_u^\delta(s)\|_{X_{1/2}}^2\\
	\leq&\,\|u\|_{L^\infty(\delta,\infty)}^2\sum_{j=0}^\infty \beta_j (b^j)^2 \left(\int_\delta^s\ee^{-\beta_j(s-\tau)}-\ee^{-\beta_j(t-\tau)}\ds{\tau}+\int_s^t\ee^{-\beta_j(t-\tau)}\ds{\tau}\right)^2\\
	\leq&\,4\|u\|_{L^\infty(\delta,\infty)}^2\sum_{j=0}^\infty  \frac{(b^j)^2}{\beta_j} \left(1-\ee^{-\beta_j (t-s)}\right)^2\\
	\le&\, 4\|u\|_{L^\infty(\delta,\infty)}^2\sum_{j=0}^\infty \frac{(b^j)^2}{\beta_j} \left(1-\ee^{-\beta_j \sigma}\right)^2.
	\end{align*}
	We can conclude from \eqref{eq:Sdef} that the series $F:(0,\infty)\to(0,S)$ with
	\[F(\sigma)\coloneqq  \sum_{j=0}^\infty\frac{(b^j)^2}{\beta_j}(1-\ee^{-\beta_j\sigma})^2\]
	converges uniformly to a strictly monotone, continuous and surjective function. Therefore, $F$ has an inverse. The function $x_u^\delta$ is thus uniformly continuous on~$[\delta,\infty)$ and by~\eqref{eq:est-sg-exp} we obtain boundedness, i.e., $x_u^\delta\in BUC([\delta,\infty);X_{1/2})$.
\end{proof}

Finally we present a consequence of the Banach-Alaoglu theorem, see e.g.~\cite[Thm.~3.15]{Rudi91}.

\begin{Lem}\label{lem:weak_convergence}
	Let $T>0$ and $Z$ be a reflexive and separable Banach space. Then
	\begin{enumerate}[(i)]
		\item every bounded sequence $(w_n)_{n\in\N}$ in $L^\infty(0,T;Z)$ has a weak$^\star$ convergent subsequence in $L^\infty(0,T;Z)$;\label{it:weak_star}
		\item every bounded sequence $(w_n)_{n\in\N}$ in $L^p(0,T;Z)$ with $p\in(1,\infty)$ has a weakly convergent subsequence in $L^p(0,T;Z)$.\label{it:weak}
	\end{enumerate}
\end{Lem}

\begin{proof}
	Let $p\in[1,\infty)$. Then $W:=L^p(0,T;Z')$ is a separable Banach space, see~\cite[Sec.~IV.1]{Dies77}. Since $Z$ is reflexive, by \cite[Cor.~III.4]{Dies77} it has the Radon-Nikodým property. Then it follows from~\cite[Thm.~IV.1]{Dies77} that $W'=L^q(0,T;Z)$ is the dual of $W$, where $q\in(1,\infty]$ such that $p^{-1}+q^{-1}=1$. Assertion~\eqref{it:weak_star} now follows from~\cite[Thm.~3.17]{Rudi91} with $p=1$ and $q=\infty$. On the other hand, statement~\eqref{it:weak} follows from~\cite[Thm.~V.2.1]{Yosi80} by further using that~$W$ is reflexive for $p\in(1,\infty)$.
\end{proof}

\end{appendices}

\section*{Acknowledgments}
The authors would like to thank Felix L.\ Schwenninger (U Twente) and Mark R.\ Opmeer (U Bath) for helpful comments on maximal regularity.

\bibliographystyle{elsarticle-harv}
%\bibliography{MST2}

\end{document}